\documentclass[aop]{imsart}

\RequirePackage{amsthm,amsmath,amsfonts,amssymb}
\RequirePackage[numbers]{natbib}

\usepackage{hyperref}
\usepackage{csquotes}
\usepackage{graphicx}
\numberwithin{equation}{section}

\usepackage{accents}

\usepackage{tikz}

\newcommand{\1}{\mathbf{1}}

 \newcommand{\eps}{\varepsilon}

\newcommand{\dis}{\displaystyle}

\newcommand{\cq}{\mathcal{Q}}

\newcommand{\C}{\mathbb C}

\newcommand{\R}{\mathbb R}

\newcommand{\Z}{\mathbb Z}
\newcommand{\T}{\mathbb T}

\newcommand{\cac}{\mathcal C}

\newcommand{\cf}{\mathcal F}
\newcommand{\cg}{\mathcal G}

\newcommand{\ci}{\mathcal I}

\newcommand{\ck}{\mathcal K}
\newcommand{\cl}{\mathcal L}
\newcommand{\cm}{\mathcal M}

\newcommand{\cp}{\mathcal P}
\newcommand{\cs}{\mathcal S}

\newcommand{\al}{\alpha}

\newcommand{\ga}{\gamma}
\newcommand{\gga}{\Gamma}
\newcommand{\ka}{\kappa}
\newcommand{\la}{\lambda}

\newcommand{\si}{\sigma}

\newcommand{\vp}{\varphi}

\newtheorem{theorem}{Theorem}[section]

\newtheorem{lemma}[theorem]{Lemma}
\newtheorem{notation}[theorem]{Notation}

\newtheorem{proposition}[theorem]{Proposition}

\theoremstyle{definition}
\newtheorem{remark}[theorem]{Remark}

\pgfdeclareshape{crosscircle}
{
  \inheritsavedanchors[from=circle] 
  \inheritanchorborder[from=circle]
  \inheritanchor[from=circle]{north}
  \inheritanchor[from=circle]{north west}
  \inheritanchor[from=circle]{north east}
  \inheritanchor[from=circle]{center}
  \inheritanchor[from=circle]{west}
  \inheritanchor[from=circle]{east}
  \inheritanchor[from=circle]{mid}
  \inheritanchor[from=circle]{mid west}
  \inheritanchor[from=circle]{mid east}
  \inheritanchor[from=circle]{base}
  \inheritanchor[from=circle]{base west}
  \inheritanchor[from=circle]{base east}
  \inheritanchor[from=circle]{south}
  \inheritanchor[from=circle]{south west}
  \inheritanchor[from=circle]{south east}
  \inheritbackgroundpath[from=circle]
  \foregroundpath{
    \centerpoint%
    \pgf@xc=\pgf@x%
    \pgf@yc=\pgf@y%
    \pgfutil@tempdima=\radius%
    \pgfmathsetlength{\pgf@xb}{\pgfkeysvalueof{/pgf/outer xsep}}%
    \pgfmathsetlength{\pgf@yb}{\pgfkeysvalueof{/pgf/outer ysep}}%
    \ifdim\pgf@xb<\pgf@yb%
      \advance\pgfutil@tempdima by-\pgf@yb%
    \else%
      \advance\pgfutil@tempdima by-\pgf@xb%
    \fi%
    \pgfpathmoveto{\pgfpointadd{\pgfqpoint{\pgf@xc}{\pgf@yc}}{\pgfqpoint{-0.707107\pgfutil@tempdima}{-0.707107\pgfutil@tempdima}}}
    \pgfpathlineto{\pgfpointadd{\pgfqpoint{\pgf@xc}{\pgf@yc}}{\pgfqpoint{0.707107\pgfutil@tempdima}{0.707107\pgfutil@tempdima}}}
    \pgfpathmoveto{\pgfpointadd{\pgfqpoint{\pgf@xc}{\pgf@yc}}{\pgfqpoint{-0.707107\pgfutil@tempdima}{0.707107\pgfutil@tempdima}}}
    \pgfpathlineto{\pgfpointadd{\pgfqpoint{\pgf@xc}{\pgf@yc}}{\pgfqpoint{0.707107\pgfutil@tempdima}{-0.707107\pgfutil@tempdima}}}
  }
}
\makeatother

\definecolor{gr}{rgb}   {0.,   0.69,   0.23 }
\definecolor{bl}{rgb}   {0.,   0.5,   1. }
\definecolor{mg}{rgb}   {0.85,  0.,    0.85}
\definecolor{yl}{rgb}   {0.8,  0.7,   0.}
\definecolor{or}{rgb}  {0.7,0.2,0.2}

\colorlet{symbols}{black!90!black}
\colorlet{symbolsb}{black!90!black}
\colorlet{testcolor}{green!60!black}

\usetikzlibrary{shapes.misc}
\usetikzlibrary{shapes.symbols}
\usetikzlibrary{decorations}
\usetikzlibrary{decorations.markings}

\def\drawx{\draw[-,solid] (-3pt,-3pt) -- (3pt,3pt);\draw[-,solid] (-3pt,3pt) -- (3pt,-3pt);}
\tikzset{
	root/.style={circle,fill=testcolor,inner sep=0pt, minimum size=2mm},
	dot/.style={circle,fill=black,inner sep=0pt, minimum size=1mm},
	var/.style={circle,fill=black!10,draw=black,inner sep=0pt, minimum size=2mm},
	dotred/.style={circle,fill=black!50,inner sep=0pt, minimum size=2mm},
	generic/.style={semithick,shorten >=1pt,shorten <=1pt},
	dist/.style={ultra thick,draw=testcolor,shorten >=1pt,shorten <=1pt},
	testfcn/.style={ultra thick,testcolor,shorten >=1pt,shorten <=1pt,<-},
	testfcnx/.style={ultra thick,testcolor,shorten >=1pt,shorten <=1pt,<-,
		postaction={decorate,decoration={markings,mark=at position 0.6 with {\drawx}}}},
	kprime/.style={semithick,shorten >=1pt,shorten <=1pt,densely dashed,->},
	kprimex/.style={semithick,shorten >=1pt,shorten <=1pt,densely dashed,->,
		postaction={decorate,decoration={markings,mark=at position 0.4 with {\drawx}}}},
	kernel/.style={semithick,shorten >=1pt,shorten <=1pt,->},
	multx/.style={shorten >=1pt,shorten <=1pt,
		postaction={decorate,decoration={markings,mark=at position 0.5 with {\drawx}}}},
	kernelx/.style={semithick,shorten >=1pt,shorten <=1pt,->,
		postaction={decorate,decoration={markings,mark=at position 0.4 with {\drawx}}}},
	kernel1/.style={->,semithick,shorten >=1pt,shorten <=1pt,postaction={decorate,decoration={markings,mark=at position 0.45 with {\draw[-] (0,-0.1) -- (0,0.1);}}}},
	kernel2/.style={->,semithick,shorten >=1pt,shorten <=1pt,postaction={decorate,decoration={markings,mark=at position 0.45 with {\draw[-] (0.05,-0.1) -- (0.05,0.1);\draw[-] (-0.05,-0.1) -- (-0.05,0.1);}}}},
	kernelBig/.style={semithick,shorten >=1pt,shorten <=1pt,decorate, decoration={zigzag,amplitude=1.5pt,segment length = 3pt,pre length=2pt,post length=2pt}},
	rho/.style={dotted,semithick,shorten >=1pt,shorten <=1pt},
	renorm/.style={shape=circle,fill=white,inner sep=1pt},
	labl/.style={shape=rectangle,fill=white,inner sep=1pt},
	xi/.style={circle,fill=symbols!10,draw=symbols,inner sep=0pt,minimum size=1.2mm},
	xiblack/.style={circle,fill=symbolsb,draw=symbolsb,inner sep=0pt,minimum size=1.2mm},
	xix/.style={crosscircle,fill=symbols!10,draw=symbols,inner sep=0pt,minimum size=1.2mm},
	xib/.style={circle,fill=symbols!10,draw=symbols,inner sep=0pt,minimum size=1.6mm},
	xibx/.style={crosscircle,fill=symbols!10,draw=symbols,inner sep=0pt,minimum size=1.6mm},
	not/.style={circle,fill=symbols,draw=symbols,inner sep=0pt,minimum size=0.5mm},
	notblack/.style={circle,fill=symbolsb,draw=symbolsb,inner sep=0pt,minimum size=0.5mm},
	>=stealth,
	}
\makeatletter
\def\DeclareSymbol#1#2#3{\expandafter\gdef\csname MH@symb@#1\endcsname{\tikz[baseline=#2,scale=0.15,draw=symbols]{#3}}\expandafter\gdef\csname MH@symb@#1s\endcsname{\scalebox{0.7}{\tikz[baseline=#2,scale=0.15,draw=symbols]{#3}}}}
\def\<#1>{\csname MH@symb@#1\endcsname}
\makeatother

\DeclareSymbol{circle}{0.5}{\draw (0,0.7) node[xi] {};}
\DeclareSymbol{line}{0.5}{\draw (0,0.2) node[not] {} -- (0,1.4) node[not] {};}

\DeclareSymbol{Psi}{0.5}{\draw (0,0) node[not] {} -- (0,1.5) node[xi] {};}
\DeclareSymbol{Psiblack}{0.5}{\draw (0,0) node[notblack] {} -- (0,1.5) node[xiblack] {};}
\DeclareSymbol{Psi2}{0.5}{\draw (-0.6,1.5) node[xi] {} -- (0,0) node[not] {} -- (0.6,1.5) node[xi] {};}
\DeclareSymbol{IPsi2}{0}{\draw (0,1) -- (0.8,2.2) node[xi] {};\draw (0,-0.25) node[not] {} -- (0,1) node[not] {} -- (-0.8,2.2) node[xi] {};}
\DeclareSymbol{PsiIPsi2}{0}{\draw (0,1) -- (0.8,2.2) node[xi] {};\draw (0,-0.25) node[not] {} -- (0,1) node[not] {} -- (-0.8,2.2) node[xi] {};\draw (0,-0.25) node[not] {} -- (0.8,1) node[xi] {};}

\DeclareSymbol{IPsi2dot}{0}{\draw[-,dashed] (0,1) -- (0.8,2.2) node[xi] {};\draw[-,dashed] (0,-0.25) node[not] {} -- (0,1) node[not] {} -- (-0.8,2.2) node[xi] {};}

\date{\today}

\begin{document}

\begin{frontmatter}
\title{Renormalization of a 1d quadratic Schr{\"o}dinger model with additive noise}
\runtitle{Renormalization of a 1d QSM with additive noise}

\begin{aug} 
\author[A]{\fnms{Aur\'elien}~\snm{Deya}\ead[label=e1]{aurelien.deya@univ-lorraine.fr}},
\author[B]{\fnms{Reika}~\snm{Fukuizumi}\ead[label=e2]{fukuizumi@waseda.jp}}
\and
\author[A]{\fnms{Laurent}~\snm{Thomann}\ead[label=e3]{laurent.thomann@univ-lorraine.fr}}
\address[A]{Department,
Institut \'Elie Cartan, Universit\' e de Lorraine, BP 70239, 54506 Vandoeuvre-l\`es-Nancy, France\printead[presep={,\ }]{e1,e3}}

\address[B]{Department,
Department of Mathematics,
School of fundamental science and engineering,
Waseda University, 3-4-1, Okubo, Shinjuku-ku, Tokyo, Japan.\printead[presep={,\ }]{e2}}
\end{aug}

\begin{abstract}
The study is devoted to the interpretation and wellposedness of the stochastic NLS model 
\begin{equation*}
(\imath \partial_t-\Delta)u=|u|^2+\dot{B}, \quad u_0=0,\quad \quad t\in \R, \ x\in \mathbb{T},
\end{equation*}
where $\dot{B}$ stands for a space-time fractional noise with index $H=(H_0,H_1)$ in a subset of $(0,1)^{2}$. \\
\indent We first establish that in the situation where $0<2H_0+H_1\leq 2$, the equation cannot be interpreted in a (classical) functional sense.\\
\indent Our investigations then focus on the rough regime corresponding to the condition $\frac74<2H_0+H_1\leq 2$. In this specific case, we exhibit an \textit{explicit} renormalization procedure allowing to restore the (local) convergence of the approximated solutions. \\
\indent We follow a pathwise-type approach emphasizing the distinction between the stochastic objects at the core of the dynamics and the general deterministic machinery.

\end{abstract}



\end{frontmatter}
\tableofcontents

\section{Introduction}

\subsection{General considerations}\label{subsec:gene-conside}

In this paper, we propose to explore a few fundamental issues related to the 1-d quadratic Schr{\"o}dinger model with additive noise, i.e. the equation: 
\begin{equation}\label{starting-equation}
(\imath \partial_t-\Delta)u=\lambda |u|^2+\dot{B}, \quad u_0=0,\quad \quad   t\in \R, \ x\in \mathbb{T},
\end{equation}
with $\mathbb{T}=\mathbb{R}/2\pi \mathbb{Z}$, $\lambda\in \{-1,1\}$ and where $\dot{B}$ stands for a space-time real valued noise, defined on a complete filtered probability space $(\Omega,\mathfrak{F},\mathbb{P})$. In this setting, we actually intend to pursue two main objectives, that can be roughly summarized as follows:

\smallskip

\noindent
$(i)$ Identify and treat situations where the  unknown $u$ cannot be treated as a function, but only as a space of general distributions, using some additional \textit{renormalization} procedure.

\smallskip

\noindent
$(ii)$ Try to go beyond  white-noise situations and handle the more sophisticated structure of a space-time \textit{fractional} noise. In other words, we would like to deal with a centered Gaussian noise $\dot{B}$ whose covariance is given for all test-functions $\vp,\psi$ on $\R^2$ by
\begin{equation}\label{cova-fractional-noise}
\mathbb{E}\big[\langle \dot B,\vp\rangle \,\langle \dot B,\psi\rangle\big]=\int_{\R^2}ds dt \int_{\R^2} dx dy\ \vp(s,x)\psi(t,y)|s-t|^{2H_0-2}|x-y|^{2H_1-2},
\end{equation}
for some indexes $H_0,H_1$ in a subset of $(0,1)$ to be specified.  Recall that the standard white noise  would correspond   to the case $H_0=H_1= \frac12$, and that the fractional noise is known for its fundamental past-dependence properties when $H_0,H_1\neq\frac12$ (see e.g. \cite{nourdin,nualart}).

\

The influence of randomness on nonlinear Schr{\"o}dinger (NLS) dynamics, and more generally on dispersive models, has been a topic of great interest and the source of an abundant literature in the past years. These investigations can essentially be structured around two main directions.

\

The first line of research, initiated by Bourgain \cite{bourgain-0,bourgain}, is concerned with the propagation of randomness within (deterministic) NLS dynamics, that is the influence of a random initial condition $\Phi$ on the general equation
\begin{equation}\label{general-deter-nls}
(\imath \partial_t-\Delta)u=\lambda\, \overline{u}^{p} u^q, \quad u(0,.)=\Phi, \quad t\in [0,T], \ x\in \R^d \ \text{or}\ \T^d,
\end{equation}
for fixed integers $p,q\geq 0$ and $\lambda\in \{-1,1\}$.  This approach has been developed by   Tzvetkov~\cite{Tzv, Tzv2} for the Schr\"odinger equation in various geometries, and extended  by Burq and Tzvetkov \cite{burq4, burq5,burq7} who considered more general random initial conditions for  the wave equation, and since then, there  have been many contributions to the study of dispersive equations with random initial conditions.  Such a randomization procedure is typically motivated by the possibility to improve the classical wellposedness results associated with~\eqref{general-deter-nls}  or other dispersive models (see e.g. \cite{burq4,thomann2,burq7,CO,BOP-1, BOP-2, BOP}), and the developments along this approach have also led to the exhibition of fundamental invariant measures for  dispersive equations, see e.g. \cite{bourgain-0,bourgain,Tzv, Tzv2,BT-W,burq5,BTT,De,TzvVi1, TzvVi2,DeTzVi,oh-thomann}.  Some progress have been made in the last few years regarding local wellposedness theory for \eqref{general-deter-nls} thanks to the works of Deng, Nahmod and Yue \cite{DNY-1,DNY-0,DNY}: the so-called theory of random tensors now provides us with an extensive covering of the   important  class of Hamiltonian equations (take $q=p+1$ in \eqref{general-deter-nls})
\begin{equation}\label{hamilton-eq}
(\imath \partial_t-\Delta)u=- |u|^{p}u, \quad u(0,.)=\Phi, \quad t\in [0,T], \ x\in \T^d,
\end{equation}
for $\Phi$ in a suitable (Gaussian-type) space of random distributions.  We also refer to the recent work \cite{BCST} on the probabilistic wellposedness of NLS on the sphere.

\

The second field of investigations, which can be traced back to works by de Bouard and Debussche, focuses on the influence of a random perturbation within the differential dynamics of the equation. This so-called stochastic forcing can take the form of an additive \cite{debou-debu-2002,debou-debu-2003}, multiplicative \cite{debou-debu-1999} or even dispersive~\cite{debou-debu-2010, debu-tsutsumi, ch-gubi, dubo-reve} noise. As far as applications are concerned, random perturbations can for instance be regarded as natural ways to model the influence of temperature variations on some molecular aggregation process~\cite{BCIRG}, or the evolution of nonlinear dispersive waves in a totally disorder medium \cite{conti,GGFDC}. In the context of optically trapped Bose-Einstein condensation, the fluctuation of laser intensity can be modeled by a temporal noise \cite{GAB}.

\smallskip

From a mathematical point of view, NLS models with stochastic forcing are part of the general SPDE theory \cite{DPZ,walsh}, just like stochastic heat equations or stochastic wave models. However, it is a well known fact that the properties of Schr{\"o}dinger propagation, and especially the relatively weak regularization effects in this setting, makes the analysis of stochastic Schrödinger equations particularly delicate. As a consequence, the literature related to this specific field turns out to be more scarce than its heat or wave counterpart.

\smallskip

Schr{\"o}dinger equations with random perturbation have been mostly investigated along a fully stochastic approach, that is in the presence of a white noise in time and through the use of sophisticated Itô-type controls. On top of the aforementioned papers by de Bouard and Debussche (see \cite[Lemma 3.1]{debou-debu-1999} or~\cite[Lemma 4.2]{debou-debu-2003} for examples of Itô-type lemmas), we can here refer to recent works by Hornung about stochastic NLS on~$\R^d$ \cite{hornung}, or to Brz{\'e}zniak-Millet \cite{BM} and Cheung-Mosincat \cite{cheung-mosincat} for results on compact domains (see e.g.  \cite[Theorem~3.10]{BM} and \cite[Lemma 3.6]{cheung-mosincat} for Itô-type estimates in such a framework). Observe that in all of these studies, the implementation of stochastic analysis is permitted by the fact that the equation can be interpreted and solved in a space of functions (typically a subspace of $L^2([0,T]\times \R^d)$ or $L^2([0,T]\times \T^d)$), leading to restrictive conditions on the noise spatial regularity.

\smallskip

Pathwise-type strategies - in the vein of  our subsequent developments - have also been developed in a few studies, especially in the multiplicative situation, where specific rescaling transformations are available. This is for instance the case in the results by de Bouard-Fukuizumi~\cite{debou-fukui}, Pinaud~\cite{pinaud} for a temporal noise, Barbu, Röckner and Zhang \cite{BRZ,zhang} for a white-in-time smooth-in-space noise on $\R^d$,  or in a recent series of papers by Debussche, Liu, Martin, Tzvetkov, Visciglia and Weber \cite{DLTV,debu-martin,debu-weber,TV} for a rough spatial noise.

\smallskip

Let us now turn to the presentation of several previous works related to the pathwise approach of NLS in the additive case, which is precisely the category our model \eqref{starting-equation} falls into - note indeed that our consideration of a fractional noise rules out any possible It{\^o}-type control, i.e. control in $L^p(\Omega)$, for the stochastic integral. 

\smallskip

First, and albeit restricted to a white noise, the two references \cite{oh-okamoto,oh-popovnicu-wang} by Oh, Okamoto, Pocovnicu and Wang offer us an insightful example of a pathwise procedure in the additive situation, with a clear separation between the stochastic arguments (see the study of the stochastic convolution therein) and the deterministic analysis of the problem. In both works \cite{oh-okamoto,oh-popovnicu-wang}, the noise under consideration happens to be sufficiently regular (in space) for the equation to be directly treated in a space of functions, with no need for any renormalization trick.

\smallskip

Such a renormalization trick can be found in the reference \cite{DST} about the NLS model on~$\R^d$ with additive fractional noise (see also \cite{schaeffer}). Although the latter setting may look similar to the one  in the present study (i.e. to the setting in \eqref{starting-equation}), it should be noted that the analysis in~\cite{DST} strongly relies on a specific local regularization effect of the Schr{\" o}dinger operator on~$\R^d$ (see~\cite{constantin-saut}), which is no longer available on the torus, making a huge difference between the two problems.

\

Last but not least, let us report on the results by Forlano, Oh and Wang in \cite{FOW} and Liu in~\cite{liu}. Although these two references are concerned with white-noise situations, their underlying dynamics are close to ours, which in particular will give us the opportunity to discuss about renormalization methods and ultimately specify our objective in the study.

\smallskip

In \cite{FOW}, the authors are interested in the following rescaled cubic NLS model with additive noise
\begin{equation}\label{sto-cubic-nls-renorm}
(\imath \partial_t-\Delta)\tilde{u}=-\Big(|\tilde{u}|^2-2\int_{\mathbb{T}}|\tilde{u}|^2\Big)\tilde{u}+\dot{B}, \quad \tilde{u}_0=\Phi\in L^2(\T), \quad (t,x)\in \R_+\times \T,
\end{equation}
while the author in \cite{liu} investigates a quadratic NLS model with random initial condition
\begin{equation}\label{liu-model}
(\imath \partial_t-\Delta)\tilde{u}=|\tilde{u}|^2-2\int_{\mathbb{T}}|\tilde{u}|^2, \quad \tilde{u}_0=\Phi\in L^2(\Omega;H^{-\al}(\T)), \quad (t,x)\in \R_+\times \T^2.
\end{equation}
In both \cite{FOW} and \cite{liu}, the roughness of the random data (either $\dot{B}$ in \eqref{sto-cubic-nls-renorm} or $\Phi$ in \eqref{liu-model}) justifies the involvement of a renormalization trick, namely
\begin{equation*} 
|u|^2 \longrightarrow |\tilde{u}|^2-2\int_{\mathbb{T}}|\tilde{u}|^2.
\end{equation*}

This transformation ensures the cancellation of specific resonant terms at the Fourier level, and it is in fact directly inspired by the renormalization used by Bourgain in \cite{bourgain} for the cubic NLS equation with irregular initial condition, that is the equation
\begin{equation}\label{cubic-nls-renorm}
(\imath \partial_t-\Delta)\tilde{u}=-\Big(|\tilde{u}|^2-2\int_{\mathbb{T}}|\Phi|^2\Big)\tilde{u}, \quad \tilde{u}_0=\Phi\in \cg, \quad (t,x)\in \R_+\times \T,
\end{equation}
where $\cg$ stands for a space of rough distributions (see also \cite{oh-thomann,DNY} for higher-order generalizations of this rescaling trick). In the analysis of the original equation \eqref{cubic-nls-renorm}, the choice of such a rescaling can be further justified through at least two reasons: 

\smallskip

\noindent
$(i)$ When working with a regular enough initial condition, say $\Phi\in L^2(\T)$, the two models~\eqref{hamilton-eq} (with $p=2$) and \eqref{cubic-nls-renorm} turn out to be equivalent. Indeed, in this case, it is readily checked that if~$\tilde{u}$ is a solution to \eqref{cubic-nls-renorm}, then $u:=e^{- 2\imath t \int_{\T}|\Phi|^2} \tilde{u}$ is a solution to the standard cubic NLS equation.

\smallskip

\noindent
$(ii)$ The renormalization term $2\int_{\mathbb{T}}|\Phi|^2$ in \eqref{cubic-nls-renorm} is (obviously) \textit{explicit}, i.e. it does not depend on the (unknown) solution $\tilde{u}$ but only on the (known) initial condition~$\Phi$, which offers a better control on the rescaling procedure.

\

With these two features in mind, recall that the main wellposedness result for \eqref{cubic-nls-renorm} can be summed up as follows (see e.g. \cite[Theorem 1]{CO}):

\begin{center}
\begin{minipage}{15cm}
\textit{Given $\Phi$ in a suitable space of random negative-order distributions on $\T$ (with full Gaussian measure), there exist a natural approximation $(\Phi^{(n)})_{n\geq 1}$ of $\Phi$ with $\Phi^{(n)}\in \cac^\infty(\T)$ and a time $\tau>0$ such that by setting $\Lambda^{(n)}:= 2\int_{\mathbb{T}}|\Phi^{(n)}|^2$, the sequence of classical solutions
\begin{equation}\label{cubic-nls-sequence}
(\imath \partial_t-\Delta)\tilde{u}^{(n)}=-\big(|\tilde{u}^{(n)}|^2- \Lambda^{(n)}\big)\tilde{u}^{(n)}, \quad \tilde{u}_0^{(n)}=\Phi^{(n)}, \quad t\in [0,\tau], \ x\in \T,
\end{equation}
is well defined and converges almost surely in $\cac([-\tau,\tau]; H^{-\al}(\T))$, for some $\al >0$.}
\end{minipage}

\end{center}

\smallskip

Due to the explicitness of the sequence $(\Lambda^{(n)})_{n\geq 1}$, the latter statement shows a strikingly similar pattern to other well-known renormalization results in the \textit{heat} or the \textit{wave} settings.

\smallskip

In the \textit{heat} setting, we can quote for instance (and among many other examples) the case of the dynamical $\Phi_3^4$ model, whose main wellposedness theorem reads as follows (see \cite[Theorem~1.15]{hairer}).

\begin{center}
\begin{minipage}{15cm}
\textit{Given a space-time white noise $\dot{B}$ on $\R_+\times \T^3$, there exist a smooth approximation $(\dot{B}^{(n)})_{n\geq 1}$ of $\dot{B}$, a time $\tau>0$, as well as \emph{an explicit sequence of constants $(\si^{(n)})_{n\geq 1}$ depending only on $(\dot{B}^{(n)})_{n\geq 1}$}, such that the sequence of  classical solutions  
\begin{equation}\label{example-heat}
(\partial_t-\Delta)\tilde{u}^{(n)}=-\big((\tilde{u}^{(n)})^2-\si^{(n)}\big)\tilde{u}^{(n)}+\dot{B}^{(n)}, \quad \tilde{u}_0^{(n)}=0, \quad t\in [0,\tau], \ x\in \T^3,
\end{equation}
is well defined and converges almost surely in $\cac([0,\tau];H^{-\al}( \T^3 ))$, for some $\al >0$.}
\end{minipage}
\end{center}

\smallskip

In the \textit{wave} framework, let us recall for instance the following result - obtained in \cite[Theorem 1.1]{gubi-koch-oh-2} - about the three-dimensional stochastic equation with quadratic nonlinearity (see also \cite{OOT1,OOT2,ORT} for similar results).

\begin{center}
\begin{minipage}{15cm}
\textit{Given a space-time white noise $\dot{B}$ on $\R_+\times \T^3$, there exist a smooth approximation $(\dot{B}^{(n)})_{n\geq 1}$ of $\dot{B}$, a time $\tau>0$, as well as \emph{a sequence of time-dependent constant $(\si^{(n)}(t))_{n\geq 1}$ explicitly described in terms of $(\dot{B}^{(n)})_{n\geq 1}$}, such that the sequence of classical solutions  
\begin{equation}\label{example-wave}
\partial^2_t \tilde{u}^{(n)}+(1-\Delta)\tilde{u}^{(n)}=-(\tilde{u}^{(n)})^2+\si^{(n)}+\dot{B}^{(n)}, \quad \tilde{u}_0^{(n)}=0, \quad t\in [0,\tau], \ x\in \T^3,
\end{equation}
is well defined and converges almost surely in $\cac([0,\tau];H^{-\al}( \T^3 ))$, for some $\al >0$.}
\end{minipage}
\end{center}

\smallskip


\begin{remark}
Although we have stuck to  white-noise-driven dynamics in both~\eqref{example-heat} and~\eqref{example-wave}, it should be noted that similar renormalization results can also be found in the fractional setting (see e.g. \cite{deya-heat,deya-wave}).
\end{remark}

\smallskip

Unfortunately, the explicitness property $(ii)$ satisfied by \eqref{cubic-nls-renorm}, and which allows such a close comparison between the {\textit{heat/wave} models \eqref{example-heat}-\eqref{example-wave} and the \textit{Schr{\"o}dinger} equation}~\eqref{cubic-nls-sequence}, is no longer true for both \eqref{sto-cubic-nls-renorm} and \eqref{liu-model}. Indeed, there is no reason for the renormalization term $2\int_{\mathbb{T}}|\tilde{u}|^2$ in \eqref{sto-cubic-nls-renorm} to be independent of the whole dynamics of the \enquote{unknown}~$\tilde{u}$.

\smallskip

In light of these elements, and going back to our fractional  Schr{\"o}dinger model \eqref{starting-equation}, we are now in a position to clarify the main objective of our study, namely: 
\begin{center}
\begin{minipage}{14cm}
\textit{In situations where the noise is too rough (i.e., $H_0,H_1$ in \eqref{cova-fractional-noise} are too small) for equation~\eqref{starting-equation} to be treated directly, we would like to find out an \emph{explicit} renormalization procedure, similar to the one in \eqref{cubic-nls-sequence}, \eqref{example-heat} or \eqref{example-wave}, allowing to offer a \emph{non-trivial interpretation} of the problem.}
\end{minipage}
\end{center}

\smallskip

This search for a renormalization term depending only on $\dot{B}^{(n)}$ can also be seen as an attempt to minimize the deformation of the model, as opposed to the fully-twisted product in \eqref{sto-cubic-nls-renorm}, or even to the trivial rescaling $(\imath \partial_t-\Delta)u^{(n)}=\lambda\, (|u^{(n)}|^2-\si^{(n)})+\dot{B}^{(n)}$ with $\si^{(n)}:=|u^{(n)}|^2$. In the same spirit, we shall obviously exclude from our considerations any renormalization method that would significantly hinder the effect of random perturbation, the most trivial of which being the example $(\imath \partial_t-\Delta)u^{(n)}=\lambda|u^{(n)}|^2-\si^{(n)}+\dot{B}^{(n)}$ with $\si^{(n)}:=\dot{B}^{(n)}$. 

\smallskip

To avoid such trivial situations, and in light of the rescaling procedures used in \eqref{cubic-nls-sequence}, \eqref{example-heat} and \eqref{example-wave}, we propose to restrict the class of admissible renormalization procedures through a few natural requirements. We will thus be able to prove the following result:

\begin{theorem}\label{main-obj}
Assume that the indexes $H_0,H_1$ satisfy the three conditions
\begin{equation}\label{cond-index-hurst}
\frac12<H_0<1, \quad \frac14<H_1<1 \quad \text{and} \quad \frac74<2H_0+H_1<2. 
\end{equation}
Then for a natural smooth approximation $(\dot{B}^{(n)})_{n\geq 1}$ of $\dot{B}$, there exist two sequences 
$$(\Lambda^{(n)})_{n\geq 1},(\si^{(n)})_{n\geq 1}$$
such that:

\smallskip

\noindent
$(i)$ For every $n\geq 1$, $\Lambda^{(n)}:\Omega \times \R \to \C $ is a time process, while $\si^{(n)}:\R \to \C$ is a deterministic time function.

\smallskip

\noindent
$(ii)$ For every $n\geq 1$, both $\Lambda^{(n)}$ and $\si^{(n)}$ depend only on $\dot{B}^{(n)}$.



\smallskip

\noindent
$(iii)$ There exists a (random) time $T_0>0$ such that the sequence $(\tilde{u}^{(n)})_{n\geq 1}$ of classical solutions to the renormalized equation
\begin{equation}\label{renorm-class}
(\imath \partial_t-\Delta)\tilde{u}^{(n)}=\lambda\Big[\overline{\tilde{u}^{(n)}}-\Lambda^{(n)}\Big] \tilde{u}^{(n)}-\si^{(n)}+\dot{B}^{(n)}
\end{equation}
converges almost surely to some limit $\tilde{u}$ in the space of distributions on $[-T_0,T_0]\times \T$.

\end{theorem}

We refer to Theorem \ref{theo:well-posed-z} and Theorem \ref{theo:seq} for more quantitative results.  

\smallskip

The above properties $(i)$-$(iii)$ can thus be read as follows. Property $(i)$ ensures the non-triviality of the renormalization method: $\Lambda^{(n)}$ is reduced to a time process (independent of $x$), while $\si^{(n)}$ is just a time function. Property $(ii)$ emphasizes the explicit nature of the procedure, a major improvement over existing methods: $\Lambda^{(n)}$ and $\si^{(n)}$ only depend on the data we are given a priori, and not on $(u^{(n)})_{n\geq 1}$. Finally, property $(iii)$ guarantees the success of the method, by restoring the (asymptotic) wellposedness of the problem.

\smallskip

As indicated in item $(iii)$, our analysis will only focus on the exhibition of a \textit{local} solution, satisfying the equation on a (possibly small) time interval $[-T_0,T_0]$, with $T_0>0$. In particular, the choice of $\lambda\in \{-1,1\}$ in \eqref{starting-equation} or \eqref{renorm-class} will not play any essential role in our analysis: for simplicity, {\it we assume from now on that $\lambda=1$}.

\begin{remark}
Observe that our condition \eqref{cond-index-hurst} on $(H_0,H_1)$ can still cover situations where $\dot{B}$ is rougher than a white noise \textit{in space}, that is for which $H_1<\frac12$. On the other hand, $\dot{B}$ must always be more regular than a white noise \textit{in time}, i.e. $H_0>\frac12$, for our subsequent arguments to hold.
\end{remark}

\begin{remark}
The result of Proposition \ref{prop:conv-psi-n} below ensures that the assumption $2H_0+H_1<2$ in \eqref{cond-index-hurst} corresponds to a \textit{singular} situation, that is, one in which the solution $u$ cannot be defined as a function, and hence where the interpretation of the nonlinearity necessarily requires a renormalization procedure.\\
\indent On the other hand, we are not in a position, at this stage, to assert that the condition $2H_0+H_1>\tfrac74$ in \eqref{cond-index-hurst} is optimal. This restriction is mainly used in the proof of Proposition~\ref{prop:ice-cream}, that is, in the construction of the second-order process arising from the dynamics (see Remark \ref{rk:optim-ice-cream} for further details).
\end{remark}

\begin{remark}

Although our subsequent study leaves the global-wellposedness issue fully open, it should be noted that the deterministic NLS equation with $|u|^2$ non-linearity is globally ill-posed for any non-trivial initial condition (see \cite{FG}). Thus the equation \eqref{starting-equation} may blow up as well.
\end{remark}

\subsection{Reformulation of the equation and main steps of the strategy}\label{subsec:reform-main-steps}

Let us now provide some details about the main steps of our strategy to achieve the above purpose. With {properties $(i)$-$(iii)$} of Theorem \ref{main-obj} in mind, we henceforth focus on the rescaled model 
\begin{equation}\label{starting-equation-resc}
(\imath \partial_t-\Delta)u=\big(\overline{u}-\Lambda\big) u-\si+\dot{B}, \quad u_0=0,\quad \quad   t\in \R, \ x\in \mathbb{T},
\end{equation}
for some fixed functions $\Lambda:\Omega\times \R\to \C$ and  {$ \si: \R \to \C$}, whose expression will be determined later at the level of an approximated equation.

\smallskip

Since we are only interested in local solutions (in time), we introduce and \textit{fix once and for all} a smooth temporal cut-off function $\chi:\R\to [0,1]$ such that

\smallskip

\noindent
$(i)$ $\chi$ is symmetric, i.e. $\chi(t)=\chi(-t)$ for all $t\in \R$, 

\smallskip

\noindent
$(ii)$ $\chi(t)=1$ for all $t\in [-1,1]$, and $\chi(t)=0$ for all $|t|\geq 2$.

\smallskip

For every time $\tau>0$, we will also denote by $\chi_\tau$ the function defined for $t\in \R$ by $\chi_\tau(t):=\chi\big(\frac{t}{\tau}\big)$.

\smallskip

With this notation, let us recast the dynamics in \eqref{starting-equation-resc} under the following (localized) mild form:
\begin{align}
&u(t)=-\imath \chi(t)\int_0^t ds\, \chi(s) e^{-\imath (t-s)\Delta}\dot{B}(s,.)-\imath \chi_\tau(t)\chi(t)\int_0^t ds\, e^{-\imath (t-s)\Delta}\chi(s)  |u(s,.)|^2\nonumber\\
&+\imath \chi_\tau(t)\chi(t)\int_0^t ds\, \chi(s)\Lambda(s)e^{-\imath (t-s)\Delta} u(s,.)+\imath \chi_\tau(t)\chi(t)\int_0^t ds\, \chi(s)e^{-\imath (t-s)\Delta} \si(s) .\label{starting-equation-chi}
\end{align}

The mild formulation \eqref{starting-equation-chi} clearly emphasizes the role of the so-called stochastic convolution in the dynamics, and accordingly let us set (at least formally, since $\dot{B}$ is only defined as a distribution)
\begin{equation}\label{defi-luxo}
\<Psi>(t,x):=-\imath \chi(t)\int_0^t ds \, \chi(s) \big(e^{\imath s\Delta}\dot{B}(s)\big)(x), \quad t\in \R, \ x\in \T,
\end{equation}
or otherwise stated
$$\<Psi>(t,x)=-\imath \sum_{k\in \Z}\bigg[  \chi(t)\int_0^t ds \, \chi(s) e^{-\imath s |k|^2}\int_{\T} dy \, e^{-\imath k y} \dot{B}(s,y)  \bigg] e^{\imath k x} , \quad t\in \R, \ x\in \T,$$
where we use the notation $\dis \int_{\T}dx\, \vp(x)=\frac1{2\pi}\int_{0}^{2\pi}dx\, \vp(x)$.

Introduce the (local) time integration operator, that is 
\begin{equation}\label{defi:ci-chi}
 \ci_\chi v(t):= - \imath \chi(t)\int_0^t ds \, \chi(s)v(s),
\end{equation}
as well as the auxiliary process
\begin{equation}\label{def-uring}
\mathring{u}(t):=e^{\imath  t \Delta} u(t).
\end{equation}
Then \eqref{starting-equation-chi} is equivalent to 
\begin{equation}\label{basic-equation0}
\mathring{u}(t)=\<Psi>(t)+\chi_\tau(t) \ci_\chi \cm(\mathring{u},\mathring{u})(t)- \chi_\tau(t)\ci_\chi\big(\Lambda \cdot\mathring{u}\big)(t)- \chi_\tau(t)\ci_\chi\big({ \si} \big)(t), \quad t\in \R,
\end{equation}
where the bilinear operator $\cm$ is defined by
\begin{equation}\label{defm}
 \cm(v,w)(t) := e^{\imath t \Delta} \big( (e^{-\imath t \Delta} v(t)) \overline{(e^{-\imath t \Delta} w(t))}   \big),
\end{equation}
and the product $\Lambda \cdot\mathring{u}$ is basically understood as $(\Lambda \cdot\mathring{u})(t,x):=\Lambda(t)\mathring{u}(t,x)$.

\smallskip

Our next transformation of the model is often referred to as the Da Prato-Debussche trick in the SPDE literature (in reference to the techniques used in \cite{daprato-debussche}), but similar ideas can also be found in Bourgain's work \cite{bourgain} about a Schr{\"o}dinger model. Namely, we introduce the \textit{remainder} process $\mathring{z}$ defined  by
\begin{equation*}
\mathring{z}(t):=\mathring{u}(t)- \<Psi>(t), \quad t\in \R,
\end{equation*} 
and thus recast equation \eqref{basic-equation0} into the \textit{remainder} equation
$$\mathring{z}=\chi_\tau \cdot \ci_\chi \cm(\mathring{z}+ \<Psi>,\mathring{z}+\<Psi>)-\chi_\tau\cdot \ci_\chi\big(\Lambda \cdot (\mathring{z}+\<Psi>)\big)-\chi_\tau \cdot \ci_\chi\big({ \si}  \big),$$
or otherwise stated 
\begin{equation}\label{equa-z}
\begin{split}
&\mathring{z}=\chi_\tau \cdot \ci_\chi \cm(\mathring{z},\mathring{z})+\chi_\tau \cdot \Big[\ci_\chi \cm(\mathring{z}, \<Psi>)+ \ci_\chi \cm(\<Psi>,\mathring{z})\Big]\\
&\hspace{2cm}-\chi_\tau\cdot \ci_\chi\big(\Lambda \cdot \mathring{z}\big)+\chi_\tau \cdot \Big[\ci_\chi \cm( \<Psi>,\<Psi>)- \ci_\chi\big(\Lambda \cdot \<Psi>\big)-\ci_\chi\big({ \si} \big)\Big]\, .
\end{split}
\end{equation}

\smallskip

At this stage, and based on the latter expansion, let us clarify the general objectives raised in Section~\ref{subsec:gene-conside} through the following two-step procedure:

\smallskip

\textit{Step 1.} Identify situations where the fractional noise is so rough, i.e. the indexes $H_0,H_1$ in~\eqref{cova-fractional-noise} are small enough, that the central process $\<Psi>$ - and accordingly the solution $u$ itself - cannot be treated as a function on $\R_+\times \mathbb{T}$ (for fixed $\omega\in \Omega$), making impossible the interpretation of the square term $|u|^2$ in~\eqref{starting-equation-chi} as a standard product of functions.

\smallskip

\textit{Step 2.} In some of the situations identified in Step 1, try to find out suitable expressions for $\Lambda$ and~$\si$, depending only on $\dot{B}$ - or equivalently on $\<Psi>$ -, so that \eqref{equa-z} admits a non-trivial interpretation and allows us to set up \textit{a fixed point argument in a space of functions}.

\

The treatment of Step 1 will actually reduce to relatively standard moments estimates on~$\<Psi>$ (see the subsequent Proposition \ref{prop:conv-psi-n} and its proof). Therefore, as the reader might expect it, the bulk of our work will consist in the analysis of Step 2.\\
\indent In brief, our first idea to achieve Step 2 will be to concentrate the renormalization procedure on the last expression into brackets in \eqref{equa-z}, that is the explicit (but a priori ill-defined) term  
\begin{equation}\label{resc-ice}
\ci_\chi \cm( \<Psi>,\<Psi>)- \ci_\chi\big(\Lambda \cdot \<Psi>\big)-\ci_\chi\big({ \si}  \big),
\end{equation}
which ultimately offers the guarantee of an explicit rescaling method (see item $(iii)$ of Theorem \ref{main-obj}). The result of this procedure is summed up in the statement of Proposition \ref{prop:ice-cream} below, and it gives birth to a well-defined renormalization $\<IPsi2>$ of $\ci_\chi \cm( \<Psi>,\<Psi>)$. In turn, we can inject this new object into \eqref{equa-z} and derive the final formulation of our problem, namely the renormalized equation
\begin{equation}\label{equa-z-re}
z=\chi_\tau \cdot \ci_\chi \cm(z,z)+\chi_\tau \cdot \Big[\ci_\chi \cm(z, \<Psi>)+ \ci_\chi \cm(\<Psi>,z)\Big]-\chi_\tau\cdot \ci_\chi\big(\Lambda \cdot z\big)+\chi_\tau \cdot \<IPsi2>\, .
\end{equation}

\smallskip

The above general rescaling strategy (i.e., the fact that we focus on $\ci_\chi \cm( \<Psi>,\<Psi>)$ only) does not significantly deviate from the one used in~\cite{gubi-koch-oh}, \cite{deya-wave} or \cite{OO} within similar \enquote{Da Prato-Debussche-trick} procedures - except for the shape of our Bourgain-Wick transformation in~\eqref{defi:ipsi2-main}. What actually makes the subsequent analysis very different from those in \cite{gubi-koch-oh,deya-wave,OO}, and what will be the source of our main efforts, is the interpretation and the control of the - \textit{non-renormalized} - product terms $\ci_\chi \cm(z, \<Psi>)$, $\ci_\chi \cm(\<Psi>,z)$ in \eqref{equa-z-re}.

\smallskip

Indeed, in contrast with the heat or the wave situation, one must here cope again with the fundamental difficulties related to the weak regularization properties of the Schr{\"o}dinger operator: since $\<Psi>$ is only a negative-order distribution, we cannot rely on any known (deterministic) effect of $\ci_\chi \cm$ to bring the product term $\ci_\chi \cm(z,\<Psi>)$ back into a space of functions (as regular as $z$ may be).

\smallskip

In order to overcome this problem, we essentially have no choice but to consider the product operations
\begin{equation*} 
 \cl^{\<Psi>,+}_\chi: z\mapsto \ci_\chi \cm(z,\<Psi>) \quad \text{and}\quad   \cl^{\<Psi>,-}_\chi:z\mapsto\ci_\chi \cm(\<Psi>,z)
\end{equation*}
as \textit{random operators}, with stochastically-controlled operator norms. At a heuristic level, our aim will thus be to highlight random quantities $\cq^{\<Psi>,\pm}_{\chi}$ whose moments can be explicitly estimated, and so that  for every function $z$, 
\begin{equation*}
\big\|\cl^{\<Psi>,\pm}_\chi (z)\big\|\leq \cq^{\<Psi>,\pm}_{\chi}\,  \|z\|,
\end{equation*}
where $\|.\|$ stands for a functional norm to be suitably chosen. These investigations will be the topic of Section \ref{sec:estim-random-op}, the main result of which is reported in Proposition \ref{prop:control-op-l-intro} below.

\smallskip

Fourier analysis will play a prominent role in our study, and for this reason we introduce the following convention for spatial Fourier transform.
\begin{notation}\label{not:four-space}
For every function $\vp:\mathbb{T}\to \R$ and every $k\in \mathbb{Z}$, we will denote
$$\vp_k:=\int_{\T}dx \, e^{-\imath\, k x}\vp(x)=\frac1{2\pi}\int_0^{2\pi}dx \, e^{-\imath\, k x}\vp(x)  .$$
\end{notation}

Using this convention, let us emphasize the expression of the bilinear operator $\cm$ in Fourier modes (we refer to the appendix for the proof).
\begin{lemma}\label{Lemma-M}
For any $v,w \in \mathcal{S}'(\T)$, we have that for all $k\in \Z$,
\begin{equation} \label{cm-fourier}
\cm(v,w)_k(t)=\sum_{k_1\in \Z} e^{\imath t\Omega_{k,k_1}} v_{k+k_1} \overline{w_{k_1}},  
\end{equation}
with $\Omega_{k,k_1}:=|k+k_1|^2-|k_1|^2-|k|^2=2k k_1$.
\end{lemma}

\section{Main results}

We are now in a position to elaborate on our main results about the model \eqref{starting-equation}, in a specific rough regime induced by the space-time fractional noise. 

\smallskip

 The consideration of such a noise - famous for its past-dependence properties (see e.g. \cite{nourdin,nualart}) - has become a standard topic in the SDE setting. It has also given birth to a quite extensive literature in the heat or the wave settings (see e.g. \cite{hu-nualart-song,deya-heat,balan-conus,deya-wave}). It should be noted, on the other hand, that this line of investigations still remains relatively unexplored for Schr{\"o}dinger equations (see e.g. \cite{DST}).

\subsection{Fractional noise on $\R\times \T$}

The space-time fractional noise is classically defined on~$\R^2$, using the covariance formula \eqref{cova-fractional-noise}. We will here handle its natural counterpart on $\R\times \T$, defined through a basic periodization procedure in space.

\smallskip

In order to motivate this definition, let us briefly go back to the classical white noise case and consider a white noise $\dot{W}$ on the whole Euclidean space $\R^2$, that is a centered Gaussian noise with covariance (formally) given for $(t,x),(t',x')\in \R^2$ by 
$$
\mathbb{E}\big[\dot{W}(t,x) \dot{W}(t',x')\big]=\delta_{\{t=t'\}} \delta_{\{x=x'\}},
$$
where $\delta$ stands for the Dirac mass distribution. Then by setting, for every $(t,x)\in \R\times \T$,
$$\dot{\mathbf{W}}(t,x) :=\sum_{k\in \Z} \bigg( \int_{\T}dy\, e^{-\imath k y}  \dot{W}(t,y)\bigg) e^{\imath k x},$$
one can check that for all test-functions $\vp,\psi$ on $\R\times \T$, 
\begin{align*}
&\mathbb{E}\Big[ \langle \dot{\mathbf{W}},\vp\rangle\, \overline{\langle \dot{\mathbf{W}},\psi\rangle}\Big]\\
&=\sum_{k,k'\in \Z}\int_{\R^2}dtdt'\int_{\T^2} dx dx'\,  \vp(t,x)\overline{\psi(t',x')}e^{\imath k x}e^{-\imath k' x'} \int_{\T^2}dydy'\, e^{-\imath k y}e^{\imath k' y'} \mathbb{E}\big[ \dot{W}(t,y)\dot{W}(t',y')\big]\\
&=\int_{\R}dt \sum_{k\in \Z}\bigg(\int_{\T} dx \, e^{\imath k x} \vp(t,x)\bigg)\bigg(\int_{\T} dx'\, \overline{\psi(t,x')}e^{-\imath k x'}\bigg)=\int_{\R\times \T} dt dx \, \vp(t,x)\overline{\psi(t,x)},
\end{align*}
which corresponds precisely to the covariance of a space-time white-noise on $\R\times \T$.

\smallskip

Keeping this general injection procedure in mind, we define, for all indexes $H_0,H_1\in (0,1)$ and all $(t,x)\in \R\times \T$,
\begin{equation}\label{dot-b-formal}
\dot{\mathbf{B}}(t,x) :=\sum_{k\in \Z} \bigg( \int_{\T}dy\, e^{-\imath k y}  \dot{B}(t,y)\bigg) e^{\imath k x},
\end{equation}
where, in the right-hand side, $\dot{B}$ is a space-time fractional noise on $\R^2$ with indexes $H_0,H_1$, as defined through the covariance formula \eqref{cova-fractional-noise}.\\
\indent In the sequel, we will make no distinction between $\dot{B}$ and its space-periodization $\dot{\mathbf{B}}$.

\begin{remark}\label{rk:compa-noise}

From a stochastic point of view, the above-introduced fractional noise is fundamentally different from the spatial regularization of white noise used for instance in \cite{FOW,OWZ}, and whose independence properties are essentially the same as those of $\dot{\mathbf{W}}$. This gap can easily be observed through the behavior of the corresponding (rescaled) linear solution, as developed in Remark \ref{rk:luxo-white} below.

\end{remark}

\subsection{Approximation of the noise}\label{subsec:approx-noise}

For a rigorous treatment of $\dot{B}$ in \eqref{dot-b-formal}, as well as a suitable interpretation of the related processes $\<Psi>$, $\ci_\chi \cm(z, \<Psi>)$, $\ci_\chi \cm(\<Psi>,z)$ and $\<IPsi2>$ in \eqref{equa-z-re}, we propose to follow an approximation procedure. In other words, we intend to define these objects as the limits of the corresponding processes above a smooth approximation $\dot{B}^{(n)}$ of~$\dot{B}$.

Owing to its convenience in the context of Fourier analysis, we will here consider the so-called harmonizable approximation of $\dot{B}$. To be more specific, we first introduce the approximation $B^{(n)}$ of a fractional sheet $B$ given for all $n\geq 1$ and $(t,x)\in \R^2$ by the formula
\begin{equation}\label{defi:approx-sheet}
B^{(n)}(t,x):= c_H\int_{|\xi|\leq 2^{2n}}^{}\int_{|\eta|\leq 2^n} \frac{e^{\imath t\xi}-1}{|\xi|^{H_0+\frac{1}{2}}}\frac{e^{\imath x\eta}-1}{|\eta|^{H_1+\frac{1}{2}}}\, \widehat{W}(d\xi,d\eta) 
\end{equation}
for some constant $c_H>0$, and where $\widehat{W}$ is the Fourier transform of a  Brownian field~$W$ on~$\R^{2}$ (defined on some complete filtered probability space $(\Omega,\mathcal{F},\mathbb{P})$).

\smallskip

Due to the restricted integration domain in \eqref{defi:approx-sheet}, it is readily checked that for every fixed $n\geq 1$, $(t,x)\mapsto B^{(n)}(t,x)$ is a smooth real-valued process on $\R^2$. Besides, it is a well-known fact that for an appropriate choice of $c_H$ (that we fix from now on), the sequence $(B^{(n)})_{n\geq 1}$ converges to $B$ in the space $\cac(K)$ of a continuous function on $K$, for any compact subset $K\subset \R^2$, where $B$ is a space-time fractional sheet of index $H=(H_0,H_1)$.

\smallskip

Then we set
\begin{equation}\label{defi:approx-b-n}
\dot{B}^{(n)}(t,x):=\frac{\partial^2 B^{(n)}}{\partial t\partial x}(t,x)= -c_H\int_{|\xi|\leq 2^{2n}}^{}\int_{|\eta|\leq 2^n} \frac{\xi\, e^{\imath t\xi}}{|\xi|^{H_0+\frac{1}{2}}}\frac{\eta\, e^{\imath x\eta}}{|\eta|^{H_1+\frac{1}{2}}}\, \widehat{W}(d\xi,d\eta) .
\end{equation}
Since $B^{(n)} \to B$, we can immediately guarantee the almost sure distributional convergence of $\dot{B}^{(n)}$ to a space-time fractional noise $\dot{B}$ of index $H=(H_0,H_1)$, as desired.

\subsection{Stochastic tree-elements}\label{subsec:main-res-sto}

With the smooth approximation $(\dot{B}^{(n)})_{n\geq 1}$ in hand, we define the approximated version of the process $\<Psi>$ at the core of \eqref{equa-z-reno} (and formally given by~\eqref{defi-luxo}) as
\begin{equation}\label{defi-luxo-n}
\<Psi>^{(n)}(t):=- \imath \chi(t)\int_0^t ds \, \chi(s) e^{\imath s\Delta}\dot{B}^{(n)}(s),
\end{equation}
or otherwise stated
$$\<Psi>^{(n)}(t,x)=- \imath \sum_{k\in \Z}\bigg[ \chi(t)\int_0^t ds \, \chi(s) e^{\imath s |k|^2}\int_{\T} dy \, e^{-\imath k y} \dot{B}^{(n)}(s,y)  \bigg] e^{\imath k x} , \quad t\in \R, \ x\in \T.$$

Our wellposedness criterion for $\<Psi>$ will then be based on the following convergence results, the proof of which will be examined in the subsequent Section \ref{subsec:first-order} (see Notation \ref{notation:x-b-c} below for the definition of the Sobolev spaces $H^{- \al }(\T)$ and $W^{-\al,p}(\T)$, $\al\in \R$, $1\leq p \leq \infty$). 

\begin{proposition}\label{prop:conv-psi-n}
Assume that $H_0,H_1\in (0,1)$ satisfy $\frac32 <2H_0+H_1< 2$.

\smallskip

Then the following assertions hold true.

\smallskip

\noindent
$(i)$ For every $\al > 2- (2H_0+H_1)$, the sequence $(\<Psi>^{(n)})_{n\geq 1}$ converges a.s. to some limit $\<Psi>$ in~$\mathcal{C}(\R;W^{- \al ,p}(\T))$ for all $1 \leq p \leq \infty$.

\smallskip

\noindent
$(ii)$ For every $\al < 2- (2H_0+H_1)$ and  every $T>0$, one has $\mathbb{E}\Big[\|\<Psi>^{(n)}\|_{L^2([0,T] ; H^{-\alpha}( \mathbb{T}))}^2\Big] \to \infty$ as $n\to\infty$.

\end{proposition}

The above statement allows us to conclude that when $2H_0+H_1< 2$, the process $\<Psi>$ - which governs the regularity of the equation - can only be defined as a distribution of negative order. The treatment of such situations is precisely the main objective of this paper, and accordingly, \textit{from now on, we will only focus on cases for which $2H_0+H_1<2$.}

\begin{remark}\label{rema}
It could be proved that when $2H_0+H_1>2$, the process $\<Psi>$ becomes a well-defined function (almost surely), and in this case we could rely on the classical Strichartz estimates to ensure straightforward wellposedness of \eqref{starting-equation} (see e.g. \cite{cheung-mosincat} or \cite[Section 3]{DST} for treatments of such situations).\\
\indent Note  that in the specific case where $2H_0+H_1=2$, the conclusions of our study (from Proposition~\ref{prop:conv-psi-n} to Theorem~\ref{theo:seq}) would essentially remain the same. However, taking this situation into account would require a number of technical adjustments in the (numerous) estimates of the paper, and so we have preferred to avoid such complications.  
\end{remark}

As a fundamental consequence of the fact that $\<Psi>$ cannot be regarded as a function, the interpretation of the \enquote{square} element $\ci_\chi\cm( \<Psi>,\<Psi>)$ in \eqref{equa-z} is not clear at first sight. In fact, such an interpretation can only be obtained through a renormalization procedure.

\smallskip

In this stochastic Schr{\"o}dinger setting, we are naturally led to the following renormalization of~$\ci_\chi \cm( \<Psi>,\<Psi>)$, defined at level of the approximation $\ci_\chi \cm( \<Psi>^{(n)},\<Psi>^{(n)})$ along two steps. 

\smallskip

We will first lean on the correction term $-\ci_\chi\big(\Lambda^{(n)} \cdot \<Psi>^{(n)}\big)$ in \eqref{equa-z} and therein choose 
\begin{equation}\label{defi-lambda-n}
\Lambda^{(n)}(t)=\int_{\T} \overline{\<Psi>^{(n)}}(t):=\int_{\T} \overline{\<Psi>^{(n)}(t,x)},
\end{equation}
so that 
\begin{equation}\label{hat-cm}
\ci_\chi \cm( \<Psi>^{(n)},\<Psi>^{(n)})-\ci_\chi\big(\Lambda^{(n)} \cdot \<Psi>^{(n)}\big)=\ci_\chi\Big(\cm( \<Psi>^{(n)},\<Psi>^{(n)})-  \<Psi>^{(n)} \int_{\T} \overline{\<Psi>^{(n)}}\Big)=: \<IPsi2dot>^{(n)}.
\end{equation}
The whole point of this transformation lies in the simplification that takes place at the Fourier level: namely, using the expression \eqref{cm-fourier} of $\cm$, we obtain that 
\small
\begin{eqnarray}
\bigg[\cm( \<Psi>^{(n)},\<Psi>^{(n)})-  \<Psi>^{(n)} \int_{\T} \overline{\<Psi>^{(n)}}\bigg]_k(t)&=&\cm( \<Psi>^{(n)},\<Psi>^{(n)})_k(t)-  \<Psi>^{(n)}_k(t) \overline{\<Psi>^{(n)}_0(t)}\nonumber\\
&=&\sum_{k_1\neq 0} e^{\imath t\Omega_{k,k_1}} \<Psi>^{(n)}_{k_1+k}(t) \overline{\<Psi>^{(n)}_{k_1}(t)},\label{defi-m-tilde-fou}
\end{eqnarray}
\normalsize
and so the only resonant term in the product, i.e. the one for which $\Omega_{k,k_1}=0$, occurs at $k=0$.

\smallskip

Secondly, following a standard scheme in any {$L^2$-setting}, we proceed with an additional Wick renormalization of the element $\<IPsi2dot>^{(n)}$ in \eqref{hat-cm}. Namely, we consider the sequence of processes defined for all $n\geq 1$ by
{\begin{equation}\label{defi:ipsi2-main}
\<IPsi2>^{(n)}(t,x):=\<IPsi2dot>^{(n)}(t,x)-\mathbb{E}\bigg[\int_{\T} dy\, \<IPsi2dot>^{(n)}(t,y)\bigg]\, .
\end{equation}
Using the extra correction term $-\ci_\chi\big(\si \big)$ in \eqref{resc-ice}, the Wick reduction \eqref{defi:ipsi2-main} can in fact be derived from the specific choice
\begin{equation}\label{defi-si-n}
\si^{(n)}(t):=\mathbb{E}\bigg[\int_{\T} dy\, \bigg(\cm( \<Psi>^{(n)},\<Psi>^{(n)})(t)-  \<Psi>^{(n)}(t) \int_{\T} \overline{\<Psi>^{(n)}}(t)\bigg)(y)\bigg].
\end{equation}
Thanks to the Fourier expression \eqref{defi-m-tilde-fou}, the latter quantity can also conveniently written as
\begin{equation}\label{defi-si-n-bis}
\si^{(n)}(t)=\mathbb{E}\bigg[\bigg(\cm( \<Psi>^{(n)},\<Psi>^{(n)})(t)-  \<Psi>^{(n)}(t) \int_{\T} \overline{\<Psi>^{(n)}}(t)\bigg)_0\bigg]=\sum_{k\neq 0} \mathbb{E}\Big[ \big| \<Psi>^{(n)}_k(t)\big|^2\Big].
\end{equation}
}

\smallskip

Let us now examine the convergence of the sequence $(\<IPsi2>^{(n)})_{n\geq 1}$ in a scale of Sobolev spaces suited to our problem.

\begin{notation}\label{notation:x-b-c}
 Recall Notation  \ref{not:four-space} for the Fourier transform in space. For $s\in \R$ and $1 \leq p \leq \infty$, we define the Sobolev space $W^{s,p}(\T)$ by
$$W^{s,p}(\T) =\big\{ u\in \mathcal{S}'(\T) : \; (1-\Delta)^{s/2}u  \in L^p(\T)\big\},$$
endowed with the natural norm. In the particular  case $p=2$ we set $H^s(\T)=W^{s,2}(\T)$ and we have
$$\| \vp \|^2_{H^s(\T)}=\sum_{k \in \Z} \langle  k \rangle ^{2s} |\vp_k|^2.$$

The natural framework for the study of \eqref{starting-equation} is given by the Bourgain spaces $X^{s,b}=X^{s,b}(\R \times \T)$ which are defined by  completion of the set of test-functions on $\R \times \mathbb{T}$ with respect to the norm
\begin{equation}\label{bourg}
\| u\|^2_{X^{s,b}}:=\sum_k \, \langle k\rangle^{2s}\int d\la \, \langle \lambda -|k|^2\rangle^{2b}  \big|\cf(u_k)(\la)\big|^2,
\end{equation}
where $\cf$ refers to the Fourier transform in time, and where we have used Notation \ref{not:four-space}. Equivalently, through the change of variable $u=e^{- \imath t \Delta}z$, we will work in the space $Z^{s,b}$ defined by 
\begin{equation}\label{norm:x-b-c}
\|z\|_{Z^{s,b}}^2= \|e^{- \imath t \Delta}z\|_{X^{s,b}}^2:=\sum_{k\in \Z} \langle k\rangle^{2s}\int_{\R}d\la \, \langle \la \rangle^{2b}\big|\cf(z_k)(\la)|^2.
\end{equation}

For a pedagogical introduction to the Bourgain spaces, we refer to \cite{Ginibre} or to \cite{Erdo-Tzi}.
\end{notation}

\begin{proposition}\label{prop:ice-cream}

Assume that $\frac12 < H_0<1$ and  $\frac14< H_1<1$ satisfy  $\frac74 <2H_0+H_1< 2$, and set $\underline{H}:=2H_0+H_1-1\in (\frac34,1)$. Consider the sequence of processes defined for every $n\geq 1$ by
\begin{equation}\label{defi-i-psi-2-n}
\<IPsi2>^{(n)}:=\ci_\chi \cm( \<Psi>^{(n)},\<Psi>^{(n)})- \ci_\chi\big(\Lambda^{(n)} \cdot \<Psi>^{(n)}\big)-\ci_\chi\big({ \si^{(n)}} \big)
\end{equation}
with $\Lambda^{(n)}$ and $\si^{(n)}$ given by \eqref{defi-lambda-n} and \eqref{defi-si-n-bis}, respectively.

\smallskip

Then the following assertions hold true.

\smallskip

\noindent
$(i)$  For every pair $(s,b)$ satisfying the conditions
\begin{equation}\label{assump-lau-statement}
0<s<\underline{H}-\frac12 \quad \text{and} \quad \frac12<b< b_s, \quad b_s:=\min\big({\underline{H}}-s,\frac34-\frac{s}{2},{2\underline{H}-1}\big),
\end{equation}
the sequence $(\<IPsi2>^{(n)})_{n\geq 1}$ converges a.s. in $Z^{s,b}$. We denote its limit by $\<IPsi2>$.

\smallskip

\noindent
$(ii)$ {For every $t>0$, the sequence $(\si^{(n)}(t))_{n\geq 1}$ is asymptotically equivalent to
$$\si^{(n)}(t) \stackrel{n\to \infty}{\sim} c_{\underline{H}}\, t\, 2^{2n(1-\underline{H})},$$
for some constant $c_{\underline{H}}>0$.}

\end{proposition}

{The proof of item $(i)$ will be the topic of Section \ref{sec:tree-elements-2}, while we have postponed the proof of item $(ii)$ to the appendix section \ref{prop-ice-cream-ii}.}

\smallskip

In light of these two properties, we thus have the following picture:

\smallskip

\noindent
$(i)$ When restricting the roughness domain of the noise $\dot{B}$ to the case $\frac74<2H_0+H_1<2$ (to be compared with the condition $\frac32<2H_0+H_1<2$ in Proposition \ref{prop:conv-psi-n}), we can indeed construct a renormalized version of the product term $\ci_\chi\cm( \<Psi>,\<Psi>)$ in a space $Z^{s,b}$ adapted to our purpose.

\smallskip

\noindent
$(ii)$ The constant (i.e., non-random) sequence $\big(\si^{(n)}\big)_{n\geq 1}$ at the core of this rescaling method {explodes as $n\to \infty$}, which provides an additional evidence of how necessary the renormalization step.

\

\begin{remark}\label{rk:optim-ice-cream}
As we shall see in Section \ref{sec:tree-elements-2}, the proof of Proposition \ref{prop:ice-cream} relies on the control of sophisticated fractional expressions arising from the covariance of $\<Psi>^{(n)}$, and at this stage we are not in a position to guarantee the optimality of the condition $2H_0+H_1>\tfrac74$ appearing in the statement.\\
\indent In fact, we believe that this condition could be improved at the cost of a more demanding renormalization procedure: a result in this direction can be found in \cite[Proposition 2.2]{deya-restric}, in the particular case $H_1=\tfrac12$.
\end{remark}

\

Before we turn to the treatment of the product terms $\ci_\chi \cm(z,\<Psi>),\ci_\chi \cm( \<Psi>,z)$ in \eqref{equa-z-re}, let us further investigate the properties of the sequence $(\Lambda^{(n)}(t))_{n\geq 1}$ defined by \eqref{defi-lambda-n}, and prove that this sequence actually converges to some process $\Lambda$ in a space of (regular enough) functions. To this end, we recall that according to \eqref{defi:approx-b-n}, the approximation $\dot{B}^{(n)}$ can be represented as $\dot{B}^{(n)}:=\partial_t\partial_x B^{(n)}$, where $B^{(n)}$ - given by~\eqref{defi:approx-sheet} - converges to a fractional sheet~$B$ on $\R^2$.

\begin{proposition}\label{prop:conv-lamb-n}
For all $\frac12<\ga<H_0$ and $T>0$, one has almost surely
\begin{equation}\label{conver-la-n}
\Lambda^{(n)} \to \Lambda \quad \text{in} \ H^\ga(\R),
\end{equation}
where $\Lambda:\Omega \times \R \to  \C$ is the process given by
\begin{equation}\label{defi-lambdada}
\Lambda(t):=\frac{\imath }{2\pi} \chi(t) \bigg[ \chi(t) B(t,2\pi)  -\int_0^t ds \, \chi'(s) B(s,2\pi) \bigg].
\end{equation}
\end{proposition}

\begin{proof}
Based on the expression \eqref{defi:approx-sheet}, it is readily checked that $B^{(n)}(0,x)=B^{(n)}(t,0)=\partial_t B^{(n)}(t,0)=0$ a.s., and so, for every $t\in \R$,
\begin{align*}
\Lambda^{(n)}(t)=\int_{\T} \overline{\<Psi>^{(n)}(t,x)} \, dx&=\imath \chi(t)\int_0^t ds \, \chi(s) \int_{ \T} dy \,  \partial_x\partial_t B^{(n)}(s,y) \\
&=\frac{\imath }{2\pi} \chi(t)\int_0^t ds \, \chi(s) \partial_t B^{(n)}(s,2\pi) \\
&=\frac{\imath }{2\pi} \chi(t) \bigg[ \chi(t) B^{(n)}(t,2\pi)  -\int_0^t ds \, \chi'(s) B^{(n)}(s,2\pi) \bigg].
\end{align*}
Now, it is a well-known fact $B^{(n)}(.,2\pi) \to B(.,2\pi)$ in the classical H{\"o}lder space $\cac^\ga([-T,T])$ (for every $T>0$), which, by using the standard Sobolev spaces identifications (see e.g. \cite[Section 1.2]{runst-sickel}), immediately entails the desired convergence $\Lambda^{(n)} \to \Lambda$ in $H^{\ga}(\R)$.

\end{proof}

\begin{remark}
The results of Proposition \ref{prop:ice-cream} (item $(ii)$) and Proposition \ref{prop:conv-lamb-n} point out the fundamental difference between the roles of $\si^{(n)}$ and $\Lambda^{(n)}$ in the renormalization procedure. 

Going back to the definition of $\<IPsi2>^{(n)}$ in \eqref{defi-i-psi-2-n} and given the convergence result in \eqref{conver-la-n}, it becomes indeed clear that the sequence $(\si^{(n)})_{n\geq 1}$ is really the one responsible for the final convergence of $(\<IPsi2>^{(n)})_{n\geq 1}$, \textit{at least as a general distribution}. On the other hand, subtracting $\ci_\chi\big(\Lambda^{(n)} \cdot \<Psi>^{(n)}\big)$ allows us to guarantee that the latter convergence actually takes place in $Z^{s,b}$ (with $s>0$ and $b>\frac12$), a property of paramount importance in our nonlinear setting. 
\end{remark}

\begin{remark}
It is interesting to note that for every fixed $t\geq 0$, the limit $\Lambda$ in \eqref{defi-lambdada} only depends on the behavior of the noise on the boundary, that is in $2\pi$. 
\end{remark}

\subsection{Stochastic product as a random operator}

Let us now turn to the interpretation and control of the two product terms 
\begin{equation}\label{cl-phi-plus-minus}
\cl^{\<Psi>,+}_\chi( z):= \ci_\chi \cm(z,\<Psi>) \quad \text{and} \quad \cl^{\<Psi>,-}_\chi(z):= \ci_\chi \cm( \<Psi>,z)\, ,
\end{equation} 
involved in the remainder equation \eqref{equa-z-re} and in equation \eqref{equa-z-reno} below. 

\smallskip

As reported in Section \ref{subsec:reform-main-steps}, the negative Sobolev regularity of $\<Psi>$ (see Proposition \ref{prop:conv-psi-n}), combined with the lack of (deterministic) regularization effect of the Schr{\"o}dinger operator, rules out the possibility of a direct pathwise interpretation, and invites us to consider these two objects as \textit{random operators} acting on $z$.

\smallskip

Just as above, our strategy to implement this idea will be based on an approximation procedure. For the sake of consistency, we shall again rely on the smooth approximation $\dot{B}^{(n)}$ of~$\dot{B}$ introduced in \eqref{defi:approx-b-n}, as well as on the associated smooth process $\<Psi>^{(n)}$ defined through~\eqref{defi-luxo-n}. Thus, for every fixed $n\geq 1$, we introduce the approximated operators $\cl^{(n),+}$ and $\cl^{(n),-}$ along the formulas
\begin{equation}\label{defi-cl-n-pm}
\cl^{(n),+}(z)=\cl^{\<Psi>^{(n)},+}_\chi(z):= \ci_\chi \cm( z,\<Psi>^{(n)}), \quad \cl^{(n),-}(z)=\cl^{\<Psi>^{(n)},-}_\chi(z):= \ci_\chi \cm( \<Psi>^{(n)}, z),
\end{equation}
which, in this smoothened situation, are easy to interpret for any function $z\in L^2(\R_+\times \mathbb{T})$.

\smallskip

Our objective in this setting will thus be to prove the convergence of $(\cl^{(n),+})_{n\geq 1}$ and $(\cl^{(n),-})_{n\geq 1}$ as sequences of random operators acting on a specific space $Z^{s,b}$ (as defined in Notation \ref{notation:x-b-c}). For a clear and efficient statement of this property, we need to introduce the following additional notation.

\begin{notation}
For all $s,b\geq 0$ and $\mu\in [0,1]$, we denote by $\mathfrak{L}_\mu(Z^{s,b},Z^{s,b})$ the space of linear operators $\cl$ from $Z^{s,b}$ to $Z^{s,b}$ such that
\begin{equation}\label{l-mu-x-b-c}
\big\|\cl\big\|_{\mathfrak{L}_\mu(Z^{s,b},Z^{s,b})}:=\sup_{0<\tau\leq 1}\sup_{z\in Z^{s,b},z\neq 0} \frac{1}{\tau^\mu \|z\|_{Z^{s,b}}} \big\| \chi_\tau \cdot \cl(z)\big\|_{Z^{s,b}} <\infty.
\end{equation}
We denote by $\overline{\mathfrak{L}}_\mu(Z^{s,b},Z^{s,b})$ the space of antilinear operators from $Z^{s,b}$ to $Z^{s,b}$ satisfying the same condition.
\end{notation}

Our main result regarding the product elements in \eqref{cl-phi-plus-minus}-\eqref{defi-cl-n-pm} now reads as follows.

\begin{proposition}\label{prop:control-op-l-intro}
Let $\frac12 < H_0<1$ and  $0< H_1<1$ satisfy  $\frac32 <2H_0+H_1< 2$, and set $\underline{H}:=2H_0+H_1-1\in (\frac12,1)$. Fix $s\geq 0$ such that
\begin{equation}
1-\underline{H}<s<\frac12.
\end{equation}
Then there exists $b_s^\ast\in (\frac12,1)$ such that for every $b\in (\frac12,b_s^\ast]$, the sequence $(\cl^{(n),+})_{n\geq 1}$, resp. $(\cl^{(n),-})_{n\geq 1}$, converges almost surely in $\mathfrak{L}_\mu(Z^{s,b},Z^{s,b})$, resp. $\overline{\mathfrak{L}}_\mu(Z^{s,b},Z^{s,b})$, for some $\mu\in (0,1)$. 

\smallskip

We denote the limit by $\cl^{\<Psi>,+}_\chi$, resp. $\cl^{\<Psi>,-}_\chi$.

\end{proposition}

This result, which will be proved in Section \ref{sec:estim-random-op} below, is thus another step toward the complete interpretation and control of the remainder equation \eqref{equa-z-re} (or \eqref{equa-z-reno} below): on top of $\<IPsi2>$, we are endowed with two natural extensions $\cl^{\<Psi>,+}_\chi(z)$ and $\cl^{\<Psi>,-}_\chi(z)$ of the elements in~\eqref{defi-cl-n-pm} such that for all $z\in Z^{s,b}$ and $\tau\in (0,1]$, 
$$\big\|\chi_\tau \cdot \cl^{\<Psi>,\pm}_\chi(z)\big\|_{Z^{s,b}} \lesssim \tau^\mu \big\|z\big\|_{Z^{s,b}}, $$
for some $s\in (0,\frac12)$, $b\in (\frac12,1)$ and $\mu\in (0,1)$.

\begin{remark}

Similar constructions of \textit{random operators} can be found in the dispersive literature: see e.g. \cite[Lemma 3.5]{OWZ} for a cubic wave model with additive noise, or \cite[Proposition 3.2]{liu} for a quadratic Schr{\"o}dinger equation with rough random initial condition. This procedure is also one of the ingredients toward the recent breakthrough results by Deng, Nahmod and Yue about random Hamiltonian NLS (see e.g. the involvement of the operator norm $\big\|\cl^{\zeta}\big\|_{X^{1-b,b}(k\to k')}$ in the statement of \cite[Proposition 5.1]{DNY}).\\
\indent We are not aware of any previous similar construction for a Schr{\"o}dinger model with additive noise. Note also that the above references \cite{DNY,liu,OWZ} all involve white noises, and so the technical random-tensor estimates therein established (see e.g. \cite[Lemma C.3]{OWZ}) do not apply to our fractional setting.

\end{remark}

\subsection{Wellposedness of the (renormalized) remainder equation}

Recall that we are interested in the (renormalized) remainder equation
\begin{equation}\label{equa-z-reno}
z=\chi_\tau \cdot \ci_\chi \cm(z,z)+\Big[\chi_{\tau} \cdot \cl^{\<Psi>,+}_\chi( z)+\chi_{\tau}\cdot \cl^{\<Psi>,-}_\chi( z)\Big]-\chi_\tau\cdot \ci_\chi\big(\Lambda \cdot z\big)+\chi_{\tau}\cdot \<IPsi2>,
\end{equation}
with $  t\in \R, \ x\in \T$, where $\cl^{\<Psi>,\pm}_\chi( z)$ is interpreted through Proposition \ref{prop:control-op-l-intro} and $\<IPsi2>$ is interpreted through Proposition \ref{prop:ice-cream}.

\smallskip

Thanks to the (almost sure) estimates contained in these propositions, we essentially have all the tools in hand to set up a \textit{fully deterministic} fixed-point argument for \eqref{equa-z-reno}. The only additional ingredient we need is a suitable control on $\ci_\chi \cm(z,z)$ and $\ci_\chi\big(\Lambda \cdot z\big)$ in $Z^{s,b}$, which will be established in Section \ref{sec:deter-con} using deterministic estimates only.

\smallskip

For more clarity in the subsequent statement and its proof, let us denote 
\begin{multline}\label{topo-tree}
\big\|\Lambda,\cl^{\<Psi>,+}_\chi,\cl^{\<Psi>,-}_\chi, \<IPsi2>\big\|_{\ga,(b,s),\mu}\\
:=\big\|\Lambda\big\|_{H^\ga(\R)}+\big\|\cl^{\<Psi>,+}_\chi\big\|_{\mathfrak{L}_\mu(Z^{s,b},Z^{s,b})}+\big\|\cl^{\<Psi>,-}_\chi\big\|_{\overline{\mathfrak{L}}_\mu(Z^{s,b},Z^{s,b})}+\big\|\<IPsi2>\big\|_{Z^{s,b}},
\end{multline}
for $\ga,b,s,\mu$ to be specified later on.

\begin{theorem}\label{theo:well-posed-z}
Let $\frac12 < H_0<1$ and  $\frac14< H_1<1$ satisfy  $\frac74 <2H_0+H_1< 2$, and set $\underline{H}:=2H_0+H_1-1\in (\frac34,1)$.  Fix $s > 0$ such that 
$$1-\underline{H}<s< \underline{H}-\frac12$$
and with the notation of Proposition \ref{prop:ice-cream} and Proposition \ref{prop:control-op-l-intro}, set $b:=b_s \wedge b_s^\ast >\frac12 $. Also, fix $\mu:=\mu_{s,b}\in (0,1)$ such that the result of Proposition \ref{prop:control-op-l-intro} applies almost surely,  that is
$$\cl^{(n),+} \to \cl^{\<Psi>,+}_\chi \   \text{in} \ \mathfrak{L}_\mu(Z^{s,b},Z^{s,b}) \quad \text{and} \quad \cl^{(n),-} \to \cl^{\<Psi>,-}_\chi \   \text{in} \ \overline{\mathfrak{L}}_\mu(Z^{s,b},Z^{s,b}) .$$

Then the following assertions hold true.

\smallskip

\noindent
$(i)$ Almost surely, there exists a time $\tau_\infty >0$ such that for every $\tau \in (0,\tau_\infty]$, the equation~\eqref{equa-z-reno} admits a unique solution $z^{(\tau)}\in Z^{s,b}$. 

\smallskip

\noindent
$(ii)$ Almost surely, there exists a sequence $(\tau^{(n)})_{n\geq 1}$ of positive times such that: 

\smallskip

$(ii$-$a)$ $\tau^{(n)} \to \tau_\infty$ as $n\to\infty$; 

\smallskip

$(ii$-$b)$ for all $n\geq 1$ and $\tau\in (0,\tau^{(n)}]$, the equation
\small
\begin{multline*} 
z^{(\tau,n)}=\chi_\tau \cdot\ci_\chi \cm(z^{(\tau,n)},z^{(\tau,n)})+\chi_\tau \cdot\cl^{(n),+}_\chi(z^{(\tau,n)})+\\
\hspace{6cm}+\chi_\tau \cdot\cl^{(n),-}_\chi(z^{(\tau,n)})-\chi_\tau\cdot \ci_\chi\big(\Lambda^{(n)} \cdot z^{(\tau,n)}\big)+\chi_\tau \cdot\<IPsi2>^{(n)}
\end{multline*}
\normalsize
admits a unique solution $z^{(\tau,n)}\in Z^{s,b}$.

\smallskip

\noindent
$(iii)$ There exists a time $\tau_{\<circle>}\leq \min\big(\inf_{n\geq 1} \tau^{(n)},\tau_\infty\big)$ such that  $\tau_{\<circle>}>0$ almost surely and for every $\tau\in (0,\tau_{\<circle>}]$, the sequence $(z^{(\tau,n)})_{n\geq 1}$ converges a.s. to $z^{(\tau)}$ in $Z^{s,b}$.

\end{theorem}

\begin{remark}
Although we have insisted on the dependence on $\tau$ through the notation~$z^{(\tau)}$, it should be noted that the following consistency property holds true (as a consequence of the uniqueness of the solution): if $0< \tau_1 < \tau_2 <\tau_\infty$, then $z^{(\tau_1)} \equiv z^{(\tau_2)}$ on the time interval $[-\tau_1, \tau_1]$.
\end{remark}

\begin{proof}[Proof of Theorem \ref{theo:well-posed-z}]

$(i)$-$(ii)$. For all $0<\tau\leq 1$ and $z\in Z^{s,b}$, set
$$
\gga_\tau(z):=\chi_\tau \cdot \ci_\chi \cm(z,z)+\chi_{\tau} \cdot \cl^{\<Psi>,+}_\chi( z)+\chi_{\tau}\cdot \cl^{\<Psi>,-}_\chi( z)-\chi_\tau\cdot \ci_\chi\big(\Lambda \cdot z\big)+\chi_{\tau}\cdot \<IPsi2>.
$$
By gathering the results of Proposition \ref{prop:ice-cream}, Proposition \ref{prop:control-op-l-intro}, Proposition \ref{prop:lam-z} and Proposition~\ref{prop:control-m-z-z}, we can guarantee that (a.s.) $\gga_\tau$ is a well-defined map from $Z^{s,b}$ to $Z^{s,b}$. Besides, there exist constants $\nu>0$ and $C>0$ such that for all $z,z^{(1)},z^{(2)}\in Z^{s,b}$,
\begin{equation}\label{control-gamma-1}
\big\| \gga_\tau(z)\big\|_{Z^{s,b}}\leq C \Big\{ \tau^\nu \|z\|_{Z^{s,b}}^2+\tau^\nu N_{\<circle>}\, \|z\|_{Z^{s,b}}+N_{\<circle>}\Big\}
\end{equation}
and
\small
\begin{align}
&\big\| \gga_\tau(z^{(1)})-\gga_\tau(z^{(2)})\big\|_{Z^{s,b}}\nonumber\\
&\leq C\Big\{\tau^\nu \big(\|z^{(1)}\|_{Z^{s,b}}+\|z^{(2)}\|_{Z^{s,b}}\big) \|z^{(1)}-z^{(2)}\|_{Z^{s,b}}+\tau^\nu N_{\<circle>}\, \|z^{(1)}-z^{(2)}\|_{Z^{s,b}}\Big\}\, ,\label{control-gamma-2}
\end{align}
\normalsize
where we have used the shorthand notation $N_{\<circle>}:=\big\|\Lambda,\cl^{\<Psi>,+}_\chi,\cl^{\<Psi>,-}_\chi, \<IPsi2>\big\|_{\ga,(b,s),\mu}$. 

\smallskip

In a similar way, by setting for every $n\geq 1$
\small
$$
\gga^{(n)}_\tau(z):=\chi_\tau \cdot \ci_\chi \cm(z,z)+\chi_{\tau} \cdot \cl^{(n),+}_\chi( z)+\chi_{\tau}\cdot \cl^{(n),-}_\chi( z)-\chi_\tau\cdot \ci_\chi\big(\Lambda^{(n)} \cdot z\big)+\chi_{\tau}\cdot \<IPsi2>^{(n)},
$$
\normalsize
we get that for the same constants $\nu,C>0$ as above,
\begin{equation}\label{control-gamma-1-n}
\big\| \gga^{(n)}_\tau(z)\big\|_{Z^{s,b}}\leq C \Big\{\tau^\nu \|z\|_{Z^{s,b}}^2+\tau^\nu N^{(n)}_{\<circle>}\, \|z\|_{Z^{s,b}}+N^{(n)}_{\<circle>}\Big\}
\end{equation}
and
\begin{align}
&\big\| \gga^{(n)}_\tau(z^{(1)})-\gga^{(n)}_\tau(z^{(2)})\big\|_{Z^{s,b}}\nonumber\\
&\leq C \Big\{\tau^\nu \big(\|z^{(1)}\|_{Z^{s,b}}+\|z^{(2)}\|_{Z^{s,b}}\big) \|z^{(1)}-z^{(2)}\|_{Z^{s,b}}+\tau^\nu N^{(n)}_{\<circle>}\, \|z^{(1)}-z^{(2)}\|_{Z^{s,b}}\Big\}\, ,\label{control-gamma-2-n}
\end{align}
where we have naturally set this time $N^{(n)}_{\<circle>}:=\big\|\Lambda^{(n)},\cl^{(n),+}_\chi,\cl^{(n),-}_\chi, \<IPsi2>^{(n)}\big\|_{\ga,(b,s),\mu}$. \\
\indent Using \eqref{control-gamma-1}-\eqref{control-gamma-2}-\eqref{control-gamma-1-n}-\eqref{control-gamma-2-n} and the fact that $N^{(n)}_{\<circle>} \to N_{\<circle>}$, both assertions $(i)$ and $(ii)$ can be derived from elementary deterministic fixed-point and control arguments (see e.g. the proof of \cite[Theorem 1.6]{deya-wave} for details in a similar situation).  Besides, it holds that
$$K_{\<circle>}:=\sup_{\tau\in (0,\min(\inf_{n\geq 1} \tau^{(n)},\tau_\infty)]}\sup_{n\geq 1}\Big(\|z^{(\tau,n)}\|_{Z^{s,b}}\vee \|z^{(\tau)}\|_{Z^{s,b}}\Big) \ < \ \infty \quad \text{almost surely}. $$
 
\smallskip

\noindent
$(iii)$ Thanks to the results of Proposition \ref{prop:ice-cream}, Proposition \ref{prop:conv-lamb-n} and Proposition \ref{prop:control-op-l-intro}, we know that a.s.
$$
M_{\<circle>}^{(n)}:=\big\|\Lambda^{(n)}-\Lambda,\cl^{(n),+}_\chi-\cl^{\<Psi>,+}_\chi,\cl^{(n),-}_\chi-\cl^{\<Psi>,-}_\chi, \<IPsi2>^{(n)}-\<IPsi2>\big\|_{\ga,(b,s),\mu} \to 0 \quad  \text{as} \ n\to \infty.
$$
For every $0<\tau\leq \min\big(\inf_{n\geq 1} \tau^{(n)},\tau_\infty\big)$, both solutions $z^{(\tau,n)}$ and $z^{(\tau)}$ are well defined in~$Z^{s,b}$, and we can thus examine the equation satisfied $z^{(\tau,n)}-z^{(\tau)}$. In fact, using again the estimates of Proposition~\ref{prop:ice-cream}, Proposition~\ref{prop:control-op-l-intro}, Proposition~\ref{prop:lam-z} and Proposition~\ref{prop:control-m-z-z}, it is easy to check the existence of constants $C,\mu>0$ such that for all $0<\tau\leq \min\big(\inf_{n\geq 1} \tau^{(n)},\tau_\infty\big)$ and $n\geq 1$,
\small
\begin{align*}
&\big\|z^{(\tau,n)}-z^{(\tau)}\big\|_{Z^{s,b}}\\
&\leq C\Big\{ \tau^\nu \big(\|z^{(\tau,n)}\|_{Z^{s,b}}+\|z^{(\tau)}\|_{Z^{s,b}}\big) \|z^{(\tau,n)}-z^{(\tau)}\|_{Z^{s,b}}+\tau^\nu M_{\<circle>}^{(n)}\, \|z^{(\tau,n)}\|_{Z^{s,b}}\\
&\hspace{7cm}+\tau^\nu N_{\<circle>}\, \|z^{(\tau,n)}-z^{(\tau)}\|_{Z^{s,b}}+M_{\<circle>}^{(n)}\Big\}\\
&\leq C \Big\{\tau^\nu \big(2K_{\<circle>}+N_{\<circle>}\big) \|z^{(\tau,n)}-z^{(\tau)}\|_{Z^{s,b}}+\tau^\nu K_{\<circle>}\,  M_{\<circle>}^{(n)}+M_{\<circle>}^{(n)}\Big\}.
\end{align*}
\normalsize
It now becomes clear that we can find $\tau_{\<circle>}\in (0,\min\big(\inf_{n\geq 1} \tau^{(n)},\tau_\infty\big)]$ small enough so that for all $0<\tau\leq \tau_{\<circle>}$ and $n\geq 1$,
\begin{align*}
&\big\|z^{(\tau,n)}-z^{(\tau)}\big\|_{Z^{s,b}}\lesssim M_{\<circle>}^{(n)}.
\end{align*}
Since $M_{\<circle>}^{(n)}\to 0$ a.s., the conclusion immediately follows.

\end{proof}

\subsection{Convergence of the (renormalized) strong solutions}

To end up with this study, we propose to turn the wellposedness properties of Theorem \ref{theo:well-posed-z} into an approximation result for the original equation \eqref{starting-equation}, and thus provide a final statement following the pattern of those in \eqref{cubic-nls-sequence}, \eqref{example-heat} or \eqref{example-wave}.

\begin{theorem}\label{theo:seq}
Assume that $\frac12<H_0<1$ and $\frac14<H_1<1$ satisfy $\frac74<2H_0+H_1<2$. Then the statement of Theorem~\ref{main-obj} holds true with $\Lambda^{(n)},\si^{(n)}$ given by \eqref{defi-lambda-n} and \eqref{defi-si-n}, i.e.
\small
\begin{equation}\label{lambda-si-n}
\Lambda^{(n)}(t):=\int_{\T} \overline{\<Psi>^{(n)}(t,x)} \, dx \quad \text{and} \quad  \si^{(n)}(t,x):=\mathbb{E}\bigg[e^{-\imath t\Delta }\bigg(\cm( \<Psi>^{(n)},\<Psi>^{(n)})(t)-  \<Psi>^{(n)}(t) \int_{\T} \overline{\<Psi>^{(n)}}(t)\bigg)(x)\bigg].
\end{equation}
\normalsize
In other words, there exists a (random) time $T_0>0$ such that the sequence $(u^{(n)})_{n\geq 1}$ of classical solutions to the renormalized equation
\begin{equation}\label{renorm-equation}
(\imath \partial_t-\Delta)u^{(n)}=\Big[\overline{u^{(n)}}-\Lambda^{(n)}\Big] u^{(n)}-\si^{(n)}+\dot{B}^{(n)}
\end{equation}
is well defined and converges a.s. in $\mathcal{C}\big([-T_0,T_0];H^{-\al}(\T)\big)$, for every $\al>2- (2H_0+H_1)$.

\smallskip

Moreover, denoting the limit by $u$, it holds that 
\begin{equation}\label{decompo-solu-u}
u-e^{-\imath . \Delta }\<Psi>\, \in\,   \mathcal{C}\big([-T_0, T_0]; H^{s}(\T)\big),
\end{equation}
for every $0<s<2H_0+H_1-\frac32$.
\end{theorem}

Observe that the fundamental decomposition of the solution $u$ in \eqref{decompo-solu-u} is an immediate consequence of the Da Prato-Debussche-trick used in the study.

\begin{proof}[Proof of Theorem \ref{theo:seq}]

Using the setting and notation of Theorem \ref{theo:well-posed-z}$(iii)$, we fix $\tau:=\tau_{\<circle>}\in (0,1]$ and consider the sequence $(z^{(\tau,n)})_{n\geq 1}$ of solutions in $Z^{s,b}$ to the equation
\small
\begin{equation}\label{equa-z-reno-approx-bis}
\begin{split}
&z^{(\tau,n)}=\chi_\tau \cdot\ci_\chi \cm(z^{(\tau,n)},z^{(\tau,n)})+\chi_\tau \cdot\cl^{(n),+}_\chi(z^{(\tau,n)})+\chi_\tau \cdot\cl^{(n),-}_\chi(z^{(\tau,n)})\\
&\hspace{6cm}-\chi_\tau\cdot \ci_\chi\big(\Lambda^{(n)} \cdot z^{(\tau,n)}\big)+\chi_\tau \cdot\<IPsi2>^{(n)}.
\end{split}
\end{equation}
\normalsize
Let us now quickly trace back the steps that led us (in Section \ref{subsec:reform-main-steps}) to this remainder equation. 

\smallskip

First, combining \eqref{hat-cm}, \eqref{defi:ipsi2-main} and \eqref{defi-cl-n-pm}, we can rephrase \eqref{equa-z-reno-approx-bis} as
\small
\begin{align*}
&z^{(\tau,n)}_k=\chi_\tau \cdot\ci_\chi \cm\big(z^{(\tau,n)},z^{(\tau,n)}\big)_k+\chi_\tau \cdot \ci_\chi \cm(z^{(\tau,n)},\<Psi>^{(n)})_{ k}+\chi_\tau \cdot \ci_\chi \cm(\<Psi>^{(n)},z^{(\tau,n)})_{k}\\
&\hspace{1cm}-\chi_\tau\cdot \ci_\chi\big(\Lambda^{(n)} \cdot z^{(\tau,n)}\big)+\chi_\tau\cdot \Big[\ci_\chi\big( \cm\big(\<Psi>^{(n)},\<Psi>^{(n)}\big)_k-\ci_\chi\big(\Lambda^{(n)} \cdot \<Psi>^{(n)}_k\big)-\ci_\chi\big( e^{-\imath |k|^2 .}\si^{(n)}_k \big)\Big],
\end{align*}
\normalsize
or in other words
\small
\begin{equation*}
\begin{split}
&z^{(\tau,n)}_k=\chi_\tau \cdot\ci_\chi \cm\big(z^{(\tau,n)}+\<Psi>^{(n)},z^{(\tau,n)}+\<Psi>^{(n)}\big)_k\\
&\hspace{4cm}-\chi_\tau\cdot \ci_\chi\big(\Lambda^{(n)} \cdot (z_k^{(\tau,n)}+ \<Psi>^{(n)}_k)\big)-\chi_\tau\cdot \ci_\chi\big( e^{-\imath |k|^2 .}\si^{(n)}_k \big).
\end{split}
\end{equation*}
\normalsize
Then, setting
\begin{equation*}
u^{(\tau,n)}_k(t):=e^{\imath |k|^2t}\big[z^{(\tau,n)}_k(t)+\<Psi>^{(n)}_k(t)\big],
\end{equation*}
we get that for every $t\in \R$,
\small
\begin{align*}
&e^{-\imath |k|^2t}u^{(\tau,n)}_k(t)-\<Psi>^{(n)}_k(t)\\
&=\chi_\tau(t) \cdot\chi(t)\int_0^t ds \, \chi(s) \sum_{k_1\in \Z} e^{\imath s\Omega_{k,k_1}} e^{-\imath s|k+k_1|^2}u^{(\tau,n)}_{k+k_1}(s) \overline{ e^{-\imath s |k_1|^2}u^{(\tau,n)}_{k_1}(s)}\\
&\hspace{0.5cm}-\chi_\tau(t) \cdot\chi(t)\bigg[\int_0^t ds \, \chi(s)e^{-\imath s|k|^2}\Lambda^{(n)}(s) u^{(\tau,n)}_k(s)+\int_0^t ds \, \chi(s)e^{-\imath s|k|^2 }\si^{(n)}_k(s)\bigg]\\
&=\chi_\tau(t) \cdot\chi(t)\int_0^t ds \, \chi(s)e^{-\imath s |k|^2 } \sum_{k_1\in \Z}  u^{(\tau,n)}_{k+k_1}(s) \overline{ u^{(\tau,n)}_{k_1}(s)}\\
&\hspace{0.5cm}-\chi_\tau(t) \cdot\chi(t)\bigg[\int_0^t ds \, \chi(s)e^{-\imath s |k|^2}\Lambda^{(n)}(s) u^{(\tau,n)}_k(s)+\int_0^t ds \, \chi(s) e^{-\imath s |k|^2 }\si^{(n)}_k(s)\bigg],
\end{align*}
\normalsize
that is, with expression \eqref{defi-luxo-n} of $\<Psi>^{(n)}$ in mind,
\small
\begin{align*}
u^{(\tau,n)}_k(t)&= \chi(t)\int_0^t ds \, \chi(s)e^{\imath (t-s) |k|^2} \dot{B}^{(n)}_k(s)+\chi_\tau(t) \cdot\chi(t)\int_0^t ds \, \chi(s)e^{\imath (t-s) |k|^2}\big(|u^{(\tau,n)}(s)|^2\big)_k \\
&-\chi_\tau(t) \cdot\chi(t)\bigg[\int_0^t ds \, \chi(s)e^{\imath (t-s) |k|^2}\Lambda^{(n)}(s) u^{(\tau,n)}_k(s)+\int_0^t ds \,\chi(s) e^{\imath (t-s) |k|^2}\si^{(n)}_k(s)\bigg].
\end{align*}
\normalsize
In particular, for every $t\in [-\tau,\tau]$, it holds that
\small
\begin{equation*}
\begin{split}
&u^{(\tau,n)}_k(t)=\int_0^t ds \, e^{\imath (t-s)|k|^2} \dot{B}^{(n)}_k(s)\\
&\hspace{3cm}+\int_0^t ds \, e^{\imath (t-s)|k|^2}\Big[\big(|u^{(\tau,n)}(s)|^2\big)_k-\Lambda^{(n)}(s) u^{(\tau,n)}_k(s) - \si^{(n)}_k(s)\Big],
\end{split}
\end{equation*}
\normalsize
which is nothing but the mild formulation of \eqref{renorm-equation} on $[-\tau,\tau]$. By uniqueness of the (mild) solution for the regularized problem, we can thus assert that for every $t\in [-\tau,\tau]$, $u^{(\tau,n)}(t,.)\equiv u^{(n)}(t,.)$, or otherwise stated: for all $t\in [-\tau,\tau]$ and $k\in \Z$,
\begin{equation*} 
u^{(n)}_k(t)=e^{\imath |k|^2t}\big[z^{(\tau,n)}_k(t)+\<Psi>^{(n)}_k(t)\big],
\end{equation*}
namely

\begin{equation}\label{repres-til-u}
u^{(n)}(t)=e^{-\imath t \Delta}\big[z^{(\tau,n)}(t)+\<Psi>^{(n)}(t)\big].
\end{equation}
We know from Proposition \ref{prop:conv-psi-n} that for $\al>2- 2H_0-H_1$, the sequence $e^{-\imath . \Delta} \<Psi>^{(n)}$ converges (almost surely) to $ e^{-\imath . \Delta}\<Psi>$ in $\mathcal{C}\big([-\tau,\tau];H^{-\al}(\T)\big)$.  Besides, by  Theorem \ref{theo:well-posed-z}$(iii)$,
${\dis  \big\|z^{(\tau,n)}-z^{(\tau)}\big\|_{Z^{s,b}}  \stackrel{n\to\infty}{\longrightarrow}0}$, thus by \eqref{norm:x-b-c}
$$\dis  \big\|e^{-\imath t \Delta} \big(z^{(\tau,n)}-z^{(\tau)}  \big)\big\|_{X^{s,b}}  \stackrel{n\to\infty}{\longrightarrow}0.$$

Now use that $X^{s,b} \subset  \mathcal{C}\big([-\tau,\tau];H^{s}(\T)\big) \subset \mathcal{C}\big([-\tau,\tau];H^{-\al}(\T)\big)$, for $b>1/2$, then from~\eqref{repres-til-u}, we deduce  the (almost sure) convergence of $u^{(n)}$ in $\mathcal{C}\big([-\tau,\tau];H^{-\al}(\T)\big)$, as well as the decomposition \eqref{decompo-solu-u} of $u$.
\end{proof}

We have thus reached our objective: in contrast with the rescaled models \eqref{sto-cubic-nls-renorm} and \eqref{liu-model}, the renormalizing sequences $(\Lambda^{(n)})_{n\geq 1},(\si^{(n)})_{n\geq 1}$ in \eqref{lambda-si-n} are explicit, i.e. they are explicitly defined in terms of the approximated noise $\dot{B}^{(n)}$, which, as far as is new for a stochastic NLS equation.

\smallskip

Let us now complement the above results with a series of remarks.

\begin{remark}
In the present study, we have preferred to highlight the impact of an additive noise and thus chosen to work with a zero initial condition in \eqref{starting-equation}. However, the consideration of a regular enough initial condition $f$ on $\mathbb{T}$ (for instance, $f\in H^s(\mathbb{T})$, with $s>0$ as in Theorem \ref{theo:well-posed-z}) could be easily included into our analysis.
\end{remark}

\begin{remark}\label{rk:non-trivial-renorm}
As explained in Section \ref{subsec:main-res-sto}, the choice of the sequences $(\Lambda^{(n)})_{n},(\si^{(n)})_{n}$ in \eqref{lambda-si-n} follows from our renormalization of $\ci_\chi \cm( \<Psi>^{(n)},\<Psi>^{(n)})$, that is from our definition of $\<IPsi2>^{(n)}$ in \eqref{hat-cm}-\eqref{defi:ipsi2-main}. 
\end{remark}

\begin{remark}

Based on the proofs of Theorems \ref{theo:well-posed-z} and \ref{theo:seq}, it should be clear to the reader that the limit-solution $u$ only depends on the noise approximation through the $4$-uplet $\big(\Lambda,\cl^{\<Psi>,+}_\chi,\cl^{\<Psi>,-}_\chi, \<IPsi2>\big)$. In other words, any alternative approximation $\dot{\tilde{B}}^{(n)}$ of $\dot{B}$ such that (with obvious notation)
$$\Big(\tilde{\Lambda}^{(n)},\cl^{\tilde{\<Psi>}^{(n)},+}_\chi,\cl^{\tilde{\<Psi>}^{(n)},-}_\chi, \tilde{\<IPsi2>}^{(n)}\Big) \stackrel{n\to\infty}{\longrightarrow} \big(\Lambda,\cl^{\<Psi>,+}_\chi,\cl^{\<Psi>,-}_\chi, \<IPsi2>\big)$$
for the topology in \eqref{topo-tree}, would ultimately lead to the same limit-solution $u$. 

\end{remark}

\begin{remark}

Instead of the Fourier-type approximation $\dot{B}^{(n)}$ in \eqref{defi:approx-b-n}, we could have chosen any mollifying approximation $\dot{B}_\rho^{(n)}:=\rho_n \ast_{t,x} \dot{B}$ of $\dot{B}$, where $\rho_n(s,x):=2^{3n}\rho(2^{2n} s,2^n x)$ for some test-function $\rho: \R^{2} \to \R_+$ such that $\int_{\R^{2}} \rho(s,x) \, ds dx=1$. In this case, it can be shown that the associated $4$-uplet $\big(\Lambda,\cl^{\<Psi>,+}_\chi,\cl^{\<Psi>,-}_\chi, \<IPsi2>\big)$, and accordingly the limit-solution $u$ itself, do not depend on the choice of the mollifying function $\rho$. 

\end{remark}

\begin{remark}
It turns out that if one replaces the nonlinearity $|u|^2$ in \eqref{starting-equation} with either $u^2$ or $\bar{u}^2$, then the (local) treatment of the equation becomes essentially trivial, as a consequence of the following classical bilinear estimates established Kenig, Ponce and Vega in~\cite{KPV}: for all $s>-\frac12$ there exists $b>\frac12$ such  that
\begin{equation}\label{KPV}
\bigg\|\chi_\tau  \int_0^. \, dre^{-\imath (.-r)\Delta}\big(u v\big)\bigg\|_{X^{s,b}}\lesssim  \big\|u\big\|_{X^{s,b}} \big\|v\big\|_{X^{s,b}},
\end{equation}
with a similar control for $\bar{u}\bar{v}$ (recall that $X^{s,b}$ refers to the Bourgain space introduced in~\eqref{bourg}). Denoting by $\Psi$ the solution of the linear problem $(\imath \partial_t-\Delta)\Psi=\dot{B}$ with $\Psi(0)=0$, one can show with similar arguments to those of Proposition \ref{prop:conv-psi-n} that $\Psi\in X^{s,\frac12}$ a.s. for every $s>2H_0+H_1-2>-\frac12$ (provided $2H_0+H_1>\frac32$), and therefore the corresponding equation
\begin{equation}\label{additive-noise-u-sq}
(\imath \partial_t-\Delta)u=u^2+\dot{B}, \quad u_0=0,\quad \quad  t\in \R, \ x\in \mathbb{T},
\end{equation}
is guaranteed to be (locally) wellposed in the Bourgain scale.\\
\indent Let us now emphasize that for the nonlinearity $|u|^2$, the bilinear estimate \eqref{KPV} is only known to be true for $s\geq 0$ (see Proposition \ref{prop:control-m-z-z}), with counterexamples at negative regularity (see \cite[Theorem 1.10$(iii)$]{KPV} or \cite{kishimoto}). In the fractional situation, the restriction thus corresponds to the case where $2H_0+H_1>2$, and so the estimate does not cover the rough regime $\frac74<2H_0+H_1\leq 2$ under consideration in the present study.
\end{remark}

\begin{remark}

It should be noted that the case of a standard space-time white noise (i.e., $H_0=H_1=\frac12$) in \eqref{starting-equation} is essentially beyond the reach of any expansion method. Indeed, in this case, and according to \cite[Proposition 2.6]{deya-restric}, the construction of the fundamental second-order tree $\<IPsi2>$ turns out to be impossible.  \\
\indent At this point, it is not clear to us whether a higher-order expansion strategy (with construction of higher-order random operators) or the consideration of some \textit{random tensor ansatz} (along the terminology of \cite{DNY}) could allow to cover the index domain  $\frac32<2H_0+H_1\leq \frac74$, which is not treated in our analysis  (see \cite{deya-restric} for further results in the fractional-in-time white-in-space situation).

\end{remark}

\

The rest of the paper is organized as follows. In Section \ref{sec:preliminaries}, we introduce a few technical ingredients related to the model and useful for its analysis. Sections \ref{sec:deter-con}-\ref{sec:tree-elements}-\ref{sec:estim-random-op} are then devoted to the interpretation and the control of the successive terms involved in the central equation~\eqref{equa-z-reno}. To be more specific, we start with \textit{fully deterministic} considerations about the control of $\ci_\chi\big(\Lambda \cdot z\big)$ and $\ci_\chi \cm(z,w)$ for $z,w$ in $Z^{s,b}$ (Section~\ref{sec:deter-con}), then continue with \textit{stochastic} arguments towards the interpretation of $(\<Psi>,\<IPsi2>)$ (Section \ref{sec:tree-elements}), and finally conclude with a \textit{mixed deterministic/stochastic} analysis of the product terms $\ci_\chi \cm(z,\<Psi>)$, $\ci_\chi \cm(\<Psi>,z)$, seen as random operators (Section \ref{sec:estim-random-op}). In the Appendix~\ref{AppenA} we gather the proofs of some deterministic results.

\smallskip

\section{Preliminary estimates and fully-deterministic controls}\label{sec:preliminaries}

We gather a few technical properties that will play a central role in our analysis of the (renormalized) equation.

\smallskip

We recall that throughout the paper, the notation $\cf$ stands for the Fourier transform in time, whereas the notation $z_k$ refers to the convention introduced in Notation \ref{not:four-space}. Besides, we have fixed (once and for all) a time cut-off function $\chi:\R\to \R$, and the operator  
\begin{equation*} 
 \ci_\chi v(t)= - \imath \chi(t)\int_0^t ds \, \chi(s)v(s),
\end{equation*}
has been introduced in \eqref{defi:ci-chi}.

\subsection{Basic controls on time integration and time localization}

\begin{lemma}{\normalfont (\cite[Lemma 4.1]{DNY})}\label{lem:integr-kernel}
For all function $v:\R\to\R$ smooth enough, one can write
\begin{equation}\label{rel-def-i-chi}
\cf\big( \ci_{\chi} v\big)(\la)=\int_{\R} d\la_1 \, \ci_\chi(\la,\la_1) \cf(v)(\la_1),
\end{equation}
for some kernel $\ci_\chi$ satisfying
\begin{equation}\label{bou-ci-chi}
\big|\ci_\chi(\la,\la_1)\big|\lesssim \frac{1}{\langle \la\rangle \langle \la-\la_1\rangle}.
\end{equation}
\end{lemma}

We will also need the following elementary lemma related to time localization, i.e. multiplication with $\chi_\tau$ in \eqref{equa-z-reno} (see for instance \cite[Lemma 4.2]{DNY}).

\begin{lemma}\label{b-prim}
For all $\frac12 <b\leq b'<1$,  $s\in \R$ and $y:\R\times \T\to \R $ such that $y(0,.)=0$, it holds that
$$\big\| \chi_\tau\cdot y\big\|_{Z^{s,b}}\lesssim \tau^{b'-b}\big\| y\big\|_{Z^{s,b'}}.$$
\end{lemma}

\subsection{Singular integrals}

We first recall the following elementary estimate (see e.g. \cite[Lemma 4.2]{FOW} or \cite[Lemma 4.2]{GTV}).

\begin{lemma}\label{lem:tec-1d}
Let $\beta\geq \ga> 0$ be such that $\ga+\beta>1$. Then it holds that
$$\max\bigg(\sum_{k\in \mathbb{Z}}\frac{1}{\langle k -\xi\rangle^\ga \langle k-\zeta\rangle^\beta},\int_{\R} \frac{d\la}{\langle \la -\xi\rangle^\ga \langle \la-\zeta\rangle^\beta} \bigg)\lesssim \frac{1}{\langle \xi-\zeta\rangle^\al}$$
where $\al$ is explicitly given by
$$\al:=
\begin{cases}
\ga+\beta-1 & \text{if} \ \beta<1\\
\ga-\varepsilon & \text{if} \ \beta=1\\
\ga & \text{if}\ \beta>1
\end{cases},$$
for any $0<\varepsilon\leq \ga$.
\end{lemma}

The second technical lemma in this direction is more specific to our fractional setting, and is directly related to the covariance function of the linear solution to the problem (see the definition of $(\gga^{(H)},\Lambda^{(H)})$ in Proposition \ref{prop:estim-cov-1} below).

\begin{lemma}\label{lem:int-singu-zero}
Consider  parameters $\mu, \nu\in [0,1)$. 

\smallskip

\noindent
$(i)$ For all $(a,b)\in \R^2$ such that $|a|\leq |b|$, and for every $0<\varepsilon<1$, it holds that
\begin{equation}\label{adapt-lem-3}
\int_{\R}\frac{d\xi}{|\xi|^\nu} \frac{1}{\langle\xi-a\rangle}\frac{1}{\langle \xi-b\rangle} \lesssim \frac{1}{\langle a\rangle^\nu}\frac{1}{\langle b-a\rangle^{1-\varepsilon}}.
\end{equation}

\smallskip

\noindent
$(ii)$ For all $p\geq 1$, $|a|\geq 1$ and $0<\varepsilon<\nu$, one has
\begin{equation}\label{adapt-lem-3-bis}
\int_{\R}\frac{d\xi}{|\xi|^\nu} \frac{1}{\langle\xi-a\rangle^p} \lesssim \frac{1}{\langle a\rangle^{\nu-\varepsilon}},
\end{equation}
and one can choose $\eps=0$ in the case $p>1$.
\smallskip

\noindent
$(iii)$ For all $(a,b)\in \R^2$ such that $|a|\leq |b|$, and for every $0<\varepsilon<1$, it holds that
\begin{equation}\label{adapt-lem-33}
\int_{\R}\frac{|\xi|^\mu d\xi}{\langle\xi-a\rangle \langle \xi-b\rangle} \lesssim \frac{\langle a\rangle^\mu}{\langle b-a\rangle^{1-\mu-\varepsilon}}.
\end{equation}
\smallskip

\noindent
$(iv)$ For all $p\geq 2$ and $|a|\geq 1$, one has
\begin{equation}\label{adapt-lem-4-bis}
\int_{\R}\frac{|\xi|^\mu d\xi}{\langle\xi-a\rangle^p} \lesssim \langle a\rangle^{\mu}.
\end{equation}
\end{lemma}

For the sake of clarity, we have postponed the proof of this technical lemma to the appendix section \ref{subsec:append-1}.

\subsection{Fully deterministic controls on the product}\label{sec:deter-con}

We conclude this preliminary section with the statement of the bilinear estimate allowing to control the  product terms  $\ci_\chi\cm(z,w)$ and  $\ci_\chi\big( \Lambda\cdot w\big)$, for all $\Lambda \in H^\ga(\R)$ (along the result of Proposition \ref{prop:conv-lamb-n}) and $z,w$ in the scale $Z^{s,b}$ at the core of our procedure. Recall that the bilinear operator $\cm$ is defined in~\eqref{defm}.

\begin{proposition}\label{prop:control-m-z-z}
For all $\frac12 < b<\frac58$ and $s \geq 0$, there exists $\mu>0$ such that for every $\tau\in (0,1]$,
\begin{equation}\label{control-m-z-z}
\big\| \chi_\tau\cdot \ci_\chi \cm(z,w)\big\|_{Z^{s,b}}\lesssim \tau^{\mu} \| z\|_{Z^{s,b}}\| w\|_{Z^{s,b}}.
\end{equation}
\end{proposition}

\begin{proof}
Inequality \eqref{control-m-z-z} is a rewriting of a classical  estimate in Bourgain spaces. More precisely, setting $u=e^{-\imath  t \Delta} z$ and $v=e^{-\imath  t \Delta} w$, then \eqref{control-m-z-z} is equivalent to 
\begin{equation}\label{KPV2}
\bigg\|\chi_\tau  \int_0^. \chi_\tau e^{-\imath (.-r)\Delta}\big(u\overline{v}\big)\, dr\bigg\|_{X^{s,b}}\lesssim \tau^\mu\big\|u\big\|_{X^{s,b}} \big\|v\big\|_{X^{s,b}} , \quad \text{for some} \ \mu>0,
\end{equation}
which we now prove. Let $\frac12 < b<\frac58$, then for $\mu>0$ small enough,   by \cite[Lemme 3.2]{Ginibre}, 
\begin{equation*} 
\bigg\|\chi_\tau  \int_0^. \chi_\tau e^{-\imath (.-r)\Delta}F\, dr\bigg\|_{X^{s,b}}\lesssim \tau^\mu\big\|F\big\|_{X^{s,b-1+\mu}}.
\end{equation*}
We can then combine the previous estimate with the Bourgain bilinear estimate
 \begin{equation}\label{bili2}
  \|u \overline{v}\|_{X^{s, b-1+\mu} }  \lesssim     \|u\|_{X^{s, b} }\|v\|_{X^{s, b} },
   \end{equation}
   which implies \eqref{KPV2}. For the sake of completeness, we have added the proof of \eqref{bili2} in the appendix (see Section \ref{sect-bili}).
\end{proof}

By applying Proposition~\ref{prop:control-m-z-z} to the particular case where $z=\Lambda\in H^\ga(\R)$ does not depend on the space variable, we immediately get: 
\begin{proposition}\label{prop:lam-z}
 Let  $\frac12 < b<\frac58$, $\gamma>\frac12$ and $s \geq 0$.  Then for all $\Lambda\in H^\ga(\R)$ and ${w\in Z^{s,b}}$, it holds that
\begin{equation}\label{lam-z}
\big\|\chi_\tau \cdot \ci_\chi\big( \Lambda\cdot w\big)\big\|_{Z^{s,b}} \lesssim \tau^{\mu} \|\Lambda\|_{H^{\ga}(\R)}\|w\|_{Z^{s,b}}
\end{equation}
for some $\mu>0$, and where the product $\Lambda\cdot w$ is naturally understood as $(\Lambda\cdot w)(t,x):=\Lambda(t) w(t,x)$.
\end{proposition}

\section{Study of the tree-elements I: first-order}\label{sec:tree-elements}

We now turn to the investigations of the explicit stochastic terms in \eqref{equa-z}, starting with the (rescaled) linear solution $\<Psi>$.

\subsection{Covariance estimates for the (rescaled) linear solution}\label{subsec:cova}

Following \eqref{defi-luxo-n}, the process~$\<Psi>^{(n)}$ is explicitly defined by the formula
\begin{equation*}
\<Psi>^{(n)}_k(t)=\chi(t)\int_0^t ds \, \chi(s) e^{-\imath s|k|^2}\dot{B}^{(n)}_k(s)=\ci_{\chi}\big( e^{-\imath .|k|^2} \dot{B}^{(n)}_k\big)(t),
\end{equation*}
where $\dot{B}^{(n)}$ is the smooth approximation of the noise $\dot{B}$ introduced in \eqref{defi:approx-b-n}. It is readily checked that the covariance of the latter approximation  is given for every $n\geq 1$ by
\begin{equation}\label{cova-b-dot-n}
\mathbb{E}\big[\dot{B}^{(n)}(t,x) \dot{B}^{(n)}(t',x') \big]=c_H^2\int_{D_n} \frac{d\xi d\eta}{|\xi|^{2H_0-1}|\eta|^{2H_1-1}} e^{-\imath \xi(t-t')} e^{-\imath \eta(x-x')},
\end{equation}
where $D_n$ stands for the integration domain $D_n:=\{(\xi,\eta)\in \R^2: \ |\xi|\leq 2^{2n} \ \text{and} \ |\eta|\leq 2^n\}$. Besides, observe that the difference $\dot{B}^{(n,n+1)}:=\dot{B}^{(n+1)}-\dot{B}^{(n)}$ between two successive elements can be represented along the same pattern, that is
\begin{equation*}
\dot{B}^{(n,n+1)}(t,x):= -c_H\int_{D_{n+1}\backslash D_n} \frac{\xi\, e^{\imath t\xi}}{|\xi|^{H_0+\frac{1}{2}}}\frac{\eta\, e^{\imath x\eta}}{|\eta|^{H_1+\frac{1}{2}}}\, \widehat{W}(d\xi,d\eta).
\end{equation*}
As a result, it holds that
\begin{equation}\label{cova-b-n-n-1}
\mathbb{E}\big[\dot{B}^{(n,n+1)}(t,x) \dot{B}^{(n,n+1)}(t',x') \big]=c_H^2\int_{D_{n+1}\backslash D_n} \frac{d\xi d\eta}{|\xi|^{2H_0-1}|\eta|^{2H_1-1}} e^{-\imath \xi(t-t')} e^{-\imath \eta(x-x')}.
\end{equation}
By setting additionally
\begin{equation}\label{addit}
\<Psi>^{(n,n+1)}_k(t):=\<Psi>^{(n+1)}_k(t)-\<Psi>^{(n)}_k(t)=\ci_{\chi}\big( e^{\imath .|k|^2} \dot{B}^{(n,n+1)}_k\big)(t).
\end{equation}
our main estimates about the covariance of $\<Psi>^{(n)}$ can be stated as follows.

\begin{proposition}\label{prop:estim-cov-1}
Fix $H_0,H_1\in (0,1)$. Then for all $n\geq 1$, $k,k'\in \Z$ and $\la,\la'\in \R$, it holds that
\begin{align}
&\Big|\mathbb{E}\Big[  \cf\big(\<Psi>^{(n)}_{k}\big)(\la)\overline{\cf\big(\<Psi>^{(n)}_{k'}\big)(\la')} \Big]\Big|\lesssim \frac{1}{\langle \la\rangle \langle \la'\rangle }\gga^{(H_0)}_{k,k'}(\la,\la')\Lambda^{(H_{1})}_{k,k'} ,\label{bou:cova-luxo-n}
\end{align}
where the proportional constant does not depend on $n$, and where we have set
\begin{equation*}
\gga^{(H)}_{k,k'}(\la,\la'):=\int_\R \frac{d\xi}{|\xi|^{2H-1}}\frac{1}{\langle \xi-(|k|^2+\la)\rangle} \frac{1}{\langle \xi-(|k'|^2+\la')\rangle},
\end{equation*}
\begin{equation*}
\Lambda^{(H)}_{k,k'}:=\int_{\R} \frac{d\eta}{|\eta|^{2H-1}} \frac{1}{\langle \eta-k\rangle} \frac{1}{\langle \eta-k'\rangle}.
\end{equation*}
Besides, for all  $0\leq \ka\leq \frac12 \min(H_0,H_1)$, $n\geq 1$, $k,k'\in \Z$ and $\la,\la'\in \R$, one has 
\begin{align}
&\Big|\mathbb{E}\Big[  \cf\big(\<Psi>^{(n,n+1)}_{k}\big)(\la)\overline{\cf\big(\<Psi>^{(n,n+1)}_{k'}\big)(\la')} \Big]\Big|\lesssim \frac{2^{-2n\ka}}{\langle \la\rangle \langle \la'\rangle }\Big[\gga^{(H_0-\ka)}_{k,k'}(\la,\la')\Lambda^{(H_{1})}_{k,k'}+\gga^{(H_0)}_{k,k'}(\la,\la')\Lambda^{(H_1-\ka)}_{k,k'}\Big],\label{bou:cova-luxo-n-n+1}
\end{align}
where the proportional constant does not depend on $n$.
\end{proposition}

\begin{proof}
We only focus on the proof of \eqref{bou:cova-luxo-n-n+1}, keeping in mind that \eqref{bou:cova-luxo-n} could be proved with similar arguments. Notice that  the condition $0\leq \ka\leq \frac12 \min(H_0,H_1)$ ensures that the right hand side of \eqref{bou:cova-luxo-n-n+1} is finite. 

\smallskip

Using \eqref{addit} and the representation in Lemma \ref{lem:integr-kernel}, we can write
\begin{align}
\cf\big(\<Psi>^{(n,n+1)}_{k}\big)(\la)&=\int d\la_1 \, \ci_\chi(\la,\la_1) \cf\big( e^{-\imath .|k|^2} \dot{B}^{(n,n+1)}_k\big)(\la_1)\nonumber\\
&=\int d\la_1 \, \ci_\chi(\la,\la_1) \int dt\, e^{-\imath t(\la_1+|k|^2) } \int_{\T} dx \, e^{-\imath x  k} \dot{B}^{(n,n+1)}(t,x).\label{repres-psi}
\end{align}
We can now apply the above covariance formula \eqref{cova-b-n-n-1} and deduce that
\small
\begin{align}
&\mathbb{E}\Big[  \cf\big(\<Psi>^{(n,n+1)}_{k}\big)(\la)\overline{\cf\big(\<Psi>^{(n,n+1)}_{k'}\big)(\la')} \Big]\nonumber\\
&=\int d\la_1d\la_1' \, \ci_\chi(\la,\la_1) \overline{\ci_\chi(\la',\la_1')}\nonumber\\
&\hspace{2cm}\int dt dt'\, e^{-\imath t(\la_1+|k|^2) }e^{\imath t'(\la'_1+|k'|^2) }  \int_{\T} dx dx'\, e^{-\imath x  k} e^{\imath x'  k'}\mathbb{E}\big[\dot{B}^{(n,n+1)}(t,x) \dot{B}^{(n,n+1)}(t',x') \big]\nonumber\\
&=c_H^2\int_{D_{n+1}\backslash D_n} \frac{d\xi d\eta}{|\xi|^{2H_0-1}|\eta|^{2H_1-1}}\nonumber\\
&\int d\la_1d\la_1' \, \ci_\chi(\la,\la_1) \overline{\ci_\chi(\la',\la_1')}\int dt dt'\, e^{-\imath t(\la_1+|k|^2) }e^{\imath t'(\la'_1+|k'|^2) }  \int_{\T} dx dx'\, e^{-\imath x k} e^{\imath x' \cdot k'}e^{-\imath \xi(t-t')} e^{-\imath \eta (x-x')}\nonumber\\
&=c_H^2\int_{D_{n+1}\backslash D_n}  \frac{d\xi d\eta}{|\xi|^{2H_0-1}|\eta|^{2H_1-1}}\int d\la_1d\la_1' \, \ci_\chi(\la,\la_1) \overline{\ci_\chi(\la',\la_1')}\nonumber\\
&\hspace{1cm}\bigg(\int dt \, e^{-\imath t(\la_1+|k|^2+\xi) }\bigg)\bigg(\int dt'\, e^{\imath t'(\la'_1+|k'|^2+\xi) }\bigg)\bigg(  \int_{\T} dx \, e^{-\imath x  (\eta+ k)}\bigg)\bigg(\int_{\T}dx' e^{\imath x'  (\eta+k')}\bigg)\nonumber\\
&=c_H^2\int_{D_{n+1}\backslash D_n}  \frac{d\xi d\eta}{|\xi|^{2H_0-1}|\eta|^{2H_1-1}}\, \ci_\chi(\la,-|k|^2-\xi) \overline{\ci_\chi(\la',-|k'|^2-\xi)}\bigg(  \int_{\T} dx \, e^{-\imath x (\eta+ k)}\bigg)\bigg(\int_{\T}dx' e^{\imath x' (\eta+k')}\bigg).\label{calcul-psi-psibar}
\end{align}
\normalsize
Observe that 
$$\Big| \int_{\T} dx \, e^{-\imath x (\eta+ k)}\Big| \lesssim \frac{1}{\langle \eta+ k\rangle},$$
which, going back to \eqref{calcul-psi-psibar}, easily yields
\small
\begin{align}
&\bigg|\mathbb{E}\Big[  \cf\big(\<Psi>^{(n,n+1)}_{k}\big)(\la)\overline{\cf\big(\<Psi>^{(n,n+1)}_{k'}\big)(\la')} \Big]\bigg|\nonumber\\
&\lesssim \int_{D_{n+1}\backslash D_n}  \frac{d\xi d\eta}{|\xi|^{2H_0-1}|\eta|^{2H_1-1}}\, \big|\ci_\chi(\la,-|k|^2-\xi)\big| \big|\ci_\chi(\la',-|k'|^2-\xi)\big|\bigg|  \int_{\T} dx \, e^{-\imath x  (\eta+ k)}\bigg|\bigg|\int_{\T}dx' e^{\imath x'  (\eta+k')}\bigg|\label{transit-1}\\
&\lesssim \int_{|\xi|\geq 2^n}\int_{\eta\in \R}  \frac{d\xi d\eta}{|\xi|^{2H_0-1}|\eta|^{2H_1-1} \langle \eta+ k\rangle  \langle \eta+ k'\rangle}\, \big|\ci_\chi(\la,-|k|^2-\xi)\big| \big|\ci_\chi(\la',-|k'|^2-\xi)\big| \nonumber\\
&+\int_{\xi\in \R}\int_{|\eta|\geq 2^n}  \frac{d\xi d\eta}{|\xi|^{2H_0-1}|\eta|^{2H_1-1} \langle \eta+ k\rangle  \langle \eta+ k'\rangle}\, \big|\ci_\chi(\la,-|k|^2-\xi)\big| \big|\ci_\chi(\la',-|k'|^2-\xi)\big| \nonumber\\
&\lesssim 2^{-2n\ka}\bigg(\int_{\R} \frac{d\xi}{|\xi|^{2(H_0-\ka)-1}} \, \big|\ci_\chi(\la,-|k|^2-\xi)\big| \big|\ci_\chi(\la',-|k'|^2-\xi)\big|\bigg)\bigg( \int_{\R}\frac{d\eta}{|\eta|^{2H_1-1} \langle \eta+ k\rangle  \langle \eta+ k'\rangle } \bigg)\nonumber\\
&+2^{-2n\ka}\bigg(\int_{\R} \frac{d\xi}{|\xi|^{2H_0-1}} \, \big|\ci_\chi(\la,-|k|^2-\xi)\big| \big|\ci_\chi(\la',-|k'|^2-\xi)\big|\bigg)\bigg( \int_{\R}\frac{d\eta}{|\eta|^{2(H_1-\ka)-1}  \langle \eta+ k\rangle  \langle \eta+ k'\rangle   }\bigg).\label{transit-2}
\end{align}
\normalsize
The bound \eqref{bou:cova-luxo-n-n+1} is now an immediate consequence of \eqref{bou-ci-chi}.
\end{proof}

With similar arguments as in the previous proof, we easily derive the following estimate for the non-conjugate covariance of $\<Psi>^{(n)}$.

\begin{proposition}\label{prop:estim-cov-2}
Fix $H_0,H_1\in (0,1)$. Then for all $n\geq 1$, $k,k'\in \Z$ and $\la,\la'\in \R$, it holds that
\begin{align}
&\Big|\mathbb{E}\Big[  \cf\big(\<Psi>^{(n)}_{k}\big)(\la)\cf\big(\<Psi>^{(n)}_{k'}\big)(\la') \Big]\Big|\lesssim \frac{1}{\langle \la\rangle \langle \la'\rangle } \widetilde{\gga}^{(H_0)}_{k,k'}(\la,\la')\widetilde{\Lambda}^{(H_{1})}_{k,k'} ,\label{bou:cova-luxo-n-2}
\end{align}
where the proportional constant does not depend on $n$, and where we have set
\begin{equation*}
\widetilde{\gga}^{(H)}_{k,k'}(\la,\la'):=\int_\R \frac{d\xi}{|\xi|^{2H-1}}\frac{1}{\langle \xi-(|k|^2+\la)\rangle} \frac{1}{\langle \xi-(-\la'-|k'|^2)\rangle},
\end{equation*}
\begin{equation*}
\widetilde{\Lambda}^{(H)}_{k,k'}:=\int_{\R} \frac{d\eta}{|\eta|^{2H-1}} \frac{1}{\langle \eta-k\rangle} \frac{1}{\langle \eta+k'\rangle}.
\end{equation*}
Besides, for all  $0\leq \ka\leq \frac12 \min(H_0,H_1)$, $n\geq 1$, $k,k'\in \Z$ and $\la,\la'\in \R$, one has
\begin{align}
&\Big|\mathbb{E}\Big[  \cf\big(\<Psi>^{(n,n+1)}_{k}\big)(\la)\cf\big(\<Psi>^{(n,n+1)}_{k'}\big)(\la') \Big]\Big|\lesssim \frac{2^{-2n\ka}}{\langle \la\rangle \langle \la'\rangle }\Big[\widetilde{\gga}^{(H_0-\ka)}_{k,k'}(\la,\la')\widetilde{\Lambda}^{(H_{1})}_{k,k'}+\widetilde{\gga}^{(H_0)}_{k,k'}(\la,\la')\widetilde{\Lambda}^{(H_1-\ka)}_{k,k'}\Big],\label{bou:cova-luxo-n-n+1-2}
\end{align}
where the proportional constant does not depend on $n$.
\end{proposition}

\

The (rescaled) linear solution $\<Psi>^{(n)}$ turns out to be the central object of our analysis, and the two estimates \eqref{bou:cova-luxo-n}-\eqref{bou:cova-luxo-n-2} for its covariance will be extensively used in the sequel.

\begin{remark}\label{rk:luxo-white}
To emphasize the specificity of our fractional situation with respect to previous white-noise models, consider the case of a (spatially regularized) white noise
$$\dot{\mathbf{W}}^{(\al)}(t,x)=\sum_k \langle k\rangle^{-\al} \dot{\beta}^{(k)}_t e^{\imath k x},$$
where the $(\beta^{(k)})$'s stand for independent Brownian motions. In this case, the covariance of the corresponding (rescaled) linear solution $\Psi$ reduces to the expression
\begin{align*}
&\mathbb{E}\Big[  \cf\big(\Psi_{k}\big)(\la)\overline{\cf\big(\Psi_{k'}\big)(\la')} \Big]=\big(\1_{\{k=k'\}}\langle k\rangle^{-2\al}\big) \int d\la_1 \, \ci_\chi(\la,\la_1)\overline{\ci_{\chi}(\la',\la_1)},
\end{align*}
which sharply contrasts with the subtle interactions observed in the right-hand side of \eqref{bou:cova-luxo-n} and echoing the past-dependence properties of the fractional noise.
\end{remark}

\subsection{Proof of Proposition \ref{prop:conv-psi-n}}\label{subsec:first-order}

\

\smallskip

\noindent
$(i)$ Using the result of Proposition \ref{prop:estim-cov-1}, we get that for all $ \al \geq 0$ and \linebreak $0<\ka \leq \frac12  \min(H_0,H_1 ,  \frac{\alpha-[2-(2H_0+H_1)]}{2})$, 
\begin{align*}
&\mathbb{E}\Big[\big\|\<Psi>^{(n,n+1)}\big\|_{L^2(\R;H^{- \al }(\T))}^2\Big]=\sum_k \frac{1}{\langle k\rangle^{2 \al }}\int_{\R} d\la \, \mathbb{E}\Big[ \big| \cf(\<Psi>^{(n,n+1)}_k)(\la)\big|^2 \Big]\\
&\hspace{1cm}\lesssim 2^{-2\ka n}\sum_k \frac{1}{\langle k\rangle^{2 \al }}\bigg[\int_{\R} \frac{d\la}{\langle \la\rangle^2}  \int \frac{d\xi}{|\xi|^{2(H_0-\ka)-1}} \frac{1}{\langle \la-(|k|^2-\xi)\rangle^2} \int\frac{d\eta}{|\eta|^{2H_1-1}} \frac{1}{\langle \eta+k\rangle^2}\\
&\hspace{3cm}+\int_{\R} \frac{d\la}{\langle \la\rangle^2}  \int \frac{d\xi}{|\xi|^{2H_0-1}} \frac{1}{\langle \la-(|k|^2-\xi)\rangle^2} \int\frac{d\eta}{|\eta|^{2(H_1-\ka)-1}} \frac{1}{\langle \eta+k\rangle^2}\bigg]\\
&\hspace{1cm}\lesssim 2^{-2\ka n} \sum_k \frac{1}{\langle k\rangle^{2 \al }}  \bigg[\int \frac{d\xi}{|\xi|^{2(H_0-\ka)-1}}\int_{\R} \frac{d\la}{\langle \la\rangle^2}  \frac{1}{\langle \la-(|k|^2-\xi)\rangle^2} \int\frac{d\eta}{|\eta|^{2H_1-1}} \frac{1}{\langle \eta+k\rangle^2}\\
&\hspace{3cm}+ \int \frac{d\xi}{|\xi|^{2H_0-1}} \int_{\R} \frac{d\la}{\langle \la\rangle^2} \frac{1}{\langle \la-(|k|^2-\xi)\rangle^2} \int\frac{d\eta}{|\eta|^{2(H_1-\ka)-1}} \frac{1}{\langle \eta+k\rangle^2}\bigg]\\
&\hspace{1cm}\lesssim 2^{-2\ka n} \sum_k  \frac{1}{\langle k\rangle^{2 \al }} \bigg[\bigg(\int \frac{d\xi}{|\xi|^{2(H_0-\ka)-1}} \frac{1}{\langle |k|^2-\xi\rangle^2}\bigg) \bigg( \int\frac{d\eta}{|\eta|^{2H_1-1}} \frac{1}{\langle \eta+k\rangle^2}\bigg)\\
&\hspace{3cm}+\bigg(\int \frac{d\xi}{|\xi|^{2H_0-1}} \frac{1}{\langle |k|^2-\xi\rangle^2}\bigg) \bigg( \int\frac{d\eta}{|\eta|^{2(H_1-\ka)-1}} \frac{1}{\langle \eta+k\rangle^2}\bigg)\bigg].
\end{align*}
Now we can apply the inequalities  \eqref{adapt-lem-3-bis} and \eqref{adapt-lem-4-bis}  to    assert that for every $\varepsilon >0$ small enough, $k\neq 0$, and  $\mu \in \{ 2H_0-1 , 2(H_0-\ka)-1 \}$,  $\nu \in \{ 2H_1-1  , 2(H_1-\ka)-1 \}$
$$\int \frac{d\xi}{|\xi|^{\mu}} \frac{1}{\langle |k|^2-\xi\rangle^2} \lesssim \frac{1}{|k|^{2\mu-\eps}}  \quad  \int\frac{d\eta}{|\eta|^{\nu}} \frac{1}{\langle \eta+k\rangle^2}  \lesssim \frac{1}{|k|^{\nu-\eps}},$$
then  
\begin{align*}
&\mathbb{E}\Big[\big\|\<Psi>^{(n,n+1)}\big\|_{L^2(\R;H^{- \al }(\T))}^2\Big]\\
&\lesssim 2^{-2\ka n}\bigg( 1+\sum_{k\neq 0} \frac{1}{| k|^{2 \al }} \bigg[ \frac{1}{| k|^{4H_0-4\ka-2-\varepsilon}}\frac{1}{| k|^{2H_1-1-\varepsilon}}+\frac{1}{| k|^{4H_0-2-\varepsilon}}\frac{1}{| k|^{2H_1-2\ka-1-\varepsilon}}\bigg]\bigg).
\end{align*}
It is readily checked that the sum into brackets is finite as soon as $2 \al +(4H_0-4\ka-2-\varepsilon)+(2H_1-1-\varepsilon)>1$, that is $ \al >2-2H_0-H_1+2\ka+\varepsilon$, where $\varepsilon>0$ is chosen small enough. As a result, we get the bound
\begin{align*}
&\mathbb{E}\Big[\big\|\<Psi>^{(n,n+1)}\big\|_{L^2(\R;H^{- \al }(\T))}^2\Big]\lesssim 2^{-2\ka n}.
\end{align*}
Once endowed with this moment estimate, the almost sure convergence of $\<Psi>^{(n)}$ in $L^2\big(\R;H^{- \al }(\T)\big)$ follows from an elementary Borel-Cantelli argument.  Using  the fact that the random variables under consideration are Gaussians, it is then possible to upgrade the convergence in  $\mathcal{C}\big(\R;W^{- \al, p}(\T)\big)$, for any $1 \leq p \leq \infty$. We refer to \cite[Paragraph 2.2]{DST} for the details of the proof.

\

\noindent
 $(ii)$ Going back to the definition of $\<Psi>^{(n)}$ in \eqref{defi-luxo-n} and since $\chi\equiv 1$ on $[0,1]$, we have
\small
\begin{align}
&\mathbb{E}\Big[\big\|\<Psi>^{(n)}\big\|_{L^2([0,T] ; H^{-\alpha}(\mathbb{T}))}^2\Big]\geq \sum_{k\geq 1} \langle k\rangle^{-2 \alpha} \int_{0}^{T\wedge 1} dt \,  \mathbb{E}\bigg[\bigg|\int_0^t ds \, e^{-\imath s|k|^2} \dot{B}^{(n)}_k(s)\bigg|^2\bigg]\nonumber\\
&\geq \sum_{k\geq 1} \langle k\rangle^{-2 \alpha}  \int_{0}^{T\wedge 1} dt \,  \int_0^t\int_0^t ds ds' \, e^{-\imath (s-s')|k|^2} \mathbb{E}\Big[\dot{B}^{(n)}_k(s)\overline{\dot{B}^{(n)}_k(s')}\Big]\label{refe-dive}\\
&\geq \sum_{k\geq 1}  \langle k\rangle^{-2 \alpha} \int_{0}^{T\wedge 1} dt \,  \int_0^t\int_0^t ds ds' \, e^{-\imath (s-s')|k|^2}\int_{\mathbb{T}}dx\int_{\mathbb{T}}dx' \, e^{-\imath k(x-x')}   \mathbb{E}\Big[\dot{B}^{(n)}(s,x)\dot{B}^{(n)}(s',x)\Big].\nonumber
\end{align}
\normalsize
We can now apply the covariance identity in \eqref{cova-b-dot-n}, which yields
\small
\begin{align}
&\mathbb{E}\Big[\big\|\<Psi>^{(n)}\big\|_{L^2([0,T] ; H^{-\alpha}(\mathbb{T}))}^2\Big]\nonumber\\
&\gtrsim \sum_{k\geq 1} \langle k\rangle^{-2 \alpha}\int_{|\xi|\leq 2^{2n}} \frac{d\xi}{|\xi|^{2H_0-1}}\int_{|\eta|\leq 2^{n}} \frac{d\eta}{|\eta|^{2H_1-1}}\nonumber\\
&\hspace{2cm}\int_{0}^{T\wedge 1} dt \,  \int_0^t\int_0^t ds ds' \, e^{-\imath (s-s')|k|^2}\int_{\mathbb{T}}dx\int_{\mathbb{T}}dx' \, e^{-\imath k(x-x')}  e^{-\imath \xi(s-s')} e^{-\imath \eta (x-x')}\nonumber\\
&\gtrsim \sum_{k\geq 1} \langle k\rangle^{-2 \alpha}\bigg[\int_{|\xi|\leq 2^{2n}} \frac{d\xi}{|\xi|^{2H_0-1}} \int_{0}^{T\wedge 1} dt\, \bigg| \int_0^t ds  \, e^{\imath s(|k|^2-\xi)}\bigg|^2\bigg] \bigg[\int_{|\eta|\leq 2^{n}}\frac{d\eta}{|\eta|^{2H_1-1}}\bigg| \int_{\mathbb{T}}dx \, e^{-\imath x(k+\eta)}\bigg|^2\bigg]\nonumber\\
&\gtrsim \sum_{1\leq k \leq 2^{n-1}}\langle k\rangle^{-2 \alpha} \bigg[\int_{\frac12 |k|^2\leq |\xi|\leq \frac32 |k|^2} \frac{d\xi}{|\xi|^{2H_0-1}} \int_{0}^{T\wedge 1} dt\, \bigg| \int_0^t ds  \, e^{\imath s(|k|^2-\xi)}\bigg|^2\bigg] \nonumber \\
&\hspace{8cm} \bigg[\int_{\frac12 |k|\leq |\eta|\leq \frac32 |k|}\frac{d\eta}{|\eta|^{2H_1-1}}\bigg| \int_{\mathbb{T}}dx \, e^{-\imath x(k+\eta)}\bigg|^2\bigg].\label{diver-l-2}
\end{align}
\normalsize
At this point, observe that for every $k\geq 1$,
\begin{align*}
&\int_{\frac12 k^2\leq |\xi|\leq \frac32 k^2} \frac{d\xi}{|\xi|^{2H_0-1}} \int_{0}^{T\wedge 1} dt\, \bigg| \int_0^t ds  \, e^{\imath s(|k|^2-\xi)}\bigg|^2\\
&=k^{4-4H_0}\int_{\frac12 \leq |\xi|\leq \frac32 } \frac{d\xi}{|\xi|^{2H_0-1}} \int_{0}^{T\wedge 1} dt\, \bigg| \int_0^t ds  \, e^{\imath s|k|^2(1-\xi)}\bigg|^2\\
&\gtrsim k^{4-4H_0}\int_{0 \leq \xi\leq \frac{1}{2k^2} } d\xi \int_{0}^{T\wedge 1} dt\, \bigg| \int_0^t ds  \, e^{\imath s|k|^2\xi }\bigg|^2\gtrsim \frac{1}{k^{4H_0-2}}\int_{0 \leq \xi\leq \frac{1}{2} } d\xi \int_{0}^{T\wedge 1} dt\, \bigg| \int_0^t ds  \, e^{\imath s\xi }\bigg|^2\gtrsim \frac{1}{k^{4H_0-2}}.
\end{align*}
In the same way, for every $k\geq 1$,
\begin{align*}
&\int_{\frac12 k\leq |\eta|\leq \frac32 k}\frac{d\eta}{|\eta|^{2H_1-1}}\bigg| \int_{\mathbb{T}}dx \, e^{-\imath x(k+\eta)}\bigg|^2\\
&=k^{2-2H_1}\int_{\frac12 \leq |\eta|\leq \frac32 }\frac{d\eta}{|\eta|^{2H_1-1}}\bigg| \int_{\mathbb{T}}dx \, e^{-\imath xk(1+\eta)}\bigg|^2\\
&\gtrsim k^{2-2H_1}\int_{0\leq \eta\leq \frac{1}{2k} }d\eta\bigg| \int_{\mathbb{T}}dx \, e^{-\imath xk\eta}\bigg|^2\gtrsim \frac{1}{k^{2H_1-1}}\int_{0\leq \eta\leq \frac{1}{2} }d\eta\bigg| \int_{\mathbb{T}}dx \, e^{-\imath x\eta}\bigg|^2\gtrsim \frac{1}{k^{2H_1-1}}.
\end{align*}
Going back to \eqref{diver-l-2}, we obtain
\begin{equation*}
\mathbb{E}\Big[\big\|\<Psi>^{(n)}\big\|_{L^2([0,T] ; H^{-\alpha}(\mathbb{T}))}^2\Big]\gtrsim \sum_{1\leq k \leq 2^{n-1}} \frac{1}{k^{2\alpha+4H_0+2H_1-3}}\gtrsim 2^{2n(2-2H_0-H_1-\alpha)}, 
\end{equation*}
where we have used the fact that $2\alpha+4H_0+2H_1-3<1$ to derive the last inequality. This concludes the proof of Proposition \ref{prop:conv-psi-n}.

\

\section{Study of the tree-elements II: second-order}\label{sec:tree-elements-2}

This section is devoted to the proof of the assertion $(i)$ in Proposition \ref{prop:ice-cream}. In other words, we are here interested in the construction of $\<IPsi2>$, one of our main objectives in the study. \\
\indent For more clarity, let us introduce the (twisted) product operator $\widetilde{\cm}$ defined for all $v,w \in \mathcal{S}'(\T)$ and $k\in \Z$ by
\begin{equation*} 
\widetilde{\cm}(v,w)_k(t)=\sum_{k_1\in \Z\backslash \{0\}} e^{\imath t\Omega_{k,k_1}} v_{k+k_1} \overline{w_{k_1}}.
\end{equation*}
With this notation, we can recast the definition \eqref{defi:ipsi2-main} of $\<IPsi2>^{(n)}$ into
\begin{equation*}
\<IPsi2>^{(n)}_k(t):=\ci_\chi \Big(\widetilde{\cm}( \<Psi>^{(n)},\<Psi>^{(n)})_k(.)-{\mathbb{E}\big[\widetilde{\cm}( \<Psi>^{(n)},\<Psi>^{(n)})_0(.)\big]}\Big)(t).
\end{equation*}

For the sake of clarity, we will in fact divide the proof of the convergence of $(\<IPsi2>^{(n)})_{n\geq 1}$ into two successive propositions.

\begin{proposition}\label{prop:ice-cream-1}
Assume that $\frac12 < H_0<1$ and  $\frac14< H_1<1$ satisfy  $\frac74 <2H_0+H_1< 2$, and set $\underline{H}:=2H_0+H_1-1\in (\frac34,1)$.  Then for every pair $(s,b)$ satisfying the conditions
\begin{equation}\label{assump-lau-statement-proof}
 0<s<\underline{H}-\frac12 \quad \text{and} \quad \frac12<b< \min\big({\underline{H}}-s,\frac34-\frac{s}{2}\big),
\end{equation}
the sequence of processes $(\mathbf{\Psi}^{(n)})_{n\geq 1}$ defined in Fourier mode by 
\begin{equation}\label{exp-til}
\mathbf{\Psi}^{(n)}_k(t):=\ci_\chi \Big(\widetilde{\cm}( \<Psi>^{(n)},\<Psi>^{(n)})_k(.)-{\mathbb{E}\big[\widetilde{\cm}( \<Psi>^{(n)},\<Psi>^{(n)})_k(.)\big]}\Big)(t) 
\end{equation}
converges almost surely in $Z^{s,b}$. 

\end{proposition}

\begin{proposition}\label{prop:ice-cream-2}
Assume that $\frac12 < H_0<1$ and  $\frac14< H_1<1$ satisfy  $\frac74 <2H_0+H_1< 2$, and set $\underline{H}:=2H_0+H_1-1\in (\frac34,1)$. Then for any $(s,b)$ satisfying
\begin{equation}\label{conditi}
0<s<\frac12 \quad \text{and} \quad \frac12< b <\min\big(2\underline{H}-1,1-s\big),
\end{equation} 
the sequence of functions $(f^{(n)})_{n\geq 1}$ defined in Fourier mode by 
$$f^{(n)}_k(t):=\ci_\chi \Big( \mathbb{E}\big[\widetilde{\cm}( \<Psi>^{(n)},\<Psi>^{(n)})_k(.)\big]-\mathbb{E}\big[\widetilde{\cm}( \<Psi>^{(n)},\<Psi>^{(n)})_0(.)\big]\Big)(t) $$
converges in $Z^{s,b}$.
\end{proposition}

\subsection{Proof of Proposition \ref{prop:ice-cream-1}} Setting 
$$\mathbf{\Psi}^{(n,n+1)}:=\mathbf{\Psi}^{(n+1)}-\mathbf{\Psi}^{(n)},$$
the following moment estimate allows us to guarantee the almost sure convergence of $(\mathbf{\Psi}^{(n)})_{n\geq 1}$ in $Z^{s,b}$ (as an easy consequence of the Borel-Cantelli lemma).

\begin{proposition}\label{pro-mom}
Under the hypotheses of Proposition \ref{prop:ice-cream-1}, it holds that
\begin{equation}\label{gene-ice-cream}
\mathbb{E}\Big[ \big\|\mathbf{\Psi}^{(n,n+1)}\big\|^2_{Z^{s,b}} \Big] \lesssim 2^{-\ka n},
\end{equation}
where $\ka:=\frac12 \min\big(H_0,H_1\big)>0$. 
\end{proposition}

\begin{proof}[Proof of Proposition \ref{pro-mom}]

We will focus on the proof of the uniform estimate
\begin{equation}\label{simpl-unif}
\mathbb{E}\Big[ \big\|\mathbf{\Psi}^{(n)}\big\|^2_{Z^{s,b}} \Big] \lesssim 1,
\end{equation}
and leave it to the reader to check that the subsequent arguments could be extended toward~\eqref{gene-ice-cream} by following the same basic estimation procedure as in the transition from \eqref{transit-1} to \eqref{transit-2}.

\smallskip

By using the representation formula \eqref{rel-def-i-chi} of $\ci_\chi$, we can first write
\begin{align*}
&\cf\big(\ci_\chi\big(\widetilde{\cm}( \<Psi>^{(n)},\<Psi>^{(n)})_k\big)\big)(\la)=\int_{\R} d\la' \, \ci_\chi(\la,\la') \cf(\widetilde{\cm}(\<Psi>^{(n)},\<Psi>^{(n)})_k)(\la')\\
&=\sum_{k_1\neq 0} \int_{\R} d\la' \, \ci_\chi(\la,\la') \int dt\, e^{-\imath \la' t}e^{\imath t \Omega_{k,k_1}} \<Psi>^{(n)}_{k+k_1}(t) \overline{\<Psi>^{(n)}_{k_1}(t)}\\
&=\sum_{k_1\neq 0}\int d\la_1 \, \overline{\cf(\<Psi>^{(n)}_{k_1})(\la_1)} \int d\la_2 \, \cf(\<Psi>^{(n)}_{k+k_1})(\la_2) \int_{\R} d\la' \, \ci_\chi(\la,\la')\int dt\, e^{-\imath \la' t}e^{\imath t\Omega_{k,k_1}} e^{\imath t\la_2}e^{-\imath t \la_1},
\end{align*}
and thus
\begin{equation*}
\cf\big(\ci_\chi\big(\widetilde{\cm}( \<Psi>^{(n)},\<Psi>^{(n)})_k\big)\big)(\la)=\sum_{k_1\neq 0}\int d\la_1d\la_2 \, \ci_\chi(\la,\Omega_{k,k_1}+\la_2-\la_1) \overline{\cf(\<Psi>^{(n)}_{k_1})(\la_1)} \cf(\<Psi>^{(n)}_{k+k_1})(\la_2).
\end{equation*}
Based on this expression, and by applying Wick's formula, we can compute
\begin{align*}
&\mathbb{E}\Big[\big| \cf \big(\ci_\chi \big(\widetilde{\cm}( \<Psi>^{(n)},\<Psi>^{(n)})_k\big)(\la)\big|^2\Big]\\
&=\sum_{k_1\neq 0}\int d\la_1d\la_2 \sum_{k'_1\neq 0}\int d\la'_1d\la'_2  \ \ci_\chi(\la,\Omega_{k,k_1}+\la_2-\la_1) \overline{\ci_\chi(\la,\Omega_{k,k'_1}+\la'_2-\la'_1)} \\
&\hspace{3cm}\mathbb{E}\Big[\overline{\cf(\<Psi>^{(n)}_{k_1})(\la_1)} \cf(\<Psi>^{(n)}_{k+k_1})(\la_2)\cf(\<Psi>^{(n)}_{k'_1})(\la'_1) \overline{\cf(\<Psi>^{(n)}_{k+k_1'})(\la'_2)}\Big]\\
&=\Big|\mathbb{E}\big[\cf \big(\ci_\chi \big(\widetilde{\cm}( \<Psi>^{(n)},\<Psi>^{(n)})_k\big)(\la)\big]\Big|^2\\
&\hspace{0.5cm}+\sum_{k_1\neq 0}\int d\la_1d\la_2 \sum_{k'_1\neq 0}\int d\la'_1d\la'_2  \ \ci_\chi(\la,\Omega_{k,k_1}+\la_2-\la_1) \overline{\ci_\chi(\la,\Omega_{k,k'_1}+\la'_2-\la'_1)} \\
&\hspace{3.5cm}\mathbb{E}\Big[\overline{\cf(\<Psi>^{(n)}_{k_1})(\la_1)} \cf(\<Psi>^{(n)}_{k'_1})(\la'_1) \Big]\mathbb{E}\Big[ \cf(\<Psi>^{(n)}_{k+k_1})(\la_2)\overline{\cf(\<Psi>^{(n)}_{k+k_1'})(\la'_2)}\Big]\\
&\hspace{0.5cm}+\sum_{k_1\neq 0}\int d\la_1d\la_2 \sum_{k'_1\neq 0}\int d\la'_1d\la'_2  \ \ci_\chi(\la,\Omega_{k,k_1}+\la_2-\la_1) \overline{\ci_\chi(\la,\Omega_{k,k'_1}+\la'_2-\la'_1)} \\
&\hspace{3.5cm}\mathbb{E}\Big[\overline{\cf(\<Psi>^{(n)}_{k_1})(\la_1)} \overline{\cf(\<Psi>^{(n)}_{k+k_1'})(\la'_2)}\Big]\mathbb{E}\Big[\cf(\<Psi>^{(n)}_{k+k_1})(\la_2)\cf(\<Psi>^{(n)}_{k'_1})(\la'_1) \Big].
\end{align*}
Thus, going back to the expression \eqref{exp-til} of $\mathbf{\Psi}^{(n)}$, we deduce the decomposition
\begin{align*}
&\mathbb{E}\Big[\big| \cf \big(\mathbf{\Psi}^{(n)}_k\big)(\la)\big|^2 \Big]=I^{(n)}_k(\la)+J^{(n)}_k(\la),
\end{align*}
with
\small
\begin{align*}
&I^{(n)}_k(\la):=\sum_{k_1\neq 0}\int d\la_1d\la_2 \sum_{k'_1\neq 0}\int d\la'_1d\la'_2  \ \ci_\chi(\la,\Omega_{k,k_1}+\la_2-\la_1) \overline{\ci_\chi(\la,\Omega_{k,k'_1}+\la'_2-\la'_1)} \\
&\hspace{5cm}\mathbb{E}\Big[\overline{\cf(\<Psi>^{(n)}_{k_1})(\la_1)} \cf(\<Psi>^{(n)}_{k'_1})(\la'_1) \Big]\mathbb{E}\Big[ \cf(\<Psi>^{(n)}_{k+k_1})(\la_2)\overline{\cf(\<Psi>^{(n)}_{k+k_1'})(\la'_2)}\Big]
\end{align*}
\normalsize
and
\small
\begin{align*}
&J^{(n)}_k(\la):=\sum_{k_1\neq 0}\int d\la_1d\la_2 \sum_{k'_1\neq  0}\int d\la'_1d\la'_2  \ \ci_\chi(\la,\Omega_{k,k_1}+\la_2-\la_1) \overline{\ci_\chi(\la,\Omega_{k,k'_1}+\la'_2-\la'_1)} \\
&\hspace{5cm}\mathbb{E}\Big[\overline{\cf(\<Psi>^{(n)}_{k_1})(\la_1)} \overline{\cf(\<Psi>^{(n)}_{k+k_1'})(\la'_2)}\Big]\mathbb{E}\Big[\cf(\<Psi>^{(n)}_{k+k_1})(\la_2)\cf(\<Psi>^{(n)}_{k'_1})(\la'_1) \Big].
\end{align*}
\normalsize

\smallskip

The following sections \ref{subse:i-n} and \ref{subse:ii-n} are devoted to the exhibition of uniform bounds for these two quantities, when evaluated in the $Z^{s,b}$-topology.

\subsubsection{Estimates related to $I^{(n)}$}\label{subse:i-n}

We shall prove that
 \begin{equation}\label{cj-fini}
\sup_{n\geq 0}\big\| I^{(n)} \big\|_{\widehat{Z}^{s,b}} <\infty,
\end{equation}
where, for more clarity, we set from now on 
\begin{equation}\label{z-chapo}
\big\|\mathfrak{f} \big\|_{\widehat{Z}^{s,b}}:=\sum_k  \langle k \rangle^{2s}\int d\la \, \langle \la \rangle^{2b}\big| \mathfrak{f}_k(\la)\big|.
\end{equation}

\smallskip

By using \eqref{bou-ci-chi} and the estimate of Proposition \ref{prop:estim-cov-1}, we can first assert that  
\small
\begin{align*}
\big|I^{(n)}_k(\la)\big|
&\lesssim \sum_{k_1,k_1'\in \Z\backslash \{0\}}\int d\la_1d\la_2 \int d\la'_1d\la'_2  \ \Big|\ci_\chi(\la,\Omega_{k,k_1}+\la_2-\la_1) \overline{\ci_\chi(\la,\Omega_{k,k'_1}+\la'_2-\la'_1)}\Big| \\
&\hspace{4cm}\Big|\mathbb{E}\Big[\overline{\cf(\<Psi>^{(n)}_{k_1})(\la_1)} \cf(\<Psi>^{(n)}_{k'_1})(\la'_1) \Big]\mathbb{E}\Big[ \cf(\<Psi>^{(n)}_{k+k_1})(\la_2)\overline{\cf(\<Psi>^{(n)}_{k+k_1'})(\la'_2)}\Big]\Big|\\
&\lesssim \frac{1}{\langle \la \rangle^2}\sum_{k_1,k_1'\in \Z\backslash \{0\}} \int \frac{d\la_1d\la_2 }{\langle \la_1\rangle  \langle \la_2\rangle \langle \la-\Omega_{k,k_1}+\la_1-\la_2\rangle}\int \frac{d\la'_1d\la'_2 }{\langle \la_1'\rangle \langle \la_2'\rangle \langle \la-\Omega_{k,k_1'}+\la'_1-\la'_2\rangle}   \\
&\hspace{3.5cm}\Big[\gga^{(H_0)}_{k_1,k_1'}(\la_1,\la_1') \Lambda^{(H_{1})}_{k_1,k_1'}\Big] \Big[\gga^{(H_0)}_{k+k_1,k+k_1'}(\la_2,\la_2') \Lambda^{(H_{1})}_{k+k_1,k+k_1'}\Big],
\end{align*}
\normalsize
where the proportional constant no longer depends on $n$. 

\smallskip

Set $\delta:=\min(2H_1,1)>\frac12$. By applying the result of Lemma \ref{lem:int-singu-zero} (with $\nu=2H_1-1>0$ in  the case $H_1>\frac 12$ and  with $\mu=1-2H_1\geq 0$ in  the case $H_1\leq \frac 12$)   to the expression of $\Lambda^{(H_1)}$, we derive that for every $\varepsilon >0$ small enough, 
\begin{equation}\label{ineq-58}
 \Lambda^{(H_{1})}_{k,\ell} \lesssim       \frac{1}{\langle |k|\wedge |\ell|\rangle^{2H_1-1} \langle k-\ell \rangle^{\delta-\varepsilon}}     . 
\end{equation}
Therefore
\small
\begin{align}
&\big|I^{(n)}_k(\la)\big|\nonumber\\
&\lesssim  \frac{1}{\langle \la \rangle^2}\sum_{k_1,k_1'\in \Z\backslash \{0\}} \int \frac{d\la_1d\la_2 }{\langle \la_1\rangle  \langle \la_2\rangle \langle \la-\Omega_{k,k_1}+\la_1-\la_2\rangle}\int \frac{d\la'_1d\la'_2 }{\langle \la_1'\rangle \langle \la_2'\rangle \langle \la-\Omega_{k,k'_1}+\la'_1-\la'_2\rangle}  \nonumber  \\
&\hspace{1cm}\gga^{(H_0)}_{k_1,k_1'}(\la_1,\la_1')\gga^{(H_0)}_{k+k_1,k+k_1'}(\la_2,\la_2')\frac{1}{\langle |k_1|\wedge |k_1'|\rangle^{2H_1-1}}\frac{1}{\langle |k+k_1|\wedge |k+k_1'|\rangle^{2H_1-1}} \frac{1}{\langle k_1'-k_1 \rangle^{2\delta-2\varepsilon}},\nonumber
\end{align}
\normalsize
and so for all $k \in \Z$
\small
\begin{align}
& \int d\la \, \langle \la\rangle^{2b} \, | I^{(n)}_k(\la)|\nonumber\\
&\lesssim \sum_{k_1,k_1'\in \Z\backslash \{0\}} \int \frac{d\la_1d\la_2 }{\langle \la_1\rangle  \langle \la_2\rangle }\int \frac{d\la'_1d\la'_2 }{\langle \la_1'\rangle \langle \la_2'\rangle } \bigg[ \int \frac{d\la}{\langle \la\rangle^{2-2b}} \frac{1}{\langle \la-\Omega_{k,k_1}+\la_1-\la_2\rangle}\frac{1}{\langle \la-\Omega_{k,k'_1}+\la'_1-\la'_2\rangle}\bigg]  \nonumber  \\
&\hspace{1cm}\gga^{(H_0)}_{k_1,k_1'}(\la_1,\la_1')\gga^{(H_0)}_{k+k_1,k+k_1'}(\la_2,\la_2')\frac{1}{\langle |k_1|\wedge |k_1'|\rangle^{2H_1-1}}\frac{1}{\langle |k+k_1|\wedge |k+k_1'|\rangle^{2H_1-1}} \frac{1}{\langle k_1'-k_1 \rangle^{2\delta-2\varepsilon}},\nonumber\\
&\lesssim  \sum_{k_1,k_1'\in \Z\backslash \{0\}}\frac{1}{\langle |k_1|\wedge |k_1'|\rangle^{2H_1-1}}\frac{1}{\langle |k+k_1|\wedge |k+k_1'|\rangle^{2H_1-1}} \frac{1}{\langle k_1'-k_1 \rangle^{2\delta -2\varepsilon}}\nonumber\\
& \int \frac{d\la_1d\la_2 }{\langle \la_1\rangle  \langle \la_2\rangle } \int \frac{d\la'_1d\la'_2 }{\langle \la_1'\rangle \langle \la_2'\rangle }\bigg[ \frac{1}{\langle \Omega_{k,k_1}-\la_1+\la_2\rangle^{1-b}}\frac{1}{\langle \Omega_{k,k'_1}-\la'_1+\la'_2\rangle^{1-b}}\bigg]  \gga^{(H_0)}_{k_1,k_1'}(\la_1,\la_1')\gga^{(H_0)}_{k+k_1,k+k_1'}(\la_2,\la_2')\label{correc-lau}. 
\end{align}
\normalsize

\smallskip

Let us now apply  \eqref{borne-00} with $\lambda= \Omega_{k,k_1}$, then combining this estimate with the subsequent Lemma~\ref{lem:lau}, we deduce 
\begin{multline*}
\int \frac{d\la_1d\la_2 }{\langle \la_1\rangle  \langle \la_2\rangle }\int \frac{d\la'_1d\la'_2 }{\langle \la_1'\rangle \langle \la_2'\rangle } \bigg[ \frac{1}{\langle \Omega_{k,k_1}-\la_1+\la_2\rangle^{1-b}}\frac{1}{\langle \Omega_{k,k'_1}-\la'_1+\la'_2\rangle^{1-b}}\bigg]  \\
\hspace{7cm}\gga^{(H_0)}_{k_1,k_1'}(\la_1,\la_1')\gga^{(H_0)}_{k+k_1,k+k_1'}(\la_2,\la_2') \\
\begin{aligned}
&\lesssim \frac{1}{\langle\Omega_{k,k_1} \rangle^{1-b}\langle\Omega_{k,k'_1} \rangle^{1-b}}\int \frac{d\la_1 d\la'_1 }{\langle \la_1\rangle^b  \langle \la_1'\rangle^b }\Big|\gga^{(H_0)}_{k_1,k_1'}(\la_1,\la_1')\Big|   \int \frac{d\la_2 d\la'_2 }{\langle \la_2\rangle^b  \langle \la_2'\rangle^b }\Big|\gga^{(H_0)}_{k+k_1,k+k_1'}(\la_2,\la_2')\Big| \\
&\lesssim \frac{1}{\langle\Omega_{k,k_1} \rangle^{1-b}\langle\Omega_{k,k'_1} \rangle^{1-b}}\frac{1}{\langle |k_1|\wedge |k'_1|\rangle^{2(2H_0-1)}}\frac{1}{\langle |k+k_1|\wedge |k+k'_1|\rangle^{2(2H_0-1)}}.
\end{aligned}
\end{multline*}

\smallskip

By injecting the above estimates into \eqref{correc-lau} and setting $\underline{H}:=2H_0+H_1-1$, we get that
\small
\begin{align}
& \int d\la \, \langle \la\rangle^{2b} \, |I^{(n)}_k(\la) |\nonumber\\
&\lesssim\sum_{k_1,k_1'\in \Z\backslash \{0\}}\frac{1}{\langle |k_1|\wedge |k_1'|\rangle^{2\underline{H}-1}}\frac{1}{\langle |k+k_1|\wedge |k+k_1'|\rangle^{2\underline{H}-1}} \frac{1}{\langle k_1'-k_1 \rangle^{2\delta-2\varepsilon}}\frac{1}{\langle \Omega_{k,k_1}\rangle^{1-b}}\frac{1}{\langle \Omega_{k,k'_1}\rangle^{1-b}}.\label{majo-imp-i}
\end{align}
\normalsize
We now treat the contributions $k \neq 0$ and $k=0$ separately. \medskip

$\bullet$ Terms with $k \neq 0$. We have
\begin{align}
&\sum_{k \neq 0} \langle k\rangle^{2s}\int d\la \, \langle \la\rangle^{2b} \, | I^{(n)}_k(\la)|\nonumber\\
&\lesssim  \sum_{k_1,k_1'}\frac{1}{\langle |k_1|\wedge |k_1'|\rangle^{2\underline{H}+1-2b}} \frac{1}{\langle k_1'-k_1 \rangle^{2\delta-2\varepsilon}}\bigg(\sum_{k} \frac{1}{\langle k\rangle^{2-2b-2s}} \frac{1}{\langle |k+k_1|\wedge |k+k_1'|\rangle^{2\underline{H}-1}}\bigg).\label{correc-lau-2}
\end{align}
\normalsize
Recall that $\frac34<\underline{H}<1$, so $\frac12 <2\underline{H}-1<1$, and by \eqref{assump-lau-statement-proof}, one has $0<2-2b-2s<1$ as well as $\big(2-2b-2s\big)+\big(2\underline{H}-1\big)=1+2\big(\underline{H}-b-s\big)>1$.
Therefore
\begin{align*}
&\sum_{k \neq 0} \frac{1}{\langle k\rangle^{2-2b-2s}} \frac{1}{\langle |k+k_1|\wedge |k+k_1'|\rangle^{2\underline{H}-1}}\\
&\lesssim \sum_{k\in \Z} \frac{1}{\langle k\rangle^{2-2b-2s}}\frac{1}{\langle k+k_1\rangle^{2\underline{H}-1}}+\sum_{k\in \Z} \frac{1}{\langle k\rangle^{2-2b-2s}}\frac{1}{\langle k+k_1'\rangle^{2\underline{H}-1}}\lesssim \frac{1}{\langle |k_1|\wedge |k_1'|\rangle^{2\underline{H}-2b-2s}}.
\end{align*}
Going back to \eqref{correc-lau-2}, we deduce
\small
\begin{align*}
&\sup_{n\geq 0}\sum_{k \neq 0} \langle k\rangle^{2s}\int d\la \, \langle \la\rangle^{2b} \, | I^{(n)}_k(\la)|\nonumber\\
&\lesssim  \sum_{k_1,k_1'}\frac{1}{\langle |k_1|\wedge |k_1'|\rangle^{1+4\underline{H}-4b-2s}} \frac{1}{\langle k_1'-k_1 \rangle^{2\delta-2\varepsilon}}\lesssim  \bigg(\sum_{k_1 \in \Z}\frac{1}{\langle k_1\rangle^{1+4\underline{H}-4b-2s}}\bigg)  \bigg( \sum_{k_1' \in \Z} \frac{1}{\langle k_1' \rangle^{2\delta-2\varepsilon}}\bigg) < \infty,
\end{align*}
\normalsize
due to $1+4\underline{H}-4b-2s=1+4(\underline{H}-b-s)+2s>1$ and $2\delta -2\varepsilon > 1$ for $\varepsilon >0$ small enough. \medskip

$\bullet$ Contribution of the case  $k = 0$. Going back to \eqref{majo-imp-i}, we have
\begin{eqnarray*}
 \int d\la \, \langle \la\rangle^{2b} \, | I^{(n)}_0(\la)| 
&\lesssim&\sum_{k_1,k_1'\in \Z\backslash \{0\}}\frac{1}{\langle |k_1|\wedge |k_1'|\rangle^{4\underline{H}-1}}   \frac{1}{\langle k_1'-k_1 \rangle^{2\delta-2\varepsilon}}\\
&\lesssim&\sum_{k_1\in   \Z }\frac{1}{\langle k_1 \rangle^{4\underline{H}-1}}   \sum_{k_1'\in \Z }\frac{1}{\langle k_1'-k_1 \rangle^{2\delta-2\varepsilon}} <\infty.
\end{eqnarray*}
As a conclusion, we get \eqref{cj-fini}.

\subsubsection{Estimates related to $J^{(n)}$}\label{subse:ii-n}

 With the notation in \eqref{z-chapo}, we must prove that 
\begin{equation}\label{control2}
\sup_{n\geq 0} \big\| J^{(n)} \big\|_{\widehat{Z}^{s,b}} <\infty.
\end{equation}

 In fact, by using the same arguments as those leading to \eqref{majo-imp-i} (replace also the estimate in Proposition \ref{prop:estim-cov-1} with the one in Proposition \ref{prop:estim-cov-2}), we get that for all $k \in \Z$,
\small
\begin{align}
& \int d\la \, \langle \la\rangle^{2b} \, |J^{(n)}_k(\la) |\nonumber\\
&\lesssim   \sum_{k_1,k_1'\in \Z\backslash \{0\}} \frac{1}{\langle |k_1|\wedge |k+k_1'|\rangle^{2\underline{H}-1}}\frac{1}{\langle |k+k_1|\wedge |k_1'|\rangle^{2\underline{H}-1}} \frac{1}{\langle k_1+k_1'+k \rangle^{2\delta-2\varepsilon}}\frac{1}{\langle \Omega_{k,k_1} \rangle^{1-b}} \frac{1}{\langle \Omega_{k,k_1'} \rangle^{1-b}}.\label{majo-j-imp}
\end{align}
\normalsize
where the proportional constant no longer depends on $n$.

\smallskip

We now treat the contributions $k \neq 0$ and $k=0$ separately. \medskip

 $\bullet$ Terms with $k \neq 0$. We have
\begin{align}
&\sum_{k \neq 0} \langle k\rangle^{2s}\int d\la \, \langle \la\rangle^{2b} \, | J^{(n)}_k(\la)|\nonumber\\
&\lesssim  \sum_{k \neq 0}  \langle k\rangle^{2s}\sum_{k_1,k_1'\in \Z\backslash \{0\}} \frac{1}{\langle |k_1|\wedge |k+k_1'|\rangle^{2\underline{H}-1}}\frac{1}{\langle |k+k_1|\wedge |k_1'|\rangle^{2\underline{H}-1}} \frac{1}{\langle k_1+k_1'+k \rangle^{2\delta-2\varepsilon}}\frac{1}{\langle \Omega_{k,k_1} \rangle^{1-b}} \frac{1}{\langle \Omega_{k,k_1'} \rangle^{1-b}}.\nonumber\\
&\lesssim \sum_{k\neq 0} \frac{1}{|k|^{2-2b-2s}}\sum_{k_1,k_1'\in \Z\backslash \{0\}} \frac{1}{\langle |k_1|\wedge |k+k_1'|\rangle^{2\underline{H}-1}}\frac{1}{\langle |k+k_1|\wedge |k_1'|\rangle^{2\underline{H}-1}} \frac{1}{\langle k_1+k_1'+k \rangle^{2\delta-2\varepsilon}} \frac{1}{|k_1|^{1-b}}\frac{1}{ |k_1'|^{1-b}}\nonumber\\
&\lesssim \sum_{k\neq 0} \frac{1}{|k|^{2-2b-2s}}\sum_{k_1\neq 0}\sum_{k_1'\neq k}  \frac{1}{\langle |k_1|\wedge |k_1'|\rangle^{2\underline{H}-1}}\frac{1}{\langle |k+k_1|\wedge |k_1'-k|\rangle^{2\underline{H}-1}} \frac{1}{\langle k_1+k_1' \rangle^{2\delta-2\varepsilon}} \frac{1}{|k_1|^{1-b}}\frac{1}{ |k_1'-k|^{1-b}}\nonumber\\
&\lesssim \sum_{k_1\neq 0}\sum_{k_1'}\frac{1}{\langle |k_1|\wedge |k_1'|\rangle^{2\underline{H}-1}} \frac{1}{\langle k_1+k_1' \rangle^{2\delta-2\varepsilon}}\frac{1}{ |k_1|^{1-b}}\sum_{k\notin \{0,k_1'\}}\frac{1}{|k|^{2-2b-2s}}  \frac{1}{\langle |k+k_1|\wedge |k_1'-k|\rangle^{2\underline{H}-1}}  \frac{1}{|k_1'-k|^{1-b}}.\label{frak-j}
\end{align}
Recall that $\underline{H}\in (\frac34,1)$. For all $k_1,k'_1\in \Z$, one has
\begin{align*}
&\sum_{k\notin \{0,k_1'\}}\frac{1}{|k|^{2-2b-2s}}  \frac{1}{\langle |k_1+k|\wedge |k_1'-k|\rangle^{2\underline{H}-1}}  \frac{1}{|k_1'-k|^{1-b}}\\
&\lesssim \sum_{k}\frac{1}{\langle k\rangle^{2-2b-2s}}  \frac{1}{\langle k_1'-k\rangle^{2\underline{H}-b}}+\sum_{k}\frac{1}{\langle k\rangle^{2-2b-2s}}  \frac{1}{\langle k+k_1\rangle^{2\underline{H}-1}}  \frac{1}{\langle k_1'-k\rangle^{1-b}}.
\end{align*}
Since $(2-2b-2s)+(2\underline{H}-b)=(2-3b-2s)+2\underline{H}>1$ and $2-2b-2s<2\underline{H}-b$,
$$\sum_{k}\frac{1}{\langle k\rangle^{2-2b-2s}}  \frac{1}{\langle k_1'-k\rangle^{2\underline{H}-b}} \lesssim \frac{1}{\langle k_1'\rangle^{\nu}},$$
where $\nu:=\min(2-2b-2s,1+2\underline{H}-3b-2s)$. Besides, 
\begin{multline*}
\sum_{k}\frac{1}{\langle k\rangle^{2-2b-2s}}  \frac{1}{\langle k+k_1\rangle^{2\underline{H}-1}}  \frac{1}{\langle k_1'-k\rangle^{1-b}}\\
\begin{aligned}
&\lesssim \bigg(\sum_{k}\frac{1}{\langle k\rangle^{2-2b-2s}}  \frac{1}{\langle k+k_1\rangle^{4\underline{H}-2}}  \bigg)^{\frac12} \bigg(\sum_{k}\frac{1}{\langle k\rangle^{2-2b-2s}}  \frac{1}{\langle k_1'-k\rangle^{2-2b}}\bigg)^{\frac12}\\
&\lesssim \frac{1}{\langle k_1\rangle^{1-b-s}}\frac{1}{\langle k_1'\rangle^{\frac32-2b-s}},
\end{aligned}
\end{multline*}
where we have used the fact that $4\underline{H}-2>1$ and $(2-2b-2s)+(2-2b)>1$.

\smallskip

As a result, we have shown that for all $k_1,k_1'\in \Z$,
\begin{align*}
&\sum_{k\notin \{0,k_1'\}}\frac{1}{|k|^{2-2b-2s}}  \frac{1}{\langle |k_1+k|\wedge |k_1'-k|\rangle^{2\underline{H}-1}}  \frac{1}{|k_1'-k|^{1-b}}\lesssim \frac{1}{\langle k_1'\rangle^{\nu}}+\frac{1}{\langle k_1\rangle^{1-b-s}}\frac{1}{\langle k_1'\rangle^{\frac32-2b-s}}.
\end{align*}
By injecting this bound into \eqref{frak-j}, we deduce that in the above setting, one has
\begin{align*}
\sum_{k \neq 0} \langle k\rangle^{2s}\int d\la \, \langle \la\rangle^{2b} \, | J^{(n)}_k(\la)|&\lesssim \mathfrak{S}_1+\mathfrak{S}_2,
\end{align*}
with
$$\mathfrak{S}_1:=\sum_{k_1\neq 0}\sum_{k_1'}\frac{1}{\langle |k_1|\wedge |k_1'|\rangle^{2\underline{H}-1}} \frac{1}{\langle k_1+k_1' \rangle^{2\delta-2\varepsilon}}\frac{1}{ |k_1|^{1-b}}\frac{1}{\langle k_1'\rangle^{\nu}}$$
and
$$\mathfrak{S}_2:=\sum_{k_1\neq 0}\sum_{k_1'}\frac{1}{\langle |k_1|\wedge |k_1'|\rangle^{2\underline{H}-1}} \frac{1}{\langle k_1+k_1' \rangle^{2\delta-2\varepsilon}}\frac{1}{ |k_1|^{1-b}}\frac{1}{\langle k_1\rangle^{1-b-s}}\frac{1}{\langle k_1'\rangle^{\frac32-2b-s}}.$$
Now, on the one hand,
\small
\begin{align*}
\mathfrak{S}_1
&\lesssim\sum_{k_1'}\frac{1}{\langle k_1'\rangle^{2\underline{H}-1+\nu}} \sum_{k_1}\frac{1}{\langle k_1+k_1' \rangle^{2\delta-2\varepsilon}}\frac{1}{ \langle k_1\rangle^{1-b}}+\sum_{k_1'}\frac{1}{\langle k_1'\rangle^{\nu}} \sum_{k_1}\frac{1}{\langle k_1+k_1' \rangle^{2\delta-2\varepsilon}}\frac{1}{ \langle k_1\rangle^{2\underline{H}-b}}\\
&\lesssim\sum_{k_1'}\frac{1}{\langle k_1'\rangle^{2\underline{H}+\nu-b}}+\sum_{k_1'}\frac{1}{\langle k_1'\rangle^{ \nu+\gamma_1}} \ < \ \infty,
\end{align*}
\normalsize
due to $2\underline{H}+\nu-b=\min(2\underline{H}+2-3b-2s,1+4\underline{H}-4b-2s)>1$ for the first series and $\gamma_1:=\min(2\delta-2\eps, 2\underline{H}-b)$ so that $\nu+\gamma_1>1$ for the second one. 

\smallskip

On the other hand, for $\varepsilon>0$ small enough,
\small
\begin{align*}
\mathfrak{S}_2
&\lesssim\sum_{k_1'}\frac{1}{\langle k_1'\rangle^{2\underline{H}+\frac12-2b-s}} \sum_{k_1}\frac{1}{\langle k_1+k_1' \rangle^{2\delta-2\varepsilon}}\frac{1}{ \langle k_1\rangle^{2-2b-s}}\\
& \hspace{5cm}+\sum_{k_1'}\frac{1}{\langle k_1'\rangle^{\frac32-2b-s}} \sum_{k_1}\frac{1}{\langle k_1+k_1' \rangle^{2\delta-2\varepsilon}}\frac{1}{ \langle k_1\rangle^{2\underline{H}+1-2b-s}}\\
&\lesssim\sum_{k_1'}\frac{1}{\langle k_1'\rangle^{2\underline{H}+\frac52-4b-2s}}+\sum_{k_1'}\frac{1}{\langle k_1'\rangle^{\frac32-2b-s+\gamma_2}} \ < \ \infty,
\end{align*}
\normalsize
because $2\underline{H}+\frac52-4b-2s>1$ and $\frac32-2b-s+\gamma_2>1$ where $\gamma_2:= \min( 2\delta-2\varepsilon, 2\underline{H}+1-2b-s)$. Hence, we have proven that
$$\sup_{n\geq 0}\sum_{k \neq 0} \langle k\rangle^{2s}\int d\la \, \langle \la\rangle^{2b} \, | J^{(n)}_k(\la)| <\infty.$$
\medskip

$\bullet$ Contribution of the case  $k = 0$. With \eqref{majo-j-imp} in mind, we have

\begin{eqnarray*}
 \int d\la \, \langle \la\rangle^{2b} \, | J^{(n)}_0(\la)| 
&\lesssim&   \sum_{k_1,k_1'\in \Z\backslash \{0\}} \frac{1}{\langle |k_1|\wedge |k_1'|\rangle^{4\underline{H}-2}}  \frac{1}{\langle k_1+k_1' \rangle^{2\delta-2\varepsilon}}\\
&\lesssim&\sum_{k_1\in   \Z }\frac{1}{\langle k_1 \rangle^{4\underline{H}-2}}   \sum_{k_1'\in \Z }\frac{1}{\langle k_1+k'_1 \rangle^{2\delta-2\varepsilon}} <\infty,
\end{eqnarray*}
since $4\underline{H}-2>1$ and $2\delta -2\varepsilon>1$.
As a conclusion, we have shown that \eqref{control2} holds true.
\end{proof}

\subsection{Proof of Proposition \ref{prop:ice-cream-2}}
Just as in the proof of Proposition \ref{pro-mom}, we will only focus on the uniform estimate
\begin{equation}\label{uniff}
\sup_{n\geq 1} \big\|f^{(n)} \big\|_{Z^{s,b}} < \infty. 
\end{equation}

To this end, let us write
\begin{align*}
&\big|\cf\big(f^{(n)}_k\big)(\la)\big|=\1_{\{k\neq 0\}} \Big|\mathbb{E}\Big[ \cf \Big(\ci_\chi \widetilde{\cm}\big(\<Psi>^{(n)},\<Psi>^{(n)}\big)_k\Big)(\la)\Big]\Big|\\
&=\1_{\{k\neq 0\}} \bigg|\sum_{k_1\neq 0}\int d\la_1d\la_2\,  \ci_\chi(\la,\Omega_{k,k_1}+\la_2-\la_1) \mathbb{E}\bigg[ \overline{\cf \<Psi>^{(n)}_{k_1}(\la_1)} \cf\<Psi>^{(n)}_{k+k_1}(\la_2)\bigg]\bigg|\\
&\lesssim \1_{\{k\neq 0\}}\frac{1}{\langle \la \rangle} \sum_{k_1\neq 0}\int \frac{d\la_1 d\la_2}{\langle \la_1 \rangle\langle \la_2 \rangle} \frac{1}{\langle \la-\Omega_{k,k_1}-\la_2+\la_1\rangle}\Gamma^{(H_0)}_{k_1,k+k_1}(\la_1,\la_2) \Lambda^{(H_1)}_{k_1,k+k_1}.
\end{align*}
Therefore by \eqref{ineq-58} (with $\delta:=\min(2H_1,1)>\frac12$),
\small
\begin{align*}
&\sum_k \langle k\rangle^{2s}\int d\la \, \langle \la\rangle^{2b} \big|\cf\big(f^{(n)}_k\big)(\la)\big|^2\\
&\lesssim \sum_{k\neq 0} \langle k\rangle^{2s} \sum_{k_1,k_1'\neq 0}\int \frac{d\la_1 d\la_2}{\langle \la_1 \rangle\langle \la_2 \rangle} \int \frac{d\la'_1 d\la'_2}{\langle \la'_1 \rangle\langle \la'_2 \rangle}\bigg(\int \frac{d\la}{\langle \la\rangle^{2-2b}}  \frac{1}{\langle \la-\Omega_{k,k_1}-\la_2+\la_1\rangle}\frac{1}{\langle \la-\Omega_{k,k'_1}-\la'_2+\la'_1\rangle}\bigg)\\
&\hspace{4cm}\Gamma^{(H_0)}_{k_1,k+k_1}(\la_1,\la_2) \Lambda^{(H_1)}_{k_1,k+k_1}\Gamma^{(H_0)}_{k'_1,k+k'_1}(\la'_1,\la'_2) \Lambda^{(H_1)}_{k'_1,k+k'_1}\\
&\lesssim \sum_{k\neq 0} \langle k\rangle^{2s} \bigg|\sum_{k_1\neq 0}\Lambda^{(H_1)}_{k_1,k+k_1} \int \frac{d\la_1 d\la_2}{\langle \la_1 \rangle\langle \la_2 \rangle}  \frac{1}{\langle \Omega_{k,k_1}+\la_2-\la_1\rangle^{1-b}}\Gamma^{(H_0)}_{k_1,k+k_1}(\la_1,\la_2) \bigg|^2\\
&\lesssim \sum_{k\neq 0} \frac{\langle k\rangle^{2s}}{\langle k\rangle^{2\delta-2\varepsilon}} \bigg|\sum_{k_1\neq 0}\frac{1}{\langle |k_1|\wedge |k+k_1|\rangle^{2H_1-1}}\frac{1}{\langle \Omega_{k,k_1}\rangle^{1-b}} \int \frac{d\la_1 d\la_2}{\langle \la_1 \rangle^b\langle \la_2 \rangle^b}  \Gamma^{(H_0)}_{k_1,k+k_1}(\la_1,\la_2) \bigg|^2\\
&\lesssim \sum_{k\neq 0} \frac{\langle k\rangle^{2s}}{\langle k\rangle^{2\delta-2\varepsilon}\langle k\rangle^{2-2b}} \bigg|\sum_{k_1\neq 0}\frac{1}{\langle |k_1|\wedge |k+k_1|\rangle^{2\underline{H}-1}}\frac{1}{\langle k_1\rangle^{1-b}}  \bigg|^2,
\end{align*}
\normalsize
where we have used Lemma \ref{lem:ajout} and Lemma \ref{lem:lau} to derive the third and fourth inequalities. Due to condition \eqref{conditi}, we can easily guarantee that for $\varepsilon >0$ small enough,
$$(1-b)+(2\underline{H}-1)>1 \quad \text{and}\quad  (2\delta-2\varepsilon)+(2-2b)-2s>1,$$
which achieves to prove \eqref{uniff}.

\subsubsection{An auxiliary lemma}

The following technical property has been used in the proof of both Proposition \ref{prop:ice-cream-1} and Proposition \ref{prop:ice-cream-2}.

\begin{lemma}\label{lem:lau}
 Let $\frac12 < b<1$ and $\frac12<H_0<1$. Then  it holds that for all $k,k' \in \Z$
\begin{equation*}
\int \frac{d\la d\la' }{\langle \la\rangle^b  \langle \la'\rangle^b }\Big|\gga^{(H_0)}_{k,k'}(\la,\la')\Big| \lesssim 
\frac{1}{\langle |k|\wedge |k'|\rangle^{2(2H_0-1)}}.
\end{equation*}
\end{lemma}

\begin{proof}
For every $\varepsilon >0$, one has
\begin{align*}
&\int \frac{d\la d\la' }{\langle \la\rangle^b  \langle \la'\rangle^b }\Big|\gga^{(H_0)}_{k,k'}(\la,\la')\Big|\nonumber\\
&\lesssim \int \frac{d\la d\la' }{\langle \la\rangle^b  \langle \la'\rangle^b } \bigg( \frac{1}{\langle ||k|^2-\la| \wedge ||k'|^2-\la'|\rangle^{2H_0-1}}\frac{1}{\langle |k|^2-\la-|k'|^2+\la'\rangle^{1-\varepsilon}}\bigg)\nonumber\\
&= \int \frac{d\la d\la' }{\langle |k|^2-\la\rangle^b  \langle |k'|^2-\la'\rangle^b } \bigg( \frac{1}{\langle |\la| \wedge |\la'|\rangle^{2H_0-1}}\frac{1}{\langle \la-\la'\rangle^{1-\varepsilon}}\bigg)\label{gain-k-inter}\\
&\lesssim \int_{\{|\la|<|\la'|\}} \frac{d\la d\la' }{\langle |k|^2-\la\rangle^b  \langle |k'|^2-\la'\rangle^b } \bigg( \frac{1}{\langle \la \rangle^{2H_0-1}}\frac{1}{\langle \la-\la'\rangle^{1-\varepsilon}}\bigg)\nonumber\\
&\hspace{1cm}+\int_{\{|\la'|<|\la|\}} \frac{d\la d\la' }{\langle |k|^2-\la\rangle^b  \langle |k'|^2-\la'\rangle^b } \bigg( \frac{1}{\langle \la' \rangle^{2H_0-1}}\frac{1}{\langle \la-\la'\rangle^{1-\varepsilon}}\bigg)\nonumber\\
&\lesssim \int\frac{d\la  }{\langle \la \rangle^{2H_0-1}  } \frac{1}{\langle |k|^2-\la\rangle^b }\int\frac{d\la' }{\langle |k'|^2-\la'\rangle^b } \frac{1}{\langle \la-\la'\rangle^{1-\varepsilon}}\nonumber
. 
\end{align*}
Now we apply Lemma~\ref{lem:tec-1d} and obtain
\begin{multline*}
\int\frac{d\la  }{\langle \la \rangle^{2H_0-1}  } \frac{1}{\langle |k|^2-\la\rangle^b }\int\frac{d\la' }{\langle |k'|^2-\la'\rangle^b } \frac{1}{\langle \la-\la'\rangle^{1-\varepsilon}} \\
\begin{aligned}
&\lesssim \int\frac{d\la  }{|\la |^{2H_0-1}  } \frac{1}{\langle |k|^2-\la\rangle^b } \frac{1}{\langle |k'|^2-\la\rangle^{b-\eps} }\\
&\lesssim \bigg(\int\frac{d\la  }{|\la |^{2H_0-1}  } \frac{1}{\langle |k|^2-\la\rangle^{2b} }\bigg)^{1/2}\bigg(\int\frac{d\la  }{|\la |^{2H_0-1}  } \frac{1}{\langle |k'|^2-\la\rangle^{2b-2\eps} }\bigg)^{1/2}.
\end{aligned}
\end{multline*}
Finally, we can apply Lemma \ref{lem:int-singu-zero} twice to get the desired estimate, provided that $\eps>0$ is chosen small enough such that $2b-2\eps>1$. 
\end{proof}

\section{Estimate on the random operator}\label{sec:estim-random-op}

It remains us to proceed with the construction and the control of the two product terms
\begin{equation*}
\cl^{\<Psi>,+}_\chi( z):= \ci_\chi \cm(z,\<Psi>) \quad \text{and} \quad \cl^{\<Psi>,-}_\chi(z):= \ci_\chi \cm( \<Psi>,z)\, ,
\end{equation*} 
involved in the remainder equation \eqref{equa-z-reno}. In other words, we will here develop the arguments of the proof of Proposition \ref{prop:control-op-l-intro}

\smallskip

To be more specific, we will focus on the convergence of $\chi_\tau\cdot \cl^{(n),+}$, but the convergence of $\chi_\tau\cdot \cl^{(n),-}$ could be derived from completely similar arguments.

\smallskip

Let us decompose the operator as
\small
\begin{align}
\cl^{(n),+}(z)_k= \ci_\chi \cm( z,\<Psi>^{(n)})_k&=\ci_{\chi}\bigg(\sum_{k_1} e^{\imath .\Omega_{k,k_1}} z_{k+k_1}(.) \overline{\<Psi>^{(n)}_{k_1}(.)}\bigg)=\cl^{\circ,(n)}(z)_k+{\cl}^{\#,(n)}(z)_k,\label{decompo-l-n-plus}
\end{align}
\normalsize
with
\begin{equation*}
\cl^{\circ,(n)}(z)_k:= \ci_{\chi}\Big( z_{k}(.) \overline{\<Psi>^{(n)}_{0}(.)}\Big), \quad 
\cl^{\#,(n)}(z)_k:= \ci_{\chi}\bigg(\sum_{k_1\neq 0} e^{\imath .\Omega_{k,k_1}} z_{k+k_1}(.) \overline{\<Psi>^{(n)}_{k_1}(.)}\bigg).
\end{equation*}
Following the convention initiated in Section \ref{subsec:cova}, we also set 
$$\cl^{(n,n+1),+}:=\cl^{(n+1),+}-\cl^{(n),+}, \quad \cl^{\circ,(n,n+1)}:=\cl^{\circ,(n+1)}-\cl^{\circ,(n)},$$
$$ {\cl}^{\#,(n,n+1)}:={\cl}^{\#,(n+1)}- {\cl}^{\#,(n)},$$
and recall that $\<Psi>^{(n,n+1)}:=\<Psi>^{(n+1)}-\<Psi>^{(n)}$.

\smallskip

For a clear presentation of our arguments, we will estimate the two operators $\cl^{\circ,(n,n+1)}$  and ${\cl}^{\#,(n,n+1)}$ separately, and then combine the results towards the desired proof (see Section \ref{subsec:proof-prop-produ}).

\subsection{Estimate for $\cl^{\circ,(n,n+1)}$}\label{subsec:flat}

We are  looking for an estimate of $\big\|\cl^{\circ,(n,n+1)}\big\|_{\cl(Z^{s,b},Z^{s,b'})}$ for $b'>b$.

\begin{proposition}\label{prop:mom-q-flat}

Fix $\frac12<H_0<1$ and $0<H_1<1$. 

\smallskip

\noindent
$(i)$ For all $0<b,b'<1$  and $s\in \R$, one has
\begin{align*}
&\big\|\cl^{\circ,(n,n+1)}(z)\big\|_{Z^{s,b'}}\lesssim  \big(\cq^{\circ,(n)}_{b,b'}\big)^{\frac12} \|z\|_{Z^{s,b}},
\end{align*}
with
$$\cq^{\circ,(n)}_{b,b'}:=\int d\la \, \langle \la\rangle^{2b'}  \int \frac{d\beta}{\langle \beta\rangle^{2b}}\bigg|\int d\la_1 \, \overline{\ci_\chi(\la,\la_1)}  \cf\big(\<Psi>^{(n,n+1)}_0 \big)(\la_1-\beta) \bigg|^2.$$
\noindent
$(ii)$ For all $\frac12<b,b'<1$ and $0<\ka\leq \frac12 \min(H_0,H_1)$, it holds that
\begin{equation*}
\mathbb{E}\Big[\big|\cq^{\circ,(n)}_{b,b'}\big|\Big]\lesssim 2^{-2\ka n}.
\end{equation*}
\end{proposition}

\begin{proof}

$(i)$ 
One has
\small
\begin{align*}
&\big\|\cl^{\circ,(n,n+1)}(z)\big\|_{Z^{s,b'}}^2=\sum_{k}\langle k\rangle^{2s}\int d\la \, \langle \la\rangle^{2b'} \bigg| \int d\la_1 \, \ci_\chi(\la,\la_1) \cf\Big( z_{k}(.) \overline{\<Psi>^{(n,n+1)}_{0}(.)}\Big)(\la_1)\bigg|^2 \\
&=\sum_{k}\langle k\rangle^{2s}\int d\la \, \langle \la\rangle^{2b'} \bigg| \int d\la_1 \, \ci_\chi(\la,\la_1) \int d\beta\, \cf\Big(\overline{\<Psi>^{(n,n+1)}_0} \Big)(\la_1-\beta) \cf\big(z_{k}\big)(\beta)\bigg|^2\\
&=\sum_{k}\langle k\rangle^{2s}\int d\la \, \langle \la\rangle^{2b'} \bigg| \int d\beta\, \Big[ \langle \beta\rangle^{b}\cf\big(z_{k}\big)(\beta)\Big] \bigg[ \frac{1}{\langle \beta\rangle^{b}}\int d\la_1 \, \ci_\chi(\la,\la_1)  \cf\Big(\overline{\<Psi>^{(n,n+1)}_0} \Big)(\la_1-\beta)\bigg] \bigg|^2 \\
&\leq \bigg( \sum_k \langle k\rangle^{2s}\int d\beta\,  \langle \beta\rangle^{2b}\big|\cf\big(z_{k}\big)(\beta)\big|^2\bigg) \bigg( \int d\la \, \langle \la\rangle^{2b'}  \int \frac{d\beta}{\langle \beta\rangle^{2b}}\bigg|\int d\la_1 \, \ci_\chi(\la,\la_1)  \cf\Big(\overline{\<Psi>^{(n,n+1)}_0} \Big)(\la_1-\beta)\bigg|^2\bigg)\, ,
\end{align*}
\normalsize
which corresponds to the desired inequality.

\smallskip

\noindent
$(ii)$ 
Let us write
\begin{align*}
&\mathbb{E}\Big[\big|\cq^{\circ,(n)}_{b,b'}\big|\Big]\\
&= \int d\la \, \langle \la\rangle^{2b'}  \int \frac{d\beta}{\langle \beta\rangle^{2b}}\mathbb{E}\bigg[\bigg|\int d\la_1 \, \overline{\ci_\chi(\la,\la_1-\beta)}  \cf\big(\<Psi>^{(n,n+1)}_0 \big)(\la_1) \bigg|^2\bigg]\\
&= \int d\la \, \langle \la\rangle^{2b'}  \int \frac{d\beta}{\langle \beta\rangle^{2b}}\int d\la_1 \int d\la_1' \, \overline{\ci_\chi(\la,\la_1-\beta)} \ci_\chi(\la,\la_1'-\beta) \mathbb{E}\Big[ \cf\big(\<Psi>^{(n,n+1)}_0 \big)(\la_1) \overline{\cf\big(\<Psi>^{(n,n+1)}_0 \big)(\la_1')} \Big],
\end{align*}
and then combine \eqref{bou-ci-chi} with the covariance estimate \eqref{bou:cova-luxo-n-n+1} to obtain that
\begin{align*}
&\mathbb{E}\Big[\big|\cq^{\circ,(n)}_{b,b'}\big|\Big]\lesssim 2^{-2n\ka}\int \frac{d\la}{\langle \la\rangle^{2-2b'}}  \int \frac{d\beta}{\langle \beta\rangle^{2b}}\\
&\hspace{2cm}\int \frac{d\la_1}{\langle\la_1\rangle} \int \frac{d\la_1'}{\langle \la_1'\rangle}\frac{1}{\langle \la-(\la_1-\beta)\rangle}\frac{1}{\langle \la-(\la_1'-\beta)\rangle}\Big[\gga^{(H_0-\ka)}_{0,0}(\la_1,\la_1')+\gga^{(H_0)}_{0,0}(\la_1,\la_1')\Big],
\end{align*}
where we have used the fact that $\Lambda^{(H_{1})}_{0,0},\Lambda^{(H_1-\ka)}_{0,0}\lesssim 1$ (by Lemma \ref{lem:int-singu-zero}, item $(iii)$). 

\smallskip

By applying the Cauchy-Schwarz inequality with respect to the pair $(\la_1,\la_1')$, we deduce
\begin{align}
&\mathbb{E}\Big[\big|\cq^{\circ,(n)}_{b,b'}\big|\Big]\lesssim 2^{-2n\ka}\bigg[\int   \frac{d\la}{\langle \la\rangle^{2-2b'}}  \int \frac{d\beta}{\langle \beta\rangle^{2b}}\int \frac{d\la_1}{\langle\la_1\rangle \langle \la-(\la_1-\beta)\rangle^2}\bigg]\nonumber\\
&\hspace{5cm}\bigg[\int \frac{d\la_1}{\langle\la_1\rangle} \int \frac{d\la_1'}{\langle \la_1'\rangle}\big|\gga^{(H_0-\ka)}_{0,0}(\la_1,\la_1')+\gga^{(H_0)}_{0,0}(\la_1,\la_1')\big|^2\bigg]^{\frac12}.\label{flat-proof}
\end{align}

\smallskip

To ensure finiteness of the first term into brackets, it suffices to observe that
\begin{align*}
\int   \frac{d\la}{\langle \la\rangle^{2-2b'}}  \int \frac{d\beta}{\langle \beta\rangle^{2b}}\int \frac{d\la_1}{\langle\la_1\rangle \langle \la-(\la_1-\beta)\rangle^2}&\lesssim \int   \frac{d\la}{\langle \la\rangle^{2-2b'}}  \int \frac{d\beta}{\langle \beta\rangle^{2b} \langle \beta+\la\rangle} \lesssim \int   \frac{d\la}{\langle \la\rangle^{3-2b'}},
\end{align*}
and since $b'<1$, the latter integral is clearly finite.

\

As for the second term into brackets in \eqref{flat-proof}, we can apply  Lemma \ref{lem:tec-1d} (recall that $H_0>\frac12$) to assert that
\begin{align*}
\int \frac{d\la_1}{\langle\la_1\rangle} \int \frac{d\la_1'}{\langle \la_1'\rangle}\big|\gga^{(H_0-\ka)}_{0,0}(\la_1,\la_1')+\gga^{(H_0)}_{0,0}(\la_1,\la_1')\big|^2&\lesssim \int \frac{d\la_1}{\langle\la_1\rangle} \int \frac{d\la_1'}{\langle \la_1'\rangle}\frac{1}{\langle \la_1-\la_1'\rangle^{2-2\varepsilon}}\\
&\lesssim \int \frac{d\la_1}{\langle\la_1\rangle^2}<\infty.
\end{align*}
\end{proof}

\subsection{Estimate for ${\cl}^{\#,(n,n+1)}$}\label{subsec:tilde}

Given the higher sophistication of ${\cl}^{\#,(n)}$ (in comparison with  $\cl^{\circ,(n)}$), let us introduce an additional notation for the Fourier kernel associated with this operator. Namely, with the notation of Lemma \ref{lem:integr-kernel}, we can write
\small
\begin{align*}
&\cf\big( {\cl}^{\#,(n)}( z)_k\big)(\la)= \cf\bigg(\ci_{\chi}\bigg(\sum_{k_1\neq 0} e^{\imath .\Omega_{k,k_1}} z_{k+k_1}(.) \overline{\<Psi>^{(n)}_{k_1}(.)}\bigg)\bigg)(\la)\\
&=\ \sum_{k_1}\1_{\{k_1\neq 0\}} \int_{\R} d\la' \, \ci_\chi(\la,\la') \int dt\, e^{-\imath \la' t}e^{\imath t \Omega_{k,k_1}} z_{k+k_1}(t) \overline{\<Psi>^{(n)}_{k_1}(t)}\\
&=\sum_{k_1} \1_{\{k_1\neq 0\}} \int d\la_1 \, \cf(z_{k+k_1})(\la_1) \int d\la_2 \, \overline{\cf\big(\<Psi>^{(n)}_{k_1}\big)(\la_2)} \int_{\R} d\la' \, \ci_\chi(\la,\la')\int dt\, e^{-\imath \la' t}e^{\imath t\Omega_{k,k_1}} e^{\imath t\la_1}e^{-\imath t \la_2} \\
&=\sum_{k_1}\int d\la_1 \, \cf(z_{k+k_1})(\la_1) \bigg[\1_{\{k_1\neq 0\}} \int d\la_2 \, \overline{\cf\big(\<Psi>^{(n)}_{k_1}\big)(\la_2)} \ci_\chi(\la,\Omega_{k,k_1}+\la_1-\la_2)\bigg]\\
&=\sum_{k_1}\int d\la_1 \, \cf(z_{k_1})(\la_1) \bigg[\1_{\{k_1\neq k\}} \int d\la_2 \, \overline{\cf\big(\<Psi>^{(n)}_{k_1-k}\big)(\la_2)} \ci_\chi(\la,\Omega_{k,k_1-k}+\la_1-\la_2)\bigg],
\end{align*}
\normalsize
that is
\begin{equation*}
\cf\big( {\cl}^{\#,(n)}( z)_k\big)(\la)=\sum_{k_1}\int d\la_1 \big(\ck^{(n)}_\chi\big)_{kk_1}(\la,\la_1) \cf(z_{k_1})(\la_1),
\end{equation*}
where we set from now on
$$\big(\ck^{(n)}_\chi\big)_{kk_1}(\la,\la_1):=\1_{\{k_1\neq  k\}}\int d\la_2 \, \overline{\cf\big(\<Psi>^{(n)}_{k_1-k}\big)(\la_2)}\ci_\chi(\la,\Omega_{k,k_1-k}+\la_1-\la_2).$$
Let us also set 
\begin{align}
\big(\ck^{(n,n+1)}_\chi\big)_{kk_1}(\la,\la_1):= \1_{\{k_1\neq  k\}}\int d\la_2 \, \overline{\cf\big(\<Psi>^{(n,n+1)}_{k_1-k}\big)(\la_2)} \ci_\chi(\la,\Omega_{k,k_1-k}+\la_1-\la_2)\label{k-n-n-1}
\end{align}
where we recall that $\<Psi>^{(n,n+1)}:=\<Psi>^{(n+1)}-\<Psi>^{(n)}$. Accordingly, it holds that
\begin{equation}\label{kernel-random-op}
\cf\big( {\cl}^{\#,(n,n+1)}( z)_k\big)(\la)=\sum_{k_1}\int d\la_1 \big(\ck^{(n,n+1)}_\chi\big)_{kk_1}(\la,\la_1) \cf(z_{k_1})(\la_1).
\end{equation}

Note that we shall later deduce the desired estimate on $\big\|\chi_\tau\cdot {\cl}^{\#,(n,n+1)}\big\|_{\cl(Z^{s,b},Z^{s,b})}$  (for some $b>\frac12$) from a suitable interpolation between 
$$\big\|\chi_\tau\cdot {\cl}^{\#,(n,n+1)}\big\|_{\cl(Z^{s,b},Z^{s,1})} \quad \text{and} \quad \big\|\chi_\tau\cdot {\cl}^{\#,(n,n+1)}\big\|_{\cl(Z^{s,b},Z^{s,a})} \ \text{for} \ 0<a<\frac12,$$
hence our subsequent investigations about these two quantities.

\begin{proposition}\label{prop:q-tilde}

Fix $\frac12<H_0<1$ and $0<H_1<1$ such that $2H_0+H_1>\frac32$. Set $\underline{H}:=2H_0+H_1-1$ and recall that the random kernel $\ck^{(n,n+1)}_\chi$ has been introduced in~\eqref{k-n-n-1}.

\smallskip

\noindent
$(i)$ For all $a,b,s\in \R$, it holds that
$$\big\|{\cl}^{\#,(n,n+1)} (z)\big\|_{Z^{s,a}} \lesssim  ({\cq}^{\#,(n)}_{a,b,s})^{\frac12} \cdot \|z\|_{Z^{s,b}},$$
where the random quantity ${\cq}^{\#,(n)}_{a,b,s}$ is given by
\small
\begin{align*}
& {\cq}^{\#,(n)}_{a,b,s}:=\nonumber\\
&\bigg( \sum_{k_1,k_1'}\frac{1}{\langle k_1\rangle^{2s} \langle k'_1\rangle^{2s}} \int \frac{d\la_1 d\la_1'}{\langle \la_1\rangle^{2b} \langle \la_1'\rangle^{2b}} \bigg| \sum_{k} \langle k\rangle^{2s}\int d\la \, \langle \la\rangle^{2a}\big(\ck^{(n,n+1)}_\chi\big)_{kk_1}(\la,\la_1) \overline{\big(\ck^{(n,n+1)}_\chi\big)_{kk_1'}(\la,\la_1')}\bigg|^2\bigg)^{\frac12}.
\end{align*}
\normalsize

\smallskip

\noindent
$(ii)$ For all $0<a<\frac12<b<1$  and  $1-\underline{H}<s<\frac12$, one has
\begin{equation*}
\mathbb{E}\Big[ \big| {\cq}^{\#,(n)}_{a,b,s}\big|^{2}\Big] \lesssim 2^{-4n\ka_s},
\end{equation*} 
where we have set $\ka_s:=\frac12 \min\big(H_0,H_1, \frac{s-(1-\underline{H})}{2}\big)>0$.
\end{proposition}

\begin{proof}

$(i)$ Using the representation \eqref{kernel-random-op}, we can write
\small
\begin{align*}
&\big\|{\cl}^{\#,(n,n+1)} (z)\big\|_{Z^{s,a}}^2=\int d\la \, \langle \la\rangle^{2a}\sum_k \langle k\rangle^{2s}\bigg|\sum_{k_1}\int d\la_1 \big(\ck^{(n,n+1)}_\chi\big)_{kk_1}(\la,\la_1) \cf(z_{k_1})(\la_1) \bigg|^2\\
&=\sum_{k_1,k_1'}\int d\la_1 d\la_1' \, \bigg(\big[\langle k_1\rangle^{s} \langle \la_1\rangle^b \cf(z_{k_1})(\la_1)\big]\big[\langle k_1'\rangle^{s}\langle \la'_1\rangle^b \overline{\cf(z_{k'_1})(\la'_1)}\big]\bigg)\\
&\hspace{0.5cm}\bigg(\langle k_1\rangle^{-s}\langle k_1'\rangle^{-s} \langle \la_1\rangle^{-b}\langle \la'_1\rangle^{-b}\sum_k  \langle k\rangle^{2s}\int d\la \, \langle \la\rangle^{2a}\big(\ck^{(n,n+1)}_\chi\big)_{kk_1}(\la,\la_1) \overline{\big(\ck^{(n,n+1)}_\chi\big)_{kk'_1}(\la,\la'_1)}\bigg),
\end{align*}
\normalsize
and we immediately derive  the desired estimate by applying the Cauchy-Schwarz inequality with respect to the $4$-uplet $(k_1,k_1',\la_1,\la_1')$.

\smallskip

\noindent
$(ii)$ For the sake of clarity, we have reported this moment estimate in Section \ref{sec:appendix} below.
\end{proof}

\begin{proposition}\label{prop:s-tilde}

Fix  $\frac12 <H_0<1$ and $0< H_1<1$ such that $\frac32 <2H_0+H_1<2$. 

\smallskip

\noindent
$(i)$ For every $0\leq s<\frac12$, it holds that
$$\big\|{\cl}^{\#,(n,n+1)} (z)\big\|_{Z^{s,1}} \lesssim \big\{ ({\cs}^{(n)})^{\frac12}+( {\cs}^{(n+1)})^{\frac12}\big\} \cdot \|z\|_{Z^{s,0}},$$
where the random quantity $ {\cs}^{(n)}$ is given by
$$ {\cs}^{(n)}:=\sup_{t\in \R, k\in \Z}\big|\langle k\rangle^2\, \<Psi>^{(n)}_{k}(t)\big|^2.$$

\smallskip

\noindent
$(ii)$ It holds that 
\begin{equation*} 
\mathbb{E}\Big[ \big| {\cs}^{(n)}\big|\Big] \lesssim 2^{7n}.
\end{equation*} 
\end{proposition}

\begin{proof}

$(i)$ Using a standard Sobolev spaces identification, we can first assert that 
\begin{align*}
\big\|{\cl}^{\#,(n)} (z)\big\|_{Z^{s,1}}^2&= \sum_k \langle k\rangle^{2s}\int d\la \, \langle \la\rangle^2 \big| \cf\big({\cl}^{\#,(n)} (z)_k\big)(\la)\big|^2\\
&\lesssim \sum_k \langle k\rangle^{2s}\Big\{ \big\| {\cl}^{\#,(n)} (z)_k\big\|_{L^2(\R)}^2+\big\| \partial_t\big({\cl}^{\#,(n)} (z)_k\big)\big)\big\|_{L^2(\R)}^2\Big\}.
\end{align*}
Recall that for any $f:\R\to\R$, $(\ci_\chi f)(t)=-i\chi(t)\int_0^t ds\, \chi(s) f(s)$, and so one has clearly
\begin{equation*}
\big\| \ci_\chi f\big\|_{L^2(\R)}^2+\big\| \partial_t\big(\ci_\chi f\big)\big\|_{L^2(\R)}^2 \lesssim \|f\|_{L^2(\R)}^2,
\end{equation*}
which yields here
\begin{align*}
\big\|{\cl}^{\#,(n)} (z)\big\|_{Z^{s,1}}^2&\lesssim \sum_{k} \langle k\rangle^{2s}\bigg\| \sum_{k_1\neq 0} e^{\imath .\Omega_{k,k_1}} z_{k+k_1}(.) \overline{\<Psi>^{(n)}_{k_1}(.)}\bigg\|_{L^2(\R)}^2\\
&\lesssim \int dt\sum_k \langle k\rangle^{2s}\bigg|\sum_{k_1} \big|\<Psi>^{(n)}_{k_1-k}(t)\big|\big| z_{k_1}(t)\big| \bigg|^2\\
&\lesssim {\cs}^{(n)}\int dt\sum_k \langle k\rangle^{2s}\bigg|\sum_{k_1} \frac{1}{\langle k_1-k\rangle^2}\frac{1}{\langle k_1\rangle^s}\big|\langle k_1\rangle^s z_{k_1}(t)\big| \bigg|^2\\
&\lesssim {\cs}^{(n)}\int dt\sum_k \langle k\rangle^{2s}\bigg(\sum_{k'_1} \frac{1}{\langle k'_1-k\rangle^2}\frac{1}{\langle k'_1\rangle^{2s}}\bigg) \bigg(\sum_{k_1} \frac{1}{\langle k_1-k\rangle^2}\big|\langle k_1\rangle^s z_{k_1}(t)\big|^2 \bigg)\\
&\lesssim  {\cs}^{(n)}\int dt\sum_k \sum_{k_1} \frac{1}{\langle k_1-k\rangle^2}\big|\langle k_1\rangle^s z_{k_1}(t)\big|^2 \lesssim {\cs}^{(n)} \, \|z\|_{Z^{s,0}}^2.
\end{align*}

\smallskip

\noindent
$(ii)$ Recall that
\begin{equation*}
\<Psi>^{(n)}_k(t)=- \imath \chi(t)\int_0^t dr \, e^{-\imath r |k|^2} \big[\chi(r) \dot{B}_k^{(n)}(r)\big],
\end{equation*}
and so, with an elementary integration-by-parts argument, we deduce that for all $(t,k)\in \R\times \Z$,
\begin{align*}
\langle k\rangle^2 \, \big|\<Psi>^{(n)}_k(t)\big|&\lesssim \big| \dot{B}^{(n)}_k(0)\big| +\int_0^2 dr\, \Big\{\big| \dot{B}^{(n)}_k(r)\big|+\big| (\partial_r\dot{B}^{(n)}_k)(r)\big|\Big\}\\
&\lesssim \int_{\T}dx\, \big| \dot{B}^{(n)}(0,x)\big| +\int_{\T} dx\int_0^2 dr\, \Big\{\big| \dot{B}^{(n)}(r,x)\big|+\big| (\partial_r\dot{B}^{(n)})(r,x)\big|\Big\}.
\end{align*}
This immediately entails that
$${\cs}^{(n)} \lesssim \int_{\T}dx\, \big| \dot{B}^{(n)}(0,x)\big|^2 +\int_{\T} dx\int_0^2 dr\, \Big\{\big| \dot{B}^{(n)}(r,x)\big|^2+\big| (\partial_r\dot{B}^{(n)})(r,x)\big|^2\Big\},$$
and accordingly
\begin{multline*}
\mathbb{E}\big[{\cs}^{(n)}\big] \lesssim \\
\lesssim \int_{\T}dx\, \mathbb{E}\Big[\big| \dot{B}^{(n)}(0,x)\big|^2\Big] +\int_{\T} dx\int_0^2 dr\, \Big\{\mathbb{E}\Big[\big| \dot{B}^{(n)}(r,x)\big|^2\Big]+\mathbb{E}\Big[\big| (\partial_r\dot{B}^{(n)})(r,x)\big|^2\Big]\Big\}.
\end{multline*}
The moments in the above right-hand side can now be readily estimated by using the representation \eqref{defi:approx-b-n} of $\dot{B}^{(n)}$. For instance, one has
$$
(\partial_r \dot{B}^{(n)})(r,x)= -\imath\, c_H\int_{|\xi|\leq 2^{2n}}^{}\int_{|\eta|\leq 2^n} \frac{\xi^2\, e^{\imath r\xi}}{|\xi|^{H_0+\frac{1}{2}}}\frac{\eta\, e^{\imath x\eta}}{|\eta|^{H_1+\frac{1}{2}}}\, \widehat{W}(d\xi,d\eta) 
$$
and accordingly
\begin{align*}
&\mathbb{E}\Big[\big| (\partial_r\dot{B}^{(n)})(r,x)\big|^2\Big]\lesssim \bigg(\int_{|\xi|\leq 2^{2n}} \frac{d\xi}{|\xi|^{2H_0-3}}\bigg)\bigg(\int_{|\eta|\leq 2^{n}} \frac{d\eta}{|\eta|^{2H_1-1}}\bigg) \lesssim 2^{2n(5-(2H_0+H_1))}\lesssim 2^{7n},
\end{align*}
which was the claim.
\end{proof}

\subsection{Proof of Proposition \ref{prop:control-op-l-intro}}\label{subsec:proof-prop-produ}

Recall that we focus on the convergence of the sequence $\big(\cl^{(n),+}\big)_{n\geq 1}$ in $\mathfrak{L}_\mu(Z^{s,b},Z^{s,b})$. This property will be deduced from a suitable combination of the estimates exhibited in the previous sections.

\smallskip

For more clarity, let us organize the arguments along three successive steps.

\smallskip

\noindent
\textit{Step 1: A first general bound.} Starting from the decomposition of $\cl^{(n),+}$ in \eqref{decompo-l-n-plus}, we can obviously write
\small
\begin{equation*}
\big\|\chi_\tau\cdot \cl^{(n,n+1),+}\big\|_{\cl(Z^{s,b},Z^{s,b})}\leq  \big\|\chi_\tau\cdot \cl^{\circ,(n,n+1)}\big\|_{\cl(Z^{s,b},Z^{s,b})}+\big\|\chi_\tau\cdot {\cl}^{\#,(n,n+1)}\big\|_{\cl(Z^{s,b},Z^{s,b})}.
\end{equation*}
\normalsize
Then, by Lemma \ref{b-prim}, we obtain that for all $0<s<\frac12<b<b'<1$,

\begin{eqnarray*}
\big\|\chi_\tau\cdot \cl^{(n,n+1),+}\big\|_{\cl(Z^{s,b},Z^{s,b})} &\lesssim & 
 \tau^{b'-b}\Big[\big\|\cl^{\circ,(n,n+1)}\big\|_{\cl(Z^{s,b},Z^{s,b'})}+\big\|{\cl}^{\#,(n,n+1)}\big\|_{\cl(Z^{s,b},Z^{s,b'})}\Big]\nonumber\\
&\lesssim& \tau^{b'-b}\Big[\big(\cq^{\circ,(n)}_{b,b'}\big)^{\frac12}+\big\|{\cl}^{\#,(n,n+1)}\big\|_{\cl(Z^{s,b},Z^{s,b'})}\Big],
\end{eqnarray*}
where we  have used the result (and notation) of Proposition \ref{prop:mom-q-flat} to get the second inequality.

\smallskip

In other words, recalling the notation \eqref{l-mu-x-b-c}, we have deduced that for all $0<s<\frac12<b<b'<1$, there exists $\mu>0$ such that
\begin{align}
&\big\| \cl^{(n,n+1),+}\big\|_{\mathfrak{L}_\mu(Z^{s,b},Z^{s,b})}\lesssim \big(\cq^{\circ,(n)}_{b,b'}\big)^{\frac12}+\big\|{\cl}^{\#,(n,n+1)}\big\|_{\cl(Z^{s,b},Z^{s,b'})}.\label{proof-pro}
\end{align}

\

\noindent
\textit{Step 2: Interpolation.} In order to bound the moments of the norm
$$\big\|{\cl}^{\#,(n,n+1)}\big\|_{\cl(Z^{s,b},Z^{s,b'})},$$
we shall now implement the interpolation procedure alluded to in Section \ref{subsec:tilde}.

\smallskip

To this end, recall that we have fixed $s\in \big(1-\underline{H}, \frac12\big)$, and introduce two parameters
\begin{equation*}\label{defi-ex-b}
a:=\frac12-\varepsilon\in (0,\frac12), \quad b':=\frac12+\varepsilon\in (\frac12,1),
\end{equation*}
for $\varepsilon >0$ to be determined later on.  Setting
\begin{equation}\label{theta}
\theta:=\frac{b'-a}{1-a}=\frac{2\varepsilon}{\frac12+\varepsilon}\in (0,1),
\end{equation}
we get by interpolation, for every $b\in (\frac12,b')$,
\begin{align}
&\big\|{\cl}^{\#,(n,n+1)}\big\|_{\cl(Z^{s,b},Z^{s,b'})}\lesssim \big\|{\cl}^{\#,(n,n+1)}\big\|_{\cl(Z^{s,b},Z^{s,a})}^{1-\theta}\big\|{\cl}^{\#,(n,n+1)}\big\|_{\cl(Z^{s,b},Z^{s,1})}^{\theta}\nonumber,
\end{align}
which, with the notation of Proposition \ref{prop:q-tilde} and Proposition \ref{prop:s-tilde}, yields 
\begin{align*}
&\big\|{\cl}^{\#,(n,n+1)}\big\|_{\cl(Z^{s,b},Z^{s,b'})}\lesssim ({\cq}^{\#,(n)}_{a,b,s})^{\frac12(1-\theta)}( {\cs}^{(n)}+{\cs}^{(n+1)})^{\frac12\theta}. 
\end{align*}
As a result, for every $b\in (\frac12,b')$,
\begin{align*}
\mathbb{E}\Big[\big\|{\cl}^{\#,(n,n+1)}\big\|_{\cl(Z^{s,b},Z^{s,b'})}^2\Big] &\lesssim \mathbb{E}\Big[ \big|{\cq}^{\#,(n)}_{a,b,s}\big|^{1-\theta}\big|{\cs}^{(n)}+{\cs}^{(n+1)}\big|^{\theta}\Big]\\
&\lesssim \mathbb{E}\Big[ \big|{\cq}^{\#,(n)}_{a,b,s}\big|\Big]^{1-\theta} \mathbb{E}\Big[ \big|{\cs}^{(n)}+{\cs}^{(n+1)}\big|\Big]^{\theta}.
\end{align*}
We are here in position to apply the moments estimates established in Proposition~\ref{prop:q-tilde} and Proposition~\ref{prop:s-tilde}, which gives, for every $b\in (\frac12,b')$,
\begin{align*}
&\mathbb{E}\Big[\big\|{\cl}^{\#,(n,n+1)}\big\|_{\cl(Z^{s,b},Z^{s,b'})}^2\Big]\lesssim 2^{-2n((1-\theta)\ka_s-\frac72\theta)}
\end{align*}
where we have set $\ka_s:=\frac12 \min\big(H_0,H_1, \frac{s-(1-\underline{H})}{2}\big)>0$.

\smallskip

Then, since $\ka_s$ does not depend on $(a,b')$, or equivalently on $\varepsilon$, and using \eqref{theta}, we can find such $\varepsilon>0$ small enough (depending on $s$) so that
$$\al_s:=(1-\theta)\ka_s-\frac72\theta >0, \quad \text{that is} \ \ \big(\frac12-\varepsilon\big) \frac{1}{2\varepsilon} > \frac{7}{2\ka_s}.$$ 
For such a choice of $\varepsilon>0$ (that we fix from now on), we have thus shown that for every $b\in (\frac12,b')$,
\begin{equation}\label{step-2-estima}
\mathbb{E}\Big[\big\|{\cl}^{\#,(n,n+1)}\big\|_{\cl(Z^{s,b},Z^{s,b'})}^2\Big]\lesssim 2^{-2n\al_s},
\end{equation}
with $\al_s>0$.

\smallskip

\noindent
\textit{Step 3: Conclusion.} Let us fix $s$ and $b'$ just as in the above Step 2, and set $b^\ast_s:=\frac12 \big(\frac12+b')\in (\frac12,b')$. By combining \eqref{proof-pro} and \eqref{step-2-estima}, we first derive that for every $b\in (\frac12,b^\ast_s]$,
\small
\begin{align*}
&\mathbb{E}\Big[\big\| \cl^{(n,n+1),+}\big\|_{\mathfrak{L}_\mu(Z^{s,b},Z^{s,b})}^2\Big]\lesssim 2^{-2n\al_s}+\mathbb{E}\Big[\big|\cq^{\circ,(n)}_{b,b'}\big|\Big],
\end{align*}
\normalsize
for some $\mu\in (0,1)$, and with $\al_s>0$. We can then inject the moment estimates obtained in Proposition \ref{prop:mom-q-flat} to assert that for any $b\in (\frac12,b^\ast_s]$, 
\begin{align}\label{fina-esti-mo}
&\mathbb{E}\Big[\big\| \cl^{(n,n+1),+}\big\|_{\mathfrak{L}_\mu(Z^{s,b},Z^{s,b})}^2\Big]\lesssim 2^{-2n\ga}, 
\end{align}
for some $\mu\in (0,1)$, and where $\ga:=\min(\al_s,\ka_s)>0$.

\smallskip

The latter estimate immediately entails the convergence of the sequence $\big(\cl^{(n),+}\big)_{n\geq 1}$ in the space $L^2\big(\Omega;\mathfrak{L}_\mu(Z^{s,b},Z^{s,b})\big)$. By combining \eqref{fina-esti-mo} with an elementary Borel-Cantelli argument, this $L^2(\Omega)$-convergence can then be turned into an almost-sure convergence in $\mathfrak{L}_\mu(Z^{s,b},Z^{s,b})$, as desired.

\

\begin{acks}
The authors received support from a Tohoku University-Universit\'e de Lorraine joint research fund. A. Deya and L. Thomann were partially supported by the ANR projet "SMOOTH" ANR-22-CE40-0017. R. Fukuizumi was supported by JSPS KAKENHI Grant Numbers JP19KK0066, JP20K03669, and JP23K25776. A. Deya would like to thank the \textit{Kyoto Research Institute for Mathematical Sciences} for their support during the workshop "Nonlinear and Random Waves", where part of this work has been carried out. The authors would like to thank Herbert Koch and Tristan Robert for enlightening discussions as well as the anonymous referee of a previous version of the paper whose comments led to fundamental improvements of our work.\\
\indent Finally, the authors are grateful to the two referees of the present version of the article for their reading work and valuable comments.
\end{acks}

\

\appendix

\section{Some proofs of deterministic estimates}\label{AppenA}

\subsection{Proof of Lemma  \ref{lem:int-singu-zero}} \label{subsec:append-1}

We divide the proof into three steps.

\smallskip

\noindent
\textit{Step 1}. Let us establish that
\begin{equation}\label{step-1-pro}
\int_{\R}\frac{d\xi}{|\xi|^\nu} \frac{1}{\langle\xi-a\rangle}\frac{1}{\langle \xi-b\rangle} \lesssim \frac{1}{\langle b-a\rangle^{1-\varepsilon}}.
\end{equation}
In fact, this estimate is an almost straightforward consequence of Lemma \ref{lem:tec-1d}. Write first
\begin{equation}\label{simple-proof}
\int_{\R}\frac{d\xi}{|\xi|^\nu} \frac{1}{\langle\xi-a\rangle}\frac{1}{\langle \xi-b\rangle} \lesssim \int_{|\xi|\leq 1}\frac{d\xi}{|\xi|^\nu} \frac{1}{\langle\xi-a\rangle}\frac{1}{\langle \xi-b\rangle}+\int_{|\xi|\geq 1}\frac{d\xi}{|\xi|^\nu} \frac{1}{\langle\xi-a\rangle}\frac{1}{\langle \xi-b\rangle}.
\end{equation}
The second integral can be trivially bounded as
$$\int_{|\xi|\geq 1}\frac{d\xi}{|\xi|^\nu} \frac{1}{\langle\xi-a\rangle}\frac{1}{\langle \xi-b\rangle}\lesssim \int_{\R}\frac{d\xi}{\langle\xi-a\rangle\langle \xi-b\rangle}\lesssim \frac{1}{\langle b-a\rangle^{1-\varepsilon}},$$
where we have used Lemma \ref{lem:tec-1d} to obtain the last inequality.

\smallskip

To control the first integral in \eqref{simple-proof}, introduce $p>1$ close enough to $1$ so that $0\leq p\nu <1$. Then for $q>1$ such that $\frac{1}{p}+\frac{1}{q}=1$, one has
$$\int_{|\xi|\leq 1}\frac{d\xi}{|\xi|^\nu} \frac{1}{\langle\xi-a\rangle}\frac{1}{\langle \xi-b\rangle}\lesssim \bigg(\int_{|\xi|\leq 1}\frac{d\xi}{|\xi|^{p\nu}} \bigg)^{\frac{1}{p}} \bigg(\int_{\R}\frac{d\xi}{\langle\xi-a\rangle^q \langle \xi-b\rangle^q}\bigg)^{\frac{1}{q}},$$
and we can again rely on Lemma \ref{lem:tec-1d} to assert that
$$\bigg(\int_{\R}\frac{d\xi}{\langle\xi-a\rangle^q \langle \xi-b\rangle^q}\bigg)^{\frac{1}{q}}\lesssim \frac{1}{\langle b-a\rangle},$$
which concludes the proof of \eqref{step-1-pro}.

\

\noindent
\textit{Step 2}. Let us now show that
\begin{equation*} 
\int_{\R}\frac{d\xi}{|\xi|^\nu} \frac{1}{\langle\xi-a\rangle}\frac{1}{\langle \xi-b\rangle} \lesssim \frac{1}{|a|^\nu}\frac{1}{\langle b-a\rangle^{1-\varepsilon}}.
\end{equation*}
To this end, write 
\begin{align*}
\int_{\R}\frac{d\xi}{|\xi|^\nu} \frac{1}{\langle\xi-a\rangle}\frac{1}{\langle \xi-b\rangle} =\int_{|\xi|\geq \frac{|a|}{2}}\frac{d\xi}{|\xi|^\nu} \frac{1}{\langle\xi-a\rangle}\frac{1}{\langle \xi-b\rangle}+\int_{|\xi|\leq \frac{|a|}{2}}\frac{d\xi}{|\xi|^\nu} \frac{1}{\langle\xi-a\rangle}\frac{1}{\langle \xi-b\rangle}.
\end{align*}
Just as in Step 1, the integral over $\{|\xi|\geq \frac{|a|}{2}\}$ can be trivially bounded as
$$\int_{|\xi|\geq \frac{|a|}{2}}\frac{d\xi}{|\xi|^\nu} \frac{1}{\langle\xi-a\rangle}\frac{1}{\langle \xi-b\rangle}\lesssim \frac{1}{|a|^\nu}\int_{\R}\frac{d\xi}{\langle\xi-a\rangle\langle \xi-b\rangle}\lesssim \frac{1}{|a|^\nu}\frac{1}{\langle b-a\rangle^{1-\varepsilon}}.$$

\smallskip

Let us now focus on the integral over $\{|\xi|\leq \frac{|a|}{2}\}$, and note first that, for symmetry reasons, we can assume that $b\geq 0$.

\smallskip

\noindent
For $|\xi|\leq \frac{|a|}{2}$ we have $\langle\xi-a\rangle\geq \frac{|a|}{2}$. Next, we observe that for all $|a| \leq b$ we have $\langle b-a\rangle\leq 4  \langle b-\frac{|a|}2\rangle  \lesssim  \langle b-\xi\rangle$.
Then 
\begin{align*}
\int_{|\xi|\leq \frac{|a|}{2}}\frac{d\xi}{|\xi|^\nu} \frac{1}{\langle\xi-a\rangle}\frac{1}{\langle \xi-b\rangle} &\lesssim  \frac1{|a| } \frac{1}{\langle b -a\rangle}  \int_{|\xi|\leq \frac{|a|}{2}}\frac{d\xi}{|\xi|^\nu} \nonumber\\
&\lesssim  \frac1{|a|^\nu } \frac{1}{\langle b -a\rangle}.
\end{align*}

\noindent
\textit{Step 3}. Let us now turn to the proof of \eqref{adapt-lem-3-bis}, which again reduces to elementary integral estimates. Namely, write
\begin{align}
\int_{\R}\frac{d\xi}{|\xi|^\nu} \frac{1}{\langle\xi-a\rangle^p} &\lesssim \int_0^\infty \frac{d\xi}{|\xi|^\nu} \frac{1}{\langle|\xi|-|a|\rangle^p}\nonumber\\
&\lesssim \int_0^{\frac{|a|}{2}} \frac{d\xi}{|\xi|^\nu} \frac{1}{\langle|\xi|-|a|\rangle^p}+\int_{\frac{3|a|}{2}}^\infty \frac{d\xi}{|\xi|^\nu} \frac{1}{\langle|\xi|-|a|\rangle^p}+\int_{\frac{|a|}{2}} ^{\frac{3|a|}{2}}\frac{d\xi}{|\xi|^\nu} \frac{1}{\langle|\xi|-|a|\rangle^p}.\label{step-three}
\end{align}
The first two integrals can be bounded as
\begin{align*}
& \int_0^{\frac{|a|}{2}} \frac{d\xi}{|\xi|^\nu} \frac{1}{\langle|\xi|-|a|\rangle^p}+\int_{\frac{3|a|}{2}}^\infty \frac{d\xi}{|\xi|^\nu} \frac{1}{\langle|\xi|-|a|\rangle^p}\lesssim \frac{1}{|a|^p} \int_0^{\frac{|a|}{2}} \frac{d\xi}{|\xi|^\nu}+\int_{\frac{3|a|}{2}}^\infty \frac{d\xi}{|\xi|^{\nu+p}} \lesssim \frac{1}{|a|^{\nu+p-1}}.
\end{align*}

As for the third integral in \eqref{step-three}, one has
\begin{align*}
&\int_{\frac{|a|}{2}} ^{\frac{3|a|}{2}}\frac{d\xi}{|\xi|^\nu} \frac{1}{\langle|\xi|-|a|\rangle^p}\lesssim    \frac{1}{|a|^{\nu}}     \int_{\frac{|a|}{2}} ^{\frac{3|a|}{2}}\frac{d\xi}{\langle|\xi|-|a|\rangle^p}\lesssim \frac{1}{|a|^{\nu}}\begin{cases} 
 \log(2+|a|) & \mbox{if $p=1$} \\
1 & \mbox{if $p>1$},
\end{cases}
\end{align*}
which completes the proof of \eqref{adapt-lem-3-bis}. \medskip

 For the proof of \eqref{adapt-lem-33}, we write $|\xi|^\mu \leq \langle \xi-a \rangle^\mu+ \langle a \rangle^\mu$, then  using similar considerations as previously,
  \begin{eqnarray*} 
\int_{\R}\frac{|\xi|^\mu d\xi}{\langle\xi-a\rangle \langle \xi-b\rangle} &\lesssim& \langle a\rangle^\mu  \int_{\R}\frac{  d\xi}{\langle\xi-a\rangle \langle \xi-b\rangle}   +   \int_{\R}\frac{  d\xi}{\langle\xi-a\rangle^{1-\mu} \langle \xi-b\rangle} \\
&\lesssim&  \frac{\langle a\rangle^\mu}{\langle b-a\rangle^{1-\varepsilon}}+ \frac{1}{\langle b-a\rangle^{1-\mu-\varepsilon}}\\
&\lesssim&\frac{\langle a\rangle^\mu}{\langle b-a\rangle^{1-\mu-\varepsilon}},
\end{eqnarray*}
  which was the claim.
    \medskip

 For the proof of \eqref{adapt-lem-4-bis}, we use again that $|\xi|^\mu \leq \langle \xi-a \rangle^\mu+ \langle a \rangle^\mu$ and   the fact that $p- \mu >1$ which implies the convergence of the integral.

\subsection{A technical result}

\begin{lemma}\label{lem:ajout}
Let $0\leq b\leq 1$. Then for all $\la_1,\la_2, \la \in \R$, it holds that
\begin{equation} \label{borne-00}
\frac{1}{ \langle \la_1\rangle  \langle \la_2\rangle  \langle     \lambda-(\la_1-\la_2)\rangle^{1-b}} \lesssim  \frac{1}{ \langle\lambda \rangle^{1-b}  \langle \la_1\rangle^b  \langle \la_2\rangle^b}.
 \end{equation}
\end{lemma}

\begin{proof}
To begin with, let us observe that for all $x,y \in \R$ we have 
\begin{equation}\label{inegxy}
\langle x\rangle \leq \sqrt{2} \langle y\rangle \langle x-y\rangle.
\end{equation}
Namely, thanks to an homogeneity argument we have $x^2 \leq 2 \big(y^2+(x-y)^2\big)$, thus 
$$1+x^2 \leq 2\big(1+y^2+(x-y)^2+ y^2(x-y)^2 \big),$$ which is the square of \eqref{inegxy}. Now, we apply \eqref{inegxy} twice to get
\begin{eqnarray*}
\langle\lambda \rangle &\leq & \sqrt{2} \langle     \lambda-(\la_1-\la_2)\rangle  \langle     \la_1-\la_2\rangle \\
 &\leq & 2 \langle     \lambda-(\la_1-\la_2)\rangle  \langle     \la_1\rangle \langle     \la_2\rangle,
\end{eqnarray*}
which in turn implies \eqref{borne-00}, since $0 \leq b \leq 1$.
\end{proof}

\subsection{Proof of Proposition \ref{prop:ice-cream}, item $(ii)$}\label{prop-ice-cream-ii}

For $t>0$ and $k\neq 0$, one has
\begin{align*}
&\mathbb{E}\Big[ \big| \<Psi>^{(n)}_k(t)\big|^2\Big]=\mathbb{E}\bigg[ \bigg|\int_0^t ds \, e^{-\imath s k^2} \dot{B}_k(s)\bigg|^2\bigg]\\
&=\int_{[0,t]^2} ds ds' \, e^{-\imath (s-s')k^2} \mathbb{E}\Big[ \dot{B}^{(n)}_k(s) \overline{\dot{B}^{(n)}_k(s')}\Big]\\
&=\int_{[0,t]^2} ds ds'\, e^{-\imath (s-s')k^2} \int_{[0,2\pi]^2} dx dx' \, e^{-\imath (x-x')k}\mathbb{E}\Big[ \dot{B}^{(n)}(s,x) \overline{\dot{B}^{(n)}(s',x')}\Big]\\
&=\int_{[0,t]^2} ds ds' \, e^{-\imath (s-s')k^2}\int_{[0,2\pi]^2} dx dx' \, e^{-\imath (x-x')k} \int_{\{|\xi|\leq 2^{2n}\}} \frac{d\xi}{|\xi|^{2H_0-1}}\int_{\{|\eta|\leq 2^n\}} \frac{d\eta}{|\eta|^{2H_1-1}} e^{-\imath \xi(s-s')} e^{-\imath \eta(x-x')}\\
&=\bigg(\int_{\{|\xi|\leq 2^{2n}\}} \frac{d\xi}{|\xi|^{2H_0-1}} \bigg|\int_0^t ds \, e^{-\imath s (\xi+k^2)}\bigg|^2\bigg)\bigg(\int_{\{|\eta|\leq 2^n\}} \frac{d\eta}{|\eta|^{2H_1-1}}\bigg|\int_0^{2\pi} dx  \, e^{-\imath x(k+\eta)}\bigg|^2\bigg)\\
&=16\bigg(\int_{\{|\xi|\leq 2^{2n}\}} \frac{d\xi}{|\xi|^{2H_0-1}}\frac{ \big| \sin\big(\frac{t}{2}(\xi-k^2)\big)\big|^2}{|\xi-k^2|^2}\bigg)\bigg(\int_{\{|\eta|\leq 2^n\}} \frac{d\eta}{|\eta|^{2H_1-1}}\frac{ \big| \sin\big(\pi(\eta-k)\big)\big|^2}{|\eta-k|^2}\bigg)=:16 I^{(n)}_k J^{(n)}_k.
\end{align*}
Let us decompose $I^{(n)}_k$ and $J^{(n)}_k$ as $I^{(n)}_k:=I^{(n),-}_k+I^{(n),+}_k$ and $J^{(n)}_k:=J^{(n),-}_k+J^{(n),+}_k$, with
\small
$$I^{(n),-}_k:= \int_{\{|\xi|\leq 2^{2n}\}\cap [\frac{k^2}{4},4k^2]} \frac{d\xi}{|\xi|^{2H_0-1}}\frac{ \big| \sin\big(\frac{t}{2}(\xi-k^2)\big)\big|^2}{|\xi-k^2|^2},$$
$$ I^{(n),+}_k:=\int_{\{|\xi|\leq 2^{2n}\}\backslash [\frac{k^2}{4},4k^2]} \frac{d\xi}{|\xi|^{2H_0-1}}\frac{ \big| \sin\big(\frac{t}{2}(\xi-k^2)\big)\big|^2}{|\xi-k^2|^2},$$
\normalsize
resp.
\small
$$J^{(n),-}_k:=\int_{\{|\eta|\leq 2^n\}\cap [\frac{k}{2},2k]} \frac{d\eta}{|\eta|^{2H_1-1}}\frac{ \big| \sin\big(\pi(\eta-k)\big)\big|^2}{|\eta-k|^2},$$
$$J^{(n),+}_k:=\int_{\{|\eta|\leq 2^n\}\backslash [\frac{k}{2},2k]} \frac{d\eta}{|\eta|^{2H_1-1}}\frac{ \big| \sin\big(\pi(\eta-k)\big)\big|^2}{|\eta-k|^2}.$$
\normalsize
As a result,
\begin{align}
\si^{(n)}(t)=\sum_{k\neq 0} \mathbb{E}\Big[ \big| \<Psi>^{(n)}(t)\big|^2\Big]&=\cp^{(n)}(t)+\mathcal{R}^{(n)}(t),\label{decompo-si-0}
\end{align}
with $\cp^{(n)}(t):=32\sum_{1\leq k<2^{n}}I^{(n),-}_k J^{(n),-}_k$ and
\begin{align*}
\mathcal{R}^{(n)}(t):=32\sum_{k\geq  2^n}I^{(n),-}_k J^{(n),-}_k+32\sum_{k=1}^{\infty} \big[I^{(n),-}_k J^{(n),+}_k+I^{(n),+}_k J^{(n),-}_k+I^{(n),+}_k J^{(n),+}_k\big].
\end{align*}

\

\

\noindent
\textit{Step 1: Estimate of $I^{(n),-}_k, J^{(n),-}_k$ for $1\leq k< 2^{n-1}$}. Note first that for $1\leq k\leq 2^{n-1}$, one has
$$I^{(n),-}_k= \int_{\frac{k^2}{4}}^{4k^2 \wedge 2^{2n}} \frac{d\xi}{|\xi|^{2H_0-1}}\frac{ \big| \sin\big(\frac{t}{2}(\xi-k^2)\big)\big|^2}{|\xi-k^2|^2}=\int_{\frac{k^2}{4}}^{4k^2 } \frac{d\xi}{|\xi|^{2H_0-1}}\frac{ \big| \sin\big(\frac{t}{2}(\xi-k^2)\big)\big|^2}{|\xi-k^2|^2}.$$
Now for every $1\leq k < 2^{n-1}$, let us write
\small
\begin{align*}
&I^{(n),-}_k=\frac{1}{k^{4H_0}}\int_{\frac14}^{4 } \frac{d\xi}{|\xi|^{2H_0-1}}\frac{ \big| \sin\big(\frac{tk^2}{2}(1-\xi)\big)\big|^2}{|1-\xi|^2}=\frac{1}{k^{4H_0}}\int_{-3 }^{\frac34} \frac{d\xi}{|1-\xi|^{2H_0-1}}\frac{ \big| \sin\big(\frac{tk^2}{2}\xi\big)\big|^2}{|\xi|^2}\\
&=\frac{t}{2k^{4H_0-2}}\int_{-\frac{3t k^2}{2}}^{\frac{3t k^2}{8}} \frac{d\xi}{\big|1-\frac{2\xi}{tk^2}\big|^{2H_0-1}}\frac{ | \sin(\xi)|^2}{|\xi|^2}\\
&=\frac{t}{2k^{4H_0-2}}\int_{-\infty}^{\infty}\frac{ | \sin(\xi)|^2}{|\xi|^2}\, d\xi+\frac{t}{2k^{4H_0-2}}\bigg[\int_{-\frac{3t k^2}{2}}^{\frac{3t k^2}{8}} \frac{d\xi}{\big|1-\frac{2\xi}{tk^2}\big|^{2H_0-1}}\frac{ | \sin(\xi)|^2}{|\xi|^2}-\int_{-\infty}^{\infty} d\xi\frac{ | \sin(\xi)|^2}{|\xi|^2}\, d\xi\bigg].
\end{align*}
\normalsize
Using the subsequent Lemma \ref{lem:tech} and setting
\begin{equation}\label{defi-a}
A:=\int_{-\infty}^{\infty}\frac{ | \sin(\xi)|^2}{|\xi|^2}\, d\xi,
\end{equation}
we deduce that for all $1\leq k <  2^{n-1}$ and $\varepsilon \in (0,1)$,
\small
\begin{align*}
\Big|I^{(n),-}_k-\frac{At}{2k^{4H_0-2}}\Big|\lesssim  \frac{t^\varepsilon}{k^{4H_0-2\varepsilon}} .
\end{align*}
\normalsize
With the same arguments, we get that for all $1\leq k < 2^{n-1}$ and $\varepsilon \in (0,1)$,
\begin{align*}
\big|J^{(n),-}_k-\frac{A\pi}{k^{2H_1-1}}\Big|\lesssim \frac{1}{k^{2H_1-\varepsilon}}.
\end{align*}
As a consequence of these two estimates, we can write for $\varepsilon >0$ small enough,
\small
\begin{align}
&\bigg| \sum_{1\leq k<  2^{n-1}} I^{(n),-}_k J^{(n),-}_k- \frac{A^2\pi t}{2}\sum_{1\leq k<  2^{n-1}}\frac{1}{k^{2\underline{H}-1}}\bigg|\leq \sum_{1\leq k<  2^{n-1}}\Big| I^{(n),-}_k J^{(n),-}_k- \frac{A^2\pi t}{2k^{2\underline{H}-1}}\Big|\nonumber\\
&\lesssim \sum_{1\leq k<  2^{n-1}}\Big| I^{(n),-}_k - \frac{A t}{2k^{4H_0-2}}\Big|\Big| J^{(n),-}_k- \frac{A\pi }{k^{2H_1-1}}\Big|+\sum_{1\leq k<  2^{n-1}}\Big| I^{(n),-}_k- \frac{A t}{2k^{4H_0-2}}\Big|\frac{1}{k^{2H_1-1}}\nonumber\\
&\hspace{2cm}+\sum_{1\leq k<  2^{n-1}}\frac{1}{k^{4H_0-2}}\Big|\Big| J^{(n),-}_k- \frac{A\pi }{k^{2H_1-1}}\Big|\nonumber\\
&\lesssim \sum_{k=1}^\infty\frac{1}{k^{2(2H_0+H_1)-3\varepsilon}}+\sum_{k=1}^\infty\frac{1}{k^{2(2H_0+H_1)-1-2\varepsilon}}+\sum_{k=1}^\infty\frac{1}{k^{2(2H_0+H_1)-2-\varepsilon}}<\infty,\label{step-1}
\end{align}
\normalsize
due to the condition $2H_0+H_1>\frac32$.

\

\noindent
\textit{Step 2: Estimate of $I^{(n),-}_k, J^{(n),-}_k$ for $2^{n-1}\leq k< 2^n$}. For $2^{n-1}\leq k< 2^n$, one has in fact
$$I^{(n),-}_k= \int_{\frac{k^2}{4}}^{4k^2 \wedge 2^{2n}} \frac{d\xi}{|\xi|^{2H_0-1}}\frac{ \big| \sin\big(\frac{t}{2}(\xi-k^2)\big)\big|^2}{|\xi-k^2|^2}=\int_{\frac{k^2}{4}}^{2^{2n} } \frac{d\xi}{|\xi|^{2H_0-1}}\frac{ \big| \sin\big(\frac{t}{2}(\xi-k^2)\big)\big|^2}{|\xi-k^2|^2}.$$
With the same changes of variables as above, we can then write
\small
\begin{align*}
&I^{(n),-}_k=\frac{1}{k^{4H_0}}\int_{\frac14}^{\frac{2^{2n}}{k^2} } \frac{d\xi}{|\xi|^{2H_0-1}}\frac{ \big| \sin\big(\frac{tk^2}{2}(1-\xi)\big)\big|^2}{|1-\xi|^2}=\frac{1}{k^{4H_0}}\int_{1-\frac{2^{2n}}{k^2}  }^{\frac34} \frac{d\xi}{|1-\xi|^{2H_0-1}}\frac{ \big| \sin\big(\frac{tk^2}{2}\xi\big)\big|^2}{|\xi|^2}\\
&=\frac{t}{2k^{4H_0-2}}\int_{\frac{t k^2}{2}(1-\frac{2^{2n}}{k^2})}^{\frac{3t k^2}{8}} \frac{d\xi}{\big|1-\frac{2\xi}{tk^2}\big|^{2H_0-1}}\frac{ | \sin(\xi)|^2}{|\xi|^2}\\
&=\frac{t}{2k^{4H_0-2}}\bigg[\int_{\frac{t k^2}{2}(1-\frac{2^{2n}}{k^2})}^{\frac{3t k^2}{8}} \frac{ | \sin(\xi)|^2}{|\xi|^2}+\int_{\frac{t k^2}{2}(1-\frac{2^{2n}}{k^2})}^{\frac{3t k^2}{8}} \bigg(\frac{1}{\big|1-\frac{2\xi}{tk^2}\big|^{2H_0-1}}-1\bigg)\frac{ | \sin(\xi)|^2}{|\xi|^2}\bigg].
\end{align*}
\normalsize
For every $2^{n-1}\leq k < 2^{n}$, one has $ 1-\frac{2^{2n}}{k^2}\in [-3,0)$, and so
\small
\begin{align*}
&\bigg|\int_{\frac{t k^2}{2}(1-\frac{2^{2n}}{k^2})}^{\frac{3t k^2}{8}} \bigg(\frac{1}{\big|1-\frac{2\xi}{tk^2}\big|^{2H_0-1}}-1\bigg)\frac{ | \sin(\xi)|^2}{|\xi|^2}\bigg|\leq \int_{-\frac{3t k^2}{2}}^{\frac{3t k^2}{8}} \bigg|\frac{1}{\big|1-\frac{2\xi}{tk^2}\big|^{2H_0-1}}-1\bigg|\frac{ | \sin(\xi)|^2}{|\xi|^2}\lesssim \frac{1}{(t k^2)^{1-\varepsilon}},
\end{align*}
\normalsize
where we have used Lemma \ref{lem:tech} to derive the second inequality.

\smallskip

Besides, for $2^{n-1}\leq k <2^{n}$,
\begin{align*}
&\int_{\frac{t k^2}{2}(1-\frac{2^{2n}}{k^2})}^{\frac{3t k^2}{8}} \frac{ | \sin(\xi)|^2}{|\xi|^2}=\int_{-\infty}^{\infty} \frac{ | \sin(\xi)|^2}{|\xi|^2}-\bigg[\int_{-\infty}^{\infty} \frac{ | \sin(\xi)|^2}{|\xi|^2}-\int_{\frac{t k^2}{2}(1-\frac{2^{2n}}{k^2})}^{\frac{3t k^2}{8}} \frac{ | \sin(\xi)|^2}{|\xi|^2}\bigg]\\
&=\int_{-\infty}^{\infty} \frac{ | \sin(\xi)|^2}{|\xi|^2}-\bigg[\int_{-\infty}^{\frac{t k^2}{2}(1-\frac{2^{2n}}{k^2})} \frac{ | \sin(\xi)|^2}{|\xi|^2}+\int_{\frac{3t k^2}{8}}^\infty \frac{ | \sin(\xi)|^2}{|\xi|^2}\bigg]\\
&=\int_{-\infty}^{\infty} \frac{ | \sin(\xi)|^2}{|\xi|^2}-\bigg[\int_{\frac{t k^2}{2}|1-\frac{2^{2n}}{k^2}|}^{\infty} \frac{ | \sin(\xi)|^2}{|\xi|^2}+\int_{\frac{3t k^2}{8}}^\infty \frac{ | \sin(\xi)|^2}{|\xi|^2}\bigg],
\end{align*}
which gives 
\begin{align*}
&\bigg|\int_{\frac{t k^2}{2}(1-\frac{2^{2n}}{k^2})}^{\frac{3t k^2}{8}} \frac{ | \sin(\xi)|^2}{|\xi|^2}-\int_{-\infty}^{\infty} \frac{ | \sin(\xi)|^2}{|\xi|^2}\bigg|\lesssim \int_{\frac{t k^2}{2}|1-\frac{2^{2n}}{k^2}|}^{\infty} \frac{ | \sin(\xi)|^2}{|\xi|^2}\\
&\lesssim \min\bigg(\frac{1}{\big|t k^2(1-\frac{2^{2n}}{k^2})\big|^{\frac12-\varepsilon}}\int_0^\infty d\xi \, \frac{ | \sin(\xi)|^2}{|\xi|^{\frac32+\varepsilon}}, \frac{1}{\big|t k^2(1-\frac{2^{2n}}{k^2})\big|^{1-\varepsilon}}\int_0^\infty d\xi \, \frac{ | \sin(\xi)|^2}{|\xi|^{1+\varepsilon}}\bigg).
\end{align*}
By combining the above estimates, we obtain that for every $2^{n-1}\leq k<  2^n$,
\begin{align*}
&\bigg|I^{(n),-}_k -\frac{At}{2k^{4H_0-2}} \bigg|\lesssim \min\bigg( \frac{t^{\frac12+\varepsilon}}{k^{4H_0-1-2\varepsilon}}\frac{1}{|1-\frac{2^{2n}}{k^2}|^{\frac12-\varepsilon}},\frac{t^{\varepsilon}}{k^{4H_0-2\varepsilon}}\frac{1}{|1-\frac{2^{2n}}{k^2}|^{1-\varepsilon}}\bigg).
\end{align*}
With the same arguments, we get that for every $2^{n-1}\leq k<  2^n$,
\begin{align*}
&\bigg|J^{(n),-}_k -\frac{A\pi}{k^{2H_1-1}}\bigg|\lesssim  \min\bigg(\frac{1}{k^{2H_1-\frac12-\varepsilon}}\frac{1}{|1-\frac{2^{n}}{k}|^{\frac12-\varepsilon}},\frac{1}{k^{2H_1-\varepsilon}}\frac{1}{|1-\frac{2^{n}}{k}|^{1-\varepsilon}}\bigg).
\end{align*}
As a consequence of these two estimates, we can write for $\varepsilon >0$ small enough,
\small
\begin{align*}
&\bigg| \sum_{2^{n-1}\leq k<  2^n} I^{(n),-}_k J^{(n),-}_k- \frac{A^2\pi t}{2}\sum_{2^{n-1}\leq k<  2^n}\frac{1}{k^{2\underline{H}-1}}\bigg|\leq \sum_{2^{n-1}\leq k<  2^n}\Big| I^{(n),-}_k J^{(n),-}_k- \frac{A^2\pi t}{2k^{2\underline{H}-1}}\Big|\\
&\lesssim \sum_{2^{n-1}\leq k<  2^n}\Big| I^{(n),-}_k - \frac{A t}{2k^{4H_0-2}}\Big|\Big| J^{(n),-}_k- \frac{A\pi }{k^{2H_1-1}}\Big|+\sum_{2^{n-1}\leq k<  2^n}\Big| I^{(n),-}_k- \frac{A t}{2k^{4H_0-2}}\Big|\frac{1}{k^{2H_1-1}}\\
&\hspace{2cm}+\sum_{2^{n-1}\leq k<  2^n}\frac{1}{k^{4H_0-2}}\Big|\Big| J^{(n),-}_k- \frac{A\pi }{k^{2H_1-1}}\Big|\\
&\lesssim \sum_{2^{n-1}\leq k<  2^n}\frac{1}{k^{2(2H_0+H_1)-\frac32-3\varepsilon}}\frac{1}{|1-\frac{2^{2n}}{k^2}|^{\frac12-\varepsilon}}\frac{1}{|1-\frac{2^{n}}{k}|^{\frac12-\varepsilon}}+\sum_{2^{n-1}\leq k<  2^n}\frac{t^{\varepsilon}}{k^{2(2H_0+H_1)-1-2\varepsilon}}\frac{1}{|1-\frac{2^{2n}}{k^2}|^{1-\varepsilon}}\\
&\hspace{2cm}+\sum_{2^{n-1}\leq k<  2^n}\frac{1}{k^{2(2H_0+H_1)-2-\varepsilon}}\frac{1}{|1-\frac{2^{n}}{k}|^{1-\varepsilon}}\\
&\lesssim \frac{1}{2^{n(2(2H_0+H_1)-\frac52-3\varepsilon)}} \bigg( \frac{1}{2^n}\sum_{2^{n-1}\leq k<  2^n}\frac{1}{(\frac{k}{2^n})^{2(2H_0+H_1)-\frac32-3\varepsilon}}\frac{1}{|1-\frac{2^{2n}}{k^2}|^{\frac12-\varepsilon}}\frac{1}{|1-\frac{2^{n}}{k}|^{\frac12-\varepsilon}}\bigg)\\
&\hspace{1cm}+\frac{1}{2^{n(2(2H_0+H_1)-2-2\varepsilon)}} \bigg(\frac{1}{2^n}\sum_{2^{n-1}\leq k<  2^n}\frac{t^{\varepsilon}}{(\frac{k}{2^n})^{2(2H_0+H_1)-1-2\varepsilon}}\frac{1}{|1-\frac{2^{2n}}{k^2}|^{1-\varepsilon}}\bigg)\\
&\hspace{2cm}+\frac{1}{2^{n(2(2H_0+H_1)-3-\varepsilon)}} \bigg(\frac{1}{2^n}\sum_{2^{n-1}\leq k<  2^n}\frac{1}{(\frac{k}{2^n})^{2(2H_0+H_1)-2-\varepsilon}}\frac{1}{|1-\frac{2^{n}}{k}|^{1-\varepsilon}}\bigg),
\end{align*}
\normalsize
and thus, for $\varepsilon >0$ small enough,
\small
\begin{align}
&\bigg| \sum_{2^{n-1}\leq k<  2^n} I^{(n),-}_k J^{(n),-}_k- \frac{A^2\pi t}{2}\sum_{2^{n-1}\leq k<  2^n}\frac{1}{k^{2\underline{H}-1}}\bigg|\nonumber\\
&\lesssim \frac{1}{2^{n(2(2H_0+H_1)-3-\varepsilon)}} \bigg[ \int_{\frac12}^1 \frac{dx}{\big|1-\frac{1}{x^2}\big|^{\frac12-\varepsilon}\big|1-\frac{1}{x}\big|^{\frac12-\varepsilon}} +\int_{\frac12}^1 \frac{dx}{\big|1-\frac{1}{x^2}\big|^{1-\varepsilon}}+\int_{\frac12}^1 \frac{dx}{\big|1-\frac{1}{x}\big|^{1-\varepsilon}} \bigg]\nonumber\\
&\lesssim \frac{1}{2^{n(2(2H_0+H_1)-3-\varepsilon)}} \int_{\frac12}^1 \frac{dx}{\big|1-x\big|^{1-\varepsilon}}\lesssim \frac{1}{2^{n(2(2H_0+H_1)-3-\varepsilon)}} \stackrel{n\to\infty}{\longrightarrow} 0,\label{step-2}
\end{align}
\normalsize
due to $2H_0+H_1>\frac32$.

\

\noindent
\textit{Step 3: Estimate of $I^{(n),-}_k, J^{(n),-}_k$ for $k=2^n$}. In this particular case, one has 
$$I^{(n),-}_k=\int_{\frac{k^2}{4}}^{2^{2n} } \frac{d\xi}{|\xi|^{2H_0-1}}\frac{ \big| \sin\big(\frac{t}{2}(\xi-k^2)\big)\big|^2}{|\xi-k^2|^2}=\int_{\frac{k^2}{4}}^{k^2 } \frac{d\xi}{|\xi|^{2H_0-1}}\frac{ \big| \sin\big(\frac{t}{2}(\xi-k^2)\big)\big|^2}{|\xi-k^2|^2}.$$
With the same changes of variables as above, we can then write
\small
\begin{align*}
I^{(n),-}_k&=\frac{t}{2k^{4H_0-2}}\int_{0}^{\frac{3t k^2}{8}} \frac{d\xi}{\big|1-\frac{2\xi}{tk^2}\big|^{2H_0-1}}\frac{ | \sin(\xi)|^2}{|\xi|^2}\\
&=\frac{t}{2k^{4H_0-2}}\bigg[\int_{0}^{\frac{3t k^2}{8}} \frac{ | \sin(\xi)|^2}{|\xi|^2}+\int_{0}^{\frac{3t k^2}{8}} \bigg(\frac{1}{\big|1-\frac{2\xi}{tk^2}\big|^{2H_0-1}}-1\bigg)\frac{ | \sin(\xi)|^2}{|\xi|^2}\bigg].
\end{align*}
\normalsize
By Lemma \ref{lem:tech}, we have
\small
\begin{align*}
&\bigg|\int_{0}^{\frac{3t k^2}{8}} \bigg(\frac{1}{\big|1-\frac{2\xi}{tk^2}\big|^{2H_0-1}}-1\bigg)\frac{ | \sin(\xi)|^2}{|\xi|^2}\bigg|\lesssim \frac{1}{(t k^2)^{1-\varepsilon}},
\end{align*}
\normalsize
and so
$$\big|I^{(n),-}_{2^n}\big| \lesssim  \frac{1}{2^{n(4H_0-2)}} \bigg[t+\frac{t^{\varepsilon}}{2^{n(2-2\varepsilon)}}\bigg].$$
In the same way,
$$\big|J^{(n),-}_{2^n}\big| \lesssim  \frac{1}{2^{n(2H_1-1)}} \bigg[1+\frac{1}{2^{n(1-\varepsilon)}}\bigg],$$
and as a result
\begin{equation}\label{step-3}
\big|I^{(n),-}_{2^n}\big|\big|J^{(n),-}_{2^n}\big| \lesssim  \frac{1}{2^{n(4H_0-2)}} \frac{1}{2^{n(2H_1-1)}}\lesssim \frac{1}{2^{n(2\underline{H}-1)}}\stackrel{n\to\infty}{\longrightarrow} 0 .
\end{equation}

\

\noindent
\textit{Step 4: Estimate of $I^{(n),-}_k, J^{(n),-}_k$ for $k> 2^{n}$}. For $k>2^n$, one has in fact
$$I^{(n),-}_k= \int_{\frac{k^2}{4}}^{4k^2 \wedge 2^{2n}} \frac{d\xi}{|\xi|^{2H_0-1}}\frac{ \big| \sin\big(\frac{t}{2}(\xi-k^2)\big)\big|^2}{|\xi-k^2|^2}=\1_{\{k\leq 2^{n+1}\}}\int_{\frac{k^2}{4}}^{2^{2n} } \frac{d\xi}{|\xi|^{2H_0-1}}\frac{ \big| \sin\big(\frac{t}{2}(\xi-k^2)\big)\big|^2}{|\xi-k^2|^2}.$$
Now for every $2^n< k \leq 2^{n+1}$, and just as in Step 2, we can write
\small
\begin{align*}
&I^{(n),-}_k=\frac{t}{2k^{4H_0-2}}\int_{\frac{t k^2}{2}(1-\frac{2^{2n}}{k^2})}^{\frac{3t k^2}{8}} \frac{d\xi}{\big|1-\frac{2\xi}{tk^2}\big|^{2H_0-1}}\frac{ | \sin(\xi)|^2}{|\xi|^2}\\
&=\frac{t}{2k^{4H_0-2}}\bigg[\int_{\frac{t k^2}{2}(1-\frac{2^{2n}}{k^2})}^{\frac{3t k^2}{8}} \frac{ | \sin(\xi)|^2}{|\xi|^2}+\int_{\frac{t k^2}{2}(1-\frac{2^{2n}}{k^2})}^{\frac{3t k^2}{8}} \bigg(\frac{1}{\big|1-\frac{2\xi}{tk^2}\big|^{2H_0-1}}-1\bigg)\frac{ | \sin(\xi)|^2}{|\xi|^2}\bigg].
\end{align*}
\normalsize
Note that for every $2^{n}< k \leq 2^{n+1}$, one has $ 1-\frac{2^{2n}}{k^2}\in (0,\frac34]$, and so
\small
\begin{align*}
&\bigg|\int_{\frac{t k^2}{2}(1-\frac{2^{2n}}{k^2})}^{\frac{3t k^2}{8}} \bigg(\frac{1}{\big|1-\frac{2\xi}{tk^2}\big|^{2H_0-1}}-1\bigg)\frac{ | \sin(\xi)|^2}{|\xi|^2}\bigg|\leq \int_{0}^{\frac{3t k^2}{8}} \bigg|\frac{1}{\big|1-\frac{2\xi}{tk^2}\big|^{2H_0-1}}-1\bigg|\frac{ | \sin(\xi)|^2}{|\xi|^2}\lesssim \frac{1}{(t k^2)^{1-\varepsilon}},
\end{align*}
\normalsize
where we have used Lemma \ref{lem:tech} to derive the second inequality. Since on the other hand
\begin{align*}
&\int_{\frac{t k^2}{2}(1-\frac{2^{2n}}{k^2})}^{\frac{3t k^2}{8}} \frac{ | \sin(\xi)|^2}{|\xi|^2} \lesssim \frac{1}{\big|t k^2(1-\frac{2^{2n}}{k^2})\big|^{\frac12-\varepsilon}}\int_0^\infty d\xi \, \frac{ | \sin(\xi)|^2}{|\xi|^{\frac32+\varepsilon}},
\end{align*} 
we deduce that for every $ k >2^n$,
\begin{align*}
\big|I^{(n),-}_k \big| \lesssim \1_{\{k\leq 2^{n+1}\}}\frac{t^{\frac12+\varepsilon}}{k^{4H_0-1-2\varepsilon}}\frac{1}{|1-\frac{2^{2n}}{k^2}|^{\frac12-\varepsilon}}.
\end{align*}
With the same arguments, we get that for every $ k> 2^{n}$,
\begin{align*}
\big|J^{(n),-}_k \big| \lesssim \1_{\{k\leq 2^{n+1}\}}\frac{1}{k^{2H_1-\frac12-\varepsilon}}\frac{1}{|1-\frac{2^{n}}{k}|^{\frac12-\varepsilon}}.
\end{align*}
As a result, for every $\varepsilon >0$ small enough,
\begin{align}
&\bigg|\sum_{k>2^n} I^{(n),-}_k J^{(n),-}_k\bigg|\lesssim t^{\frac12+\varepsilon}\sum_{2^n <k\leq 2^{n+1}}\frac{1}{k^{4H_0+2H_1-\frac32-3\varepsilon}}\frac{1}{|1-\frac{2^{2n}}{k^2}|^{\frac12-\varepsilon}} \frac{1}{|1-\frac{2^{n}}{k}|^{\frac12-\varepsilon}}\nonumber\\
&\lesssim t^{\frac12+\varepsilon}\frac{2^n}{2^{n(4H_0+2H_1-\frac32-3\varepsilon)}}   \bigg( \frac{1}{2^n}\sum_{2^n <k\leq 2^{n+1}}\frac{1}{(\frac{k}{2^n})^{4H_0+2H_1-\frac32-3\varepsilon}}\frac{1}{|1-\frac{2^{2n}}{k^2}|^{\frac12-\varepsilon}} \frac{1}{|1-\frac{2^{n}}{k}|^{\frac12-\varepsilon}}\bigg)\nonumber\\
&\lesssim \frac{t^{\frac12+\varepsilon}}{2^{n(4H_0+2H_1-\frac52-3\varepsilon)}} \int_1^2 \frac{dx}{|x|^{4H_0+2H_1-\frac32-3\varepsilon}} \frac{1}{|1-\frac{1}{x^2}|^{\frac12-\varepsilon}} \frac{1}{|1-\frac{1}{x}|^{\frac12-\varepsilon}}\nonumber\\
&\lesssim t^{\frac12+\varepsilon} \int_1^2 \frac{dx}{|x-1|^{1-2\varepsilon}}\lesssim t^{\frac12+\varepsilon},\label{step-4}
\end{align}
where we have used the fact that $4H_0+2H_1-\frac52>3-\frac52>0$.

\

\noindent
\textit{Step 5: Estimate of $I^{(n),+}_k, J^{(n),+}_k$}.

\smallskip

If $\xi\notin [\frac{k^2}{4},4k^2]$, then of course $|\xi-k^2| \gtrsim k^2$, and so for every $\varepsilon \in (0,1)$,
\begin{align*}
\big|I^{(n),+}_k\big|&\lesssim \frac{1}{(k^2)^{1-\varepsilon}}\int_{\R} \frac{d\xi}{|\xi|^{2H_0-1}}\frac{ \big| \sin\big(\frac{t}{2}(\xi-k^2)\big)\big|^2}{|\xi-k^2|^{1+\varepsilon}}\\
&\lesssim \frac{1}{(k^2)^{1-\varepsilon}}\int_{\R} \frac{d\xi}{|\xi|^{2H_0-1}}\frac{ 1}{\langle \xi-k^2\rangle^{1+\varepsilon}}\lesssim \frac{1}{(k^2)^{1-\varepsilon}}\frac{1}{(k^2)^{2H_0-1}},
\end{align*}
where we have used Lemma \ref{lem:int-singu-zero} to derive the last inequality. In other words, one has
\begin{equation}\label{step-5-i}
\big|I^{(n),+}_k\big|\lesssim \frac{1}{k^{4H_0-\varepsilon}} ,
\end{equation}
and with the same arguments, we obtain that
\begin{equation}\label{step-5-j}
\big|J^{(n),+}_k\big|\lesssim \frac{1}{k^{2H_1-\varepsilon}} .
\end{equation}

\

\noindent
\textit{Step 6: Conclusion}.

\smallskip

By combining \eqref{step-1} and \eqref{step-2}, we immediately deduce that
\begin{align*}
&\sup_{n\geq 1}\bigg|\cp^{(n)}(t)- 16A^2\pi t\sum_{1\leq k<  2^{n}}\frac{1}{k^{2\underline{H}-1}}\bigg| <\infty.
\end{align*}
As far as $\mathcal{R}^{(n)}(t)$ is concerned, we can first combine \eqref{step-3} and \eqref{step-4} to ensure that
$$\sup_{n\geq 1} \bigg|\sum_{k\geq 2^n} I^{(n),-}_k J^{(n),-}_k\bigg| <\infty.$$
On the other hand, by gathering the estimates obtained in Steps 1 to 4, we easily see that for every $k\geq 1$,
$$\big| I^{(n),-}_k(t)\big| \lesssim \frac{1}{k^{4H_0-2}} \quad \text{and} \quad \big| J^{(n),-}_k(t)\big| \lesssim \frac{1}{k^{2H_1-1}},$$
which, together with \eqref{step-5-i}-\eqref{step-5-j}, yields
\begin{align*}
&\bigg|\sum_{k=1}^{\infty} \big[I^{(n),-}_k J^{(n),+}_k+I^{(n),+}_k J^{(n),-}_k+I^{(n),+}_k J^{(n),+}_k\big]\bigg|\\
&\lesssim \sum_{k=1}^{\infty}\frac{1}{k^{4H_0-2}} \frac{1}{k^{2H_1-\varepsilon}}+\sum_{k=1}^{\infty}\frac{1}{k^{4H_0-\varepsilon}} \frac{1}{k^{2H_1-1}}+ \sum_{k=1}^{\infty}\frac{1}{k^{4H_0-\varepsilon}} \frac{1}{k^{2H_1-\varepsilon}}.
\end{align*}
Since $2H_0+H_1>\frac32$, it is readily checked that these three quantities are finite for $\varepsilon >0$ small enough. We can therefore conclude that
$$\sup_{n\geq 1} \big|\mathcal{R}^{(n)}(t)\big| <\infty.$$

\

Going back to the decomposition \eqref{decompo-si-0}, we have thus shown that
$$\sup_{n\geq 1} \bigg|\si^{(n)}(t) - 16A^2\pi t\sum_{1\leq k<  2^{n}}\frac{1}{k^{2\underline{H}-1}}\bigg| <\infty,$$
and we can finally observe that
\begin{align*}
&\sum_{1\leq k<  2^{n}}\frac{1}{k^{2\underline{H}-1}}=2^{2n(1-\underline{H})} \bigg(\frac{1}{2^n} \sum_{1\leq k<  2^{n}}\frac{1}{(\frac{k}{2^n})^{2\underline{H}-1}}\bigg)\stackrel{n\to\infty}{\sim} 2^{2n(1-\underline{H})} \int_0^1 \frac{dx}{|x|^{2\underline{H}-1}}.
\end{align*}

\

\

It only remains us to deal with the following technical lemma, which has been used during the proof.


\begin{lemma}\label{lem:tech}
Fix $a,b>0$, $\al\in (0,2)$ and $\varepsilon \in (0,1)$. Then for all $r>0$, one has 
\begin{align*}
&\int_{-ar}^{br}d\xi\,  \bigg| \frac{1}{|1-\frac{\xi}{r}|^{\al -1}}-1\bigg| \frac{|\sin(\xi)|^2}{\xi^2}\leq \frac{c_{a,b,\varepsilon,\al}}{r^{1-\varepsilon}}.
\end{align*}
As an immediate consequence, it holds that
\begin{equation}\label{conse}
\bigg| \int_{-ar}^{br} \frac{d\xi}{|1-\frac{\xi}{r}|^{\al -1}} \frac{|\sin(\xi)|^2}{\xi^2}-\int_{-\infty}^{\infty } d\xi\,  \frac{|\sin(\xi)|^2}{\xi^2} \bigg| \leq \frac{c_{a,b,\al,\varepsilon}}{r^{1-\varepsilon}} \ .
\end{equation}
\end{lemma}

\begin{proof}
One has
\begin{align*}
&\int_{-ar}^{br}d\xi\,  \bigg| \frac{1}{|1-\frac{\xi}{r}|^{\al -1}}-1\bigg| \frac{|\sin(\xi)|^2}{\xi^2}
=\frac{1}{r}\int_{-a}^{b}d\xi\,  \frac{| 1-|1-\xi|^{\al-1}|}{|1-\xi|^{\al -1}}  \frac{|\sin(r\xi)|^2}{\xi^2}\\
&\lesssim \frac{1}{r^{1-\varepsilon}}\int_{-a}^{b}\frac{d\xi}{|\xi|^{2-\varepsilon}}\,  \frac{| 1-|1-\xi|^{\al-1}|}{|1-\xi|^{\al -1}}\\
&\lesssim \frac{1}{r^{1-\varepsilon}}\bigg[\int_{-(a\wedge \frac12)}^{b\wedge \frac12}\frac{d\xi}{|\xi|^{2-\varepsilon}}\,  \frac{| 1-|1-\xi|^{\al-1}|}{|1-\xi|^{\al -1}} +c_{a,b,\varepsilon}\int_{-a}^{b}d\xi\,  \frac{| 1-|1-\xi|^{\al-1}|}{|1-\xi|^{\al -1}} \bigg] \\
&\lesssim \frac{c_{a,b,\varepsilon,\al}}{r^{1-\varepsilon}}\bigg[\int_{-(a\wedge \frac12)}^{b\wedge \frac12}\frac{d\xi}{|\xi|^{1-\varepsilon}} +\int_{-a}^{b}d\xi\,  \bigg(1+\frac{1}{|1-\xi|^{\al-1}} \bigg)\bigg] ,
\end{align*}
hence the first inequality.

\smallskip

As for \eqref{conse}, we can of course bound the quantity under consideration by
$$\int_{-ar}^{br}d\xi\,  \bigg| \frac{1}{|1-\frac{\xi}{r}|^{\al -1}}-1\bigg| \frac{|\sin(\xi)|^2}{\xi^2}+\int_{\R\backslash [-ar,br]} d\xi \, \frac{|\sin(\xi)|^2}{\xi^2} \ ,$$
and it only remains us to observe that
$$\int_{\R\backslash [-ar,br]} d\xi  \, \frac{|\sin(\xi)|^2}{\xi^2}  \lesssim \frac{1}{r^{1-\varepsilon}}\frac{1}{(a\wedge b)^{1-\varepsilon}} \int_{- \infty}^{\infty } d\xi\,  \frac{|\sin(\xi)|^2}{|\xi|^{1+\varepsilon}} .$$

\end{proof}

\subsection{Proof of Proposition \ref{prop:q-tilde}, item $(ii)$}    \label{sec:appendix}

\subsubsection{A preliminary estimate}

The following technical bound related to the random kernel $\ck^{(n,n+1)}_\chi$ will prove useful in the sequel.

\begin{lemma}\label{lem:estim-cov-k}
Assume that  $\frac12<H_0<1$ and $0<H_1<1$, pick $ 0<\ka \leq \frac12\min(H_0,H_1)$. Besides, let $0<a<\frac12<b<1$ and set $\eta:=\min(2b-1,1-2a ,  2-2H_0)>0$. Then for all $k, k_1\in \Z$, one has
\begin{equation*}
\int\frac{d\la_1}{\langle \la_1\rangle^{2b}} \int d\la \, \langle \la\rangle^{2a} \mathbb{E} \Big[\big|\big(\ck^{(n,n+1)}_\chi\big)_{kk_1}(\la,\la_1)\big|^2\Big]\lesssim \1_{\{k\neq k_1\}} \frac{2^{-2n\ka} }{\langle k(k-k_1)\rangle^{1+\eta}{\langle k-k_1\rangle^{2\underline{H}-1-4\kappa}}},
\end{equation*}
where we have set $\underline{H}:=2H_0+H_1-1$.
\end{lemma}

\noindent
\textit{Proof}. Recall the definition \eqref{k-n-n-1} of $\ck^{(n,n+1)}_\chi$. By combining \eqref{bou-ci-chi} with \eqref{bou:cova-luxo-n}, we get 
\small
\begin{align}
&\mathbb{E} \Big[\big|\big(\ck^{(n,n+1)}_\chi\big)_{kk_1}(\la,\la_1)\big|^2\Big]\nonumber\\
&\leq \int d\la_2d\la_2' \,  \big|\ci_\chi(\la,\Omega_{k,k_1-k}+\la_2-\la_1)\big|\cdot \big|\ci_\chi(\la,\Omega_{k,k_1-k}+\la'_2-\la_1)\big|\cdot \Big|\mathbb{E}\Big[ \overline{ \cf\big(\<Psi>^{(n,n+1)}_{k_1-k}\big)(\la_2)}\cf\big(\<Psi>^{(n,n+1)}_{k_1-k}\big)(\la_2')\Big]\Big|\nonumber \\
&\lesssim \frac{2^{-2n\ka}}{\langle\la\rangle^2 }\bigg[\Lambda^{(H_{ 1})}_{k_1-k,k_1-k}\int \frac{d\la_2}{\langle \la_2\rangle}\frac{d\la_2'}{\langle\la_2'\rangle} \, \frac{1}{\langle \la-\Omega_{k,k_1-k}-\la_2+\la_1\rangle} \frac{1}{\langle \la-\Omega_{k,k_1-k}-\la'_2+\la_1\rangle} \cdot \gga^{(H_0-\ka)}_{k_1-k,k_1-k}(\la_2,\la_2')\nonumber\\
&\hspace{0.5cm}+\Lambda^{(H_1-\ka)}_{k_1-k,k_1-k}\int \frac{d\la_2}{\langle \la_2\rangle}\frac{d\la_2'}{\langle\la_2'\rangle} \, \frac{1}{\langle \la-\Omega_{k,k_1-k}-\la_2+\la_1\rangle} \frac{1}{\langle\la-\Omega_{k,k_1-k}-\la'_2+\la_1\rangle} \cdot \gga^{(H_0)}_{k_1-k,k_1-k}(\la_2,\la_2')\bigg].\label{interm-cova-k}
\end{align}
\normalsize
Let us only bound the first term in \eqref{interm-cova-k}, since the second one is similar.

Thanks to Lemma \ref{lem:int-singu-zero}, we know that the following uniform bounds hold true:
\begin{equation}\label{eqLam}
 \Lambda^{(H_{1})}_{k_1-k,k_1-k} \lesssim \frac{1}{\langle k-k_1\rangle^{2H_1-1}}.
 \end{equation}
and 
\begin{multline}
 \gga^{(H_{0}-\ka)}_{k_1-k,k_1-k}(\la_2,\la_2') \lesssim \\
 \begin{aligned}\label{eqGam}
 & \lesssim \frac{1}{\langle \big| |k-k_1|^2+\lambda_2 \big| \wedge \big| |k-k_1|^2+\lambda_2'\big|   \rangle^{2H_0-1-2\kappa}     \langle \lambda_2-\lambda_2'\rangle^{1-\varepsilon}    }\\
 & \lesssim  \frac{1}{\langle  |k-k_1|^2+\lambda_2     \rangle^{2H_0-1-2\kappa}     \langle \lambda_2-\lambda_2'\rangle^{1-\varepsilon}     }+ \frac{1}{\langle  |k-k_1|^2+\lambda'_2     \rangle^{2H_0-1-2\kappa}     \langle \lambda_2-\lambda_2'\rangle^{1-\varepsilon}     }.
  \end{aligned}
  \end{multline}
Going back to~\eqref{interm-cova-k}, and using \eqref{eqLam} and \eqref{eqGam}, this yields for any small $\varepsilon>0$
\small
\begin{align*}
&\mathbb{E} \Big[\big|\big(\ck^{(n,n+1)}_\chi\big)_{kk_1}(\la,\la_1)\big|^2\Big]\\
&\lesssim \frac{2^{-2n\ka}}{\langle\la\rangle^2 \langle k-k_1\rangle^{2H_1-1} }\int \frac{d\la_2}{\langle \la_2\rangle}\frac{d\la_2'}{\langle\la_2'\rangle} \, \frac{1}{\langle \la-\Omega_{k,k_1-k}-\la_2+\la_1\rangle} \frac{1}{\langle \la-\Omega_{k,k_1-k}-\la'_2+\la_1\rangle}\\
&  \hspace{7cm}  \frac{1}{\langle  |k-k_1|^2+\lambda_2     \rangle^{2H_0-1-2\kappa}     \langle \lambda_2-\lambda_2'\rangle^{1-\varepsilon}     } \\
&\lesssim \frac{2^{-2n\ka}}{\langle\la\rangle^2   \langle k-k_1\rangle^{2\underline{H}-1-4\kappa }  }\int \frac{d\la_2}{\langle \la_2\rangle^{2-2H_0+2\kappa}  } \, \frac{1}{\langle \la-\Omega_{k,k_1-k}-\la_2+\la_1\rangle}\\
&\hspace{6cm} \int \frac{d\la_2'}{\langle\la_2'\rangle}\frac{1}{\langle \la-\Omega_{k,k_1-k}-\la'_2+\la_1\rangle  \langle \lambda_2-\lambda_2'\rangle^{1-\varepsilon}  }
\end{align*}
\normalsize
where we have used the fact that $\langle \la_2\rangle \langle |k-k_1|^2+\la_2\rangle \gtrsim \langle k-k_1\rangle^2$ (see \eqref{inegxy}) to derive the second inequality.

\smallskip

Then, still by \eqref{inegxy}, we have   
$$\langle\la_2'\rangle  \langle \la-\Omega_{k,k_1-k}-\la'_2+\la_1\rangle \gtrsim \langle \la-\Omega_{k,k_1-k}+\la_1\rangle,$$
hence 
\small
\begin{align*}
&\int \frac{d\la_2'}{\langle\la_2'\rangle}\frac{1}{\langle \la-\Omega_{k,k_1-k}-\la'_2+\la_1\rangle  \langle \lambda_2-\lambda_2'\rangle^{1-\varepsilon}}\\
  &\lesssim \frac{1}{\langle \la-\Omega_{k,k_1-k}+\la_1\rangle^{1-\varepsilon}}\int \frac{d\la_2'}{\langle \lambda_2-\lambda_2'\rangle^{1-\varepsilon}}\frac{1}{\langle\la_2'\rangle^\varepsilon\langle \la-\Omega_{k,k_1-k}-\la'_2+\la_1\rangle^\varepsilon} \\
	  &\lesssim \frac{1}{\langle \la-\Omega_{k,k_1-k}+\la_1\rangle^{1-\varepsilon}}\bigg(\int \frac{d\la_2'}{\langle \lambda_2-\lambda_2'\rangle^{1-\varepsilon}\langle\la_2'\rangle^{2\varepsilon}}\bigg)^{\frac12}\bigg(\int \frac{d\la_2''}{\langle \lambda_2-\lambda_2''\rangle^{1-\varepsilon}}\frac{1}{\langle \la-\Omega_{k,k_1-k}-\la''_2+\la_1\rangle^{2\varepsilon}} \bigg)^{\frac12} \\
&\lesssim \frac{1}{\langle \la-\Omega_{k,k_1-k}+\la_1\rangle^{1-\varepsilon}}.
\end{align*}
\normalsize
On the other hand, by  Lemma \ref{lem:tec-1d},
$$\int \frac{d\la_2}{\langle \la_2\rangle^{2-2H_0+2\kappa}  } \, \frac{1}{\langle \la-\Omega_{k,k_1-k}-\la_2+\la_1\rangle}\lesssim \frac{1}{    \langle \la-\Omega_{k,k_1-k}+\la_1\rangle^{2-2H_0+2\kappa-\varepsilon}    } .$$
Putting all the previous estimates together we get
\begin{equation*}
\mathbb{E} \Big[\big|\big(\ck^{(n,n+1)}_\chi\big)_{kk_1}(\la,\la_1)\big|^2\Big] \lesssim \frac{2^{-2n\ka}}{\langle\la\rangle^2   \langle k-k_1\rangle^{2\underline{H}-1-4\kappa }  } \frac{1}{\langle \la-\Omega_{k,k_1-k}+\la_1\rangle^{3-2H_0+2\kappa-2\varepsilon}} .
\end{equation*}
Then, since $H_0<1$ and $a <\frac12$, we obtain that $\delta:=\min(2-2a, 3-2H_0)>1$ and so 
\begin{align*}
&\int\frac{d\la_1}{\langle \la_1\rangle^{2b}} \int d\la \, \langle \la\rangle^{2a} \mathbb{E}\Big[\big|\big(\ck^{(n,n+1)}_\chi\big)_{kk_1}(\la,\la_1)\big|^2\Big]\\
 &\lesssim \frac{2^{-2n\ka}}{   \langle k-k_1\rangle^{2\underline{H}-1-4\kappa }  } \int\frac{d\la_1}{\langle \la_1\rangle^{2b}} \int\frac{d\la}{\langle\la\rangle^{2-2a} }\frac{1}{\langle \la- \Omega_{k,k_1-k} +\la_1\rangle^{3-2H_0+2\kappa-2\varepsilon}}\\
&\lesssim  \frac{2^{-2n\ka}}{   \langle k-k_1\rangle^{2\underline{H}-1-4\kappa }  } \int\frac{d\la_1}{\langle \la_1\rangle^{2b}} \frac{1}{\langle \la_1- \Omega_{k,k_1-k}\rangle^{\delta}}\lesssim \frac{2^{-2n\ka}}{   \langle k-k_1\rangle^{2\underline{H}-1-4\kappa }\langle \Omega_{k,k_1-k}\rangle^{1+\eta}},
\end{align*}
where we recall that $\eta:=\min(2b-1,1-2a ,  2-2H_0)>0$. Finally recall  that $\Omega_{k,k_1}=2k k_1$  and this  prove the claim.

\subsubsection{Proof of Proposition \ref{prop:q-tilde}, item $(ii)$}

\

\smallskip

Let us write
\small
\begin{equation*}
\mathbb{E}\Big[ \big|{\cq}^{\#,(n)}_{a,b,s}\big|^{2}\Big] 
\leq \sum_{k_1,k_1'}\frac{1}{\langle k_1\rangle^{2s} \langle k'_1\rangle^{2s}} \sum_{k,k'} \langle k\rangle^{2s} \langle k'\rangle^{2s}\int \frac{d\la_1 d\la_1'}{\langle \la_1\rangle^{2b} \langle \la_1'\rangle^{2b}} \int d\la d\la' \, \langle \la\rangle^{2a}\langle \la'\rangle^{2a}\mathcal{J}  
\end{equation*}
\normalsize
where 
\small
$$\mathcal{J}:= \mathbb{E}\Big[\big|\big(\ck^{(n,n+1)}_\chi\big)_{kk_1}(\la,\la_1)\big| \big| \big(\ck^{(n,n+1)}_\chi\big)_{kk_1'}(\la,\la_1')\big| \big|\big(\ck^{(n,n+1)}_\chi\big)_{k'k_1}(\la',\la_1)\big|  \big| \big(\ck^{(n,n+1)}_\chi\big)_{k'k_1'}(\la',\la_1')\big|\Big]. $$
\normalsize
By the H\"older inequality and the Gaussianity of the variables under consideration which implies that 
$$\mathbb{E}\Big[\big|\big(\ck^{(n,n+1)}_\chi\big)_{j \ell }(\la,\la')\big|^4\Big]^{\frac14} \lesssim \mathbb{E}\Big[\big|\big(\ck^{(n,n+1)}_\chi\big)_{j \ell }(\la,\la')\big|^2\Big]^{\frac12},$$
we deduce that $\mathcal{J} \lesssim \mathcal{G}$, where 
\small
\begin{eqnarray*}
\mathcal{G} 
 & :=&\mathbb{E}\Big[\big|\big(\ck^{(n,n+1)}_\chi\big)_{kk_1}(\la,\la_1)\big|^2\Big]^{\frac12} \mathbb{E}\Big[\big| \big(\ck^{(n,n+1)}_\chi\big)_{kk_1'}(\la,\la_1')\big|^2\Big]^{\frac12}\\
 && \hspace{6cm} \mathbb{E}\Big[ \big|\big(\ck^{(n,n+1)}_\chi\big)_{k'k_1}(\la',\la_1)\big|^2\Big]^{\frac12}\mathbb{E}\Big[ \big| \big(\ck^{(n,n+1)}_\chi\big)_{k'k_1'}(\la',\la_1')\big|^2\Big]^{\frac12}.
\end{eqnarray*}
\normalsize
As a consequence,
\small
\begin{eqnarray*} 
\mathbb{E}\Big[ \big|{\cq}^{\#,(n)}_{a,b,s}\big|^{2}\Big]  &\lesssim& \sum_{k_1,k_1'}\frac{1}{\langle k_1\rangle^{2s} \langle k'_1\rangle^{2s}} \sum_{k,k'} \langle k\rangle^{2s} \langle k'\rangle^{2s}\int \frac{d\la_1 d\la_1'}{\langle \la_1\rangle^{2b} \langle \la_1'\rangle^{2b}} \int d\la d\la' \, \langle \la\rangle^{2a}\langle \la'\rangle^{2a} \mathcal{G}\\
&\lesssim&  \sum_{k,k'} \langle k\rangle^{2s} \langle k'\rangle^{2s}  \int d\la d\la' \, \langle \la\rangle^{2a}\langle \la'\rangle^{2a} \\
&& \hspace{1cm}\Big(\sum_{k_1} \int \frac{d \lambda_1}{\langle k_1\rangle^{2s}\langle \lambda_1\rangle^{2b}} \mathbb{E}\Big[\big|\big(\ck^{(n,n+1)}_\chi\big)_{kk_1}(\la,\la_1)\big|^2\Big]^{\frac12} \mathbb{E}\Big[\big|\big(\ck^{(n,n+1)}_\chi\big)_{k'k_1}(\la',\la_1)\big|^2\Big]^{\frac12} \Big)^2\\
&\lesssim&\sum_{k,k'} \langle k\rangle^{2s} \langle k'\rangle^{2s}  \int d\la d\la' \, \langle \la\rangle^{2a}\langle \la'\rangle^{2a} \bigg(\sum_{k_1} \int \frac{d \lambda_1}{\langle k_1\rangle^{2s}\langle \lambda_1\rangle^{2b}} \mathbb{E}\Big[\big|\big(\ck^{(n,n+1)}_\chi\big)_{kk_1}(\la,\la_1)\big|^2\Big] \bigg)\\
&&\hspace{4cm} \bigg(\sum_{k_1} \int \frac{d \lambda_1}{\langle k_1\rangle^{2s}\langle \lambda_1\rangle^{2b}} \mathbb{E}\Big[\big|\big(\ck^{(n,n+1)}_\chi\big)_{k'k_1}(\la',\la_1)\big|^2\Big] \bigg)\\
&\lesssim &
 \bigg(\sum_k  \langle k\rangle^{2s}    \int d\la   \, \langle \la\rangle^{2a}    \sum_{k_1}    \int \frac{d \lambda_1}{\langle k_1\rangle^{2s}\langle \lambda_1\rangle^{2b}} \mathbb{E}\Big[\big|\big(\ck^{(n,n+1)}_\chi\big)_{kk_1}(\la,\la_1)\big|^2\Big] \bigg)^2,
\end{eqnarray*}
\normalsize
where we used the Cauchy-Schwarz inequality in the variables $k_1, \lambda_1$ to derive the third inequality.
\smallskip

Thanks to the bound obtained in Lemma \ref{lem:estim-cov-k} and recalling that $s<\frac12$, we deduce that
\begin{align}
\mathbb{E}\Big[ \big|{\cq}^{\#,(n)}_{a,b,s}\big|^{2}\Big]^{1/2}&
\lesssim 2^{-2n\ka}\bigg(\sum_{k  \in \Z}\frac{1}{\langle k\rangle^{1+\eta-2s}}\sum_{k_1  \neq k}\frac{1}{\langle k_1\rangle^{2s}}\frac{1}{\langle k_1-k\rangle^{1+\eta}}+ \sum_{k_1\in \Z} \frac{1}{\langle k_1\rangle^{2s+2\underline{H}-1-4\kappa}}  \bigg)\nonumber\\
&\lesssim 2^{-2n\ka}\bigg(\sum_{k  \in \Z}\frac{1}{\langle k\rangle^{1+\eta}}+\sum_{k_1\in \Z} \frac{1}{\langle k_1\rangle^{2s+2\underline{H}-1-4\kappa}} \bigg)\nonumber\\
&\lesssim 2^{-2n\ka}\bigg(1+\sum_{k_1\in \Z} \frac{1}{\langle k_1\rangle^{2s+2\underline{H}-1-4\kappa}} \bigg).\label{last-ser}
\end{align}
At this point, recall that $s>1-\underline{H}$, which guarantees the convergence of the series in \eqref{last-ser} for $\dis \ka=\ka_s:=\frac12\min\big(H_0,H_1, \frac{s+\underline{H}-1}{2}\big)$. Thus we have shown that 
$$\mathbb{E}\Big[ \big|{\cq}^{\#,(n)}_{a,b,s}\big|^{2}\Big] \lesssim 2^{-4n\ka_s},$$ 
which corresponds to the claim.

\subsection{Proof of Lemma \ref{Lemma-M}}
Recall the definition of $\cm$ in \eqref{defm}, then 
\begin{align*}
 \cm(v,w)_k(t)&=    e^{-\imath t |k|^2} \big( (e^{-\imath t \Delta} v(t)) \overline{(e^{-it \Delta} w(t))}   \big)_k\\
 &=e^{-\imath t |k|^2} \sum_{k_1\in \Z}\big( e^{-\imath t \Delta} v(t)\big)_{k-k_1}(t) \overline{\big(e^{-\imath t \Delta} w(t)   \big)_{-k_1}}\\
 &=e^{-\imath t |k|^2} \sum_{k_1\in \Z}e^{\imath t |k-k_1|^2} v_{k-k_1}(t) \overline{e^{\imath t |k_1|^2} w_{-k_1}(t) }\\
&=\sum_{k_1\in \Z} e^{\imath t(|k+k_1|^2-|k_1|^2-|k|^2)} v_{k+k_1}(t) \overline{w_{k_1}(t)} ,
\end{align*}
 hence the result.

\subsection{Proof of a bilinear estimate}\label{sect-bili}

 For completeness, we write the outline of the proof of Bourgain's bilinear estimate for $s\geq 0$. Namely, let us show that for all $s\geq 0$ and $1/2<b<5/8$, there exists $\mu>0$ such that  we have 
 \begin{equation}\label{bili}
  \|u \overline{v}\|_{X^{s, b-1+\mu}(\R \times \T)}  \lesssim     \|u\|_{X^{s, b}(\R \times \T)}\|v\|_{X^{s, b}(\R \times \T)}.  
   \end{equation}
To begin with, let us recall the  fundamental estimate: for any $\varepsilon>0$
$$\|u\|_{L^4([0,1] \times \T)} \lesssim \|u\|_{X^{0, 3/8+\varepsilon}(\R \times \T)},$$
which was established by Bourgain in \cite{bourgain-00}.

\

We interpolate the previous inequality with the identity $\|u\|_{L^2([0,1] \times \T)} = \|u\|_{X^{0, 0}(\R \times \T)}$ and obtain that for all $\varepsilon >0$
$$\|u\|_{L^3([0,1] \times \T)} \lesssim \|u\|_{X^{0, 1/4+\varepsilon}(\R \times \T)}.$$
Set $b=1/2+\eps$ with $\eps<1/8$. As a consequence for all $s \in \R$ and for $\mu>0$ small enough, we get 
\begin{equation}\label{borne}
 \|u\|_{L^3([0,1] ; W^{s,3} (\T))} \lesssim \|u\|_{X^{s, 1/4+\varepsilon}(\R \times \T)} \lesssim \|u\|_{X^{s, 1/2-\varepsilon-\mu}(\R \times \T)}.
\end{equation}
Now, by duality and \eqref{borne}, for $s\geq 0$ we can write
\begin{eqnarray*}
 \|u \overline{v}\|_{X^{s, -1/2+\varepsilon+\mu}(\R \times \T)} & \lesssim & \|u \overline{v}\|_{L^{3/2}([0,1] ; {W}^{s,3/2} (\T))} \\
  & \lesssim & \|u  \|_{L^{3}([0,1] ; {W}^{s,3} (\T))}  \|v  \|_{L^{3}([0,1] ; {W}^{s,3} (\T))}\\
  & \lesssim &    \|u\|_{X^{s, 1/2+\varepsilon}(\R \times \T)}\|v\|_{X^{s, 1/2+\varepsilon}(\R \times \T)},
\end{eqnarray*}
which proves \eqref{bili}.

\

\


\begingroup
\begin{thebibliography}{99}\setlength{\itemsep}{0mm}
 
\bibitem{balan-conus}
 R. M. Balan and D. Conus: Intermittency for the wave and heat equations with fractional noise in time. {\it Ann. Probab.} {\bf 44}, no. 2 (2016), 1488–1534.


\smallskip

\bibitem{BCIRG}
O. Bang, P. L. Christiansen, F. If, K. O. Rasmussen and Y. B. Gaididei: Temperature effects in a nonlinear model of monolayer Scheibe aggregates. {\it Phys. Rev. E} {\bf 49} (1994), 4627-4636.

\smallskip

\bibitem{BRZ}
V. Barbu, M. Röckner and D. Zhang: Stochastic nonlinear Schrödinger equations. {\it Nonlinear Analysis: Theory, Methods and Applications} {\bf 136} (2016), 168-194.

\smallskip
 


   \bibitem{BOP-1} 
A. B\'enyi, T. Oh, and O. Pocovnicu:  Wiener randomization on unbounded domains and an application to almost sure well-posedness of NLS.    Excursions in harmonic analysis, Volume 4. The February Fourier talks at the Norbert Wiener Center, College Park, MD, USA, 2002-2013. Cham: Birkh\"auser/Springer. Applied and Numerical Harmonic Analysis, 3--25 (2015).
  
  \smallskip
  
  \bibitem{BOP-2} 
A. B\'enyi, T. Oh, and O. Pocovnicu:  On the probabilistic Cauchy theory of the cubic nonlinear Schr\"odinger equation on $\R^d$, $d\geq 3$. {\em  Trans. Am. Math. Soc.}, Ser. B 2, 1--50 (2015). 
  
  \smallskip
  
  \bibitem{BOP}
A. Bényi, T. Oh, and O. Pocovnicu: Higher order expansions for the probabilistic local Cauchy theory of the cubic nonlinear Schrödinger equation on $\R^3$. {\it Trans. Am. Math. Soc.} {\bf 6} (2019), 114–160.




 
  \smallskip

\bibitem{bourgain-00}
J. Bourgain: Fourier transform restriction phenomena for certain lattice subsets and applications to nonlinear evolution equations. I. Schr\"odinger equations. {\it Geom. Funct. Anal.} 3 (1993), no. 2, 107-156.

\smallskip

\bibitem{bourgain-0}
J. Bourgain: Periodic nonlinear Schrödinger equation and invariant measures. {\it Commun. Math. Phys.} {\bf 166} (1994), 1-26.

\smallskip

\bibitem{bourgain}
J. Bourgain: Invariant measures for the 2D-defocusing nonlinear Schrödinger equation. {\it Commun. Math. Phys.} {\bf 176} (1996), 421-445.

\smallskip

\bibitem{BM}
Z. Brz{\'e}zniak and A. Millet: On the stochastic Strichartz estimates and the stochastic nonlinear Schrödinger equation on a compact Riemannian manifold. {\it Potential Analysis} {\bf 41} (2013), no. 2, 269-315.

\smallskip


  \bibitem{BCST} 
  N. Burq, N. Camps, C. Sun,  and N. Tzvetkov:  Probabilistic well-posedeness for the nonlinear Schr\"odinger equation on the 2d sphere I: positive regularities.
{\it arXiv:  2404.18229}. 
 
\smallskip
 

  \bibitem{BTT} 
N. Burq, L. Thomann,  and N. Tzvetkov:   Long time dynamics for the one dimensional non linear Schrödinger equation. 
{\it Ann. Inst. Fourier (Grenoble)} {\bf 63} (2013), no. 6, 2137-2198.
  
  \smallskip
  
 \bibitem{burq4} 
N. Burq  and N. Tzvetkov:    Random data Cauchy theory for supercritical wave equations. I. Local theory.  {\it Invent. Math.} {\bf 173} (2008), no. 3, 449-475.
  
 \smallskip
 
\bibitem{burq5} 
N. Burq  and N. Tzvetkov:   Random data Cauchy theory for supercritical wave equations. II. A global existence result.  {\it Invent. Math.}  {\bf 173} (2008), no. 3, 477-496.
  
   \smallskip
  
  \bibitem{BT-W} 
  N. Burq  and N. Tzvetkov:   Invariant measure for a three dimensional nonlinear wave equation. 
  {\it Int. Math. Res. Not. IMRN }(2007), no. 22, Art. ID rnm108, 26 pp.

   \smallskip
   
 \bibitem{burq7} 
N. Burq  and N. Tzvetkov:   Probabilistic well-posedness for the cubic wave equation.  {\it J. Eur. Math. Soc. (JEMS)} {\bf 16} (2014), no. 1, 1-30.
  
  
  \smallskip
 

\bibitem{cheung-mosincat}
K. Cheung and R. Mosincat: Stochastic nonlinear Schr{\"o}dinger equations on tori. {\it Stoch PDE: Anal Comp} {\bf 7} (2019), 169-208.

\smallskip

\bibitem{ch-gubi} 
K. Chouk and M. Gubinelli: Nonlinear PDEs with modulated dispersion
I: Nonlinear Schr{\"o}dinger equations. {\it Comm. Partial Differ. Equations} {\bf 40} (2015),
2047–2081. 

\smallskip
\bibitem{CO}
J. Colliander and T. Oh: Almost sure well-posedness of the cubic nonlinear Schr{\"o}dinger equation below $L^2(\mathbb{T})$, Duke Math. J. 161(3) (2012), 367–414.

\smallskip

\bibitem{constantin-saut} 
P. Constantin and J.-C. Saut: Local smoothing properties of dispersive equations. \textit{J. Amer. Math.} {\bf 1} (1998), 413-439.

\smallskip


\bibitem{conti}
C. Conti: Solitonization of the Anderson localization. {\it Phys. Rev. A} {\bf 86} (2012).


\smallskip

\bibitem{daprato-debussche}
G. Da Prato and A. Debussche: Strong solutions to the stochastic quantization equations. {\it Ann. Probab.} {\bf 31} (2003), no. 4, 1900-1916.

\smallskip

\bibitem{DPZ}
G. Da Prato and J. Zabczyk: Stochastic equations in infinite dimensions. In: Encyclopedia of Mathematics and its Applications, vol. 44. Cambridge University Press, Cambridge (1992).

\smallskip

\bibitem{debou-debu-1999}
A. de Bouard and A. Debussche: A stochastic nonlinear Schr{\"o}dinger equation with multiplicative noise. {\it Comm. Math. Phys.} {\bf 205} (1999), 161-181.

\smallskip

\bibitem{debou-debu-2002}
A. de Bouard and A. Debussche: On the effect of a noise on the solutions of the focusing supercritical nonlinear Schr{\"o}dinger equation. {\it Probab. Theory Related Fields}, {\bf 123} (2002), 76-96.

\smallskip

\bibitem{debou-debu-2003}
A. de Bouard and A. Debussche: The stochastic nonlinear Schr{\"o}dinger equation in $H^1$. {\it  Stoch. Anal. Appl.} {\bf 21} (2003), no. 1, 97-126.

\smallskip

\bibitem{debou-debu-2010}
A. de Bouard and A. Debussche: The nonlinear Schrödinger equation with white noise dispersion. {\it J. Funct. Anal.} {\bf 259} (2010) 1300-1321.

\smallskip


\bibitem{debou-fukui}
A. de Bouard and R. Fukuizumi: Representation formula for stochastic Schrödinger evolution equations and applications. {\it Nonlinearity} {\bf 25} (2012), 2993-3022.





\smallskip

\bibitem{DLTV}
 A. Debussche, R. Liu, N. Tzvetkov and N. Visciglia: Global well-posedness of the 2D nonlinear Schrödinger equation with multiplicative spatial white noise on the full space. {\it Arxiv preprint} (2023).

\smallskip

\bibitem{debu-martin}
A. Debussche and J. Martin: Solution to the stochastic Schr{\"o}dinger equation on the full space. {\it Nonlinearity} {\bf 32} (2019), no. 4.

\smallskip

\bibitem{debu-tsutsumi}
A. Debussche and Y. Tsutsumi: 1D quintic nonlinear Schrödinger equation with white noise dispersion.
 {\it J. Math. Pures Appl.} {\bf 96} (2011), 363-376.

\smallskip

\bibitem{debu-weber}
A. Debussche and H. Weber: The Schr{\"o}dinger equation with spatial white noise potential. {\it Electron. J. Probab.} {\bf 23} (2018), 1-16.


\smallskip



\bibitem{De}
 Y. Deng:   Two-dimensional nonlinear Schrödinger equation with random radial data. {\it  Anal. PDE} 5 (2012), no. 5, 913-960.
   






\smallskip

\bibitem{DNY-1}
Y. Deng, A. Nahmod and H. Yue: Optimal local well-posedness for the periodic derivative nonlinear Schrödinger equation. {\it Comm. Math. Phys.} {\bf 384} (2021), no. 2, 1061-1107.

\smallskip

\bibitem{DNY-0}
Y. Deng, A. Nahmod and H. Yue: Invariant Gibbs measures and global strong solutions for nonlinear Schrödinger equations in dimension two. {\it Arxiv preprint}.

\smallskip

\bibitem{DNY}
Y. Deng, A. Nahmod and H. Yue: Random tensors, propagation of randomness, and nonlinear dispersive equations. {\it  Invent. math.} {\bf 228} (2022), 539-686.

   \smallskip
   

\bibitem{DeTzVi}
Y. Deng, N. Tzvetkov, and N. Visciglia:
Invariant measures and long time behaviour for the Benjamin-Ono equation III.
{\it Comm. Math. Phys.} {\bf 339} (2015), no. 3, 815-857.


\smallskip

\bibitem{deya-heat}
A. Deya: On a modelled rough heat equation. {\it Probab. Theory Related Fields} {\bf 166} (2016), no. 1, 1-65.

\smallskip

\bibitem{deya-wave}
A. Deya: A nonlinear wave equation with fractional perturbation. {\it Ann. Probab.} {\bf 47} (2019), no. 3, 1775-1810.

\smallskip


\bibitem{deya-restric}
A. Deya: On the 1d stochastic Schr{\"o}dinger product. To appear in {\it Stoch. Partial Differ. Equ. Anal. Comput.}.

\smallskip

\bibitem{DST}
A. Deya, N. Schaeffer and L. Thomann: A nonlinear Schr{\"o}dinger equation with fractional noise. {\it Trans. Amer. Math. Soc.} {\bf 374} (2021), no. 6, 4375-4422.


\smallskip


\bibitem{dubo-reve}
R. Duboscq and A. R\'eveillac: On a stochastic Hardy-Littlewood-Sobolev inequality with application to Strichartz estimates for a noisy dispersion. {\it Ann. Henri Lebesgue} {\bf 5} (2022) 263-274.


\smallskip

\bibitem{Erdo-Tzi}
B. Erdogan and N. Tzirakis: Dispersive partial differential equations. Wellposedness and applications. {\it London Mathematical Society Student Texts, 86. Cambridge University Press, Cambridge}, 2016. xvi+186 pp. 

\smallskip

\bibitem{FOW}
J. Forlano, T. Oh and Y. Wang: Stochastic nonlinear Schr{\"o}dinger equation with almost space-time white noise. {\it J. Aust. Math. Soc.} {\bf 10} (2020), no. 1, 44-67.

\smallskip

\bibitem{FG}

K. Fujiwara and V. Georgiev: On global existence of $L^2$ solutions for $1D$ periodic NLS with quadratic nonlinearity. {\it J. Math. Phys.}, {\bf 62} (2021), no. 9.


\smallskip



\bibitem{GAB}
J. Garnier, F. Kh. Abdullaev and B. B. Baizakov: Collapse of a Bose-Einstein condensate induced by fluctuations of the laser intensity. {\it Phys. Rev. A} {\bf 69} (2004).



\smallskip

\bibitem{GRR} 
A. M. Garsia, E. Rodemich, H. Rumsey Jr. and M. Rosenblatt: A Real Variable Lemma and the Continuity of Paths of Some Gaussian Processes, Indiana University Mathematics Journal, Vol. 20, No. 6 (December, 1970), pp. 565-578.

\smallskip


\bibitem{GGFDC}
N. Ghofraniha, S. Gentilini, V. Folli, E. DelRe and C. Conti: Shock waves in disordered media. {\it Phys. Rev. Letter} {\bf 109} (2012).


\smallskip


\bibitem{Ginibre}
J. Ginibre: Le probl\`eme de Cauchy pour des EDP semi-lin\'eaires p\'eriodiques en variables d'espace (d'apr\`es Bourgain). (French) [The Cauchy problem for periodic semilinear PDE in space variables (after Bourgain)] {\it S\'eminaire Bourbaki, Vol. 1994/95. Ast\'erisque} No. 237 (1996), Exp. No. 796, 4, 163-187.

 \smallskip
 
\bibitem{GTV}
J. Ginibre, Y. Tsutsumi and G. Velo: On the Cauchy problem for the Zakharov system. {\it J. Funct. Anal.} {\bf 151} (1997), no. 2, 384-436.

 
\smallskip

\bibitem{gubi-imke-perk}
M. Gubinelli, P. Imkeller and N. Perkowski: Paracontrolled distributions and singular PDEs. {\it Forum Math. Pi} {\bf 3} (2015), e6, 75pp.

\smallskip

\bibitem{gubi-koch-oh} 
M. Gubinelli, H. Koch and T. Oh: Renormalization of the two-dimensional stochastic nonlinear wave equation. {\it Trans. Amer. Math. Soc.} {\bf 370} (2018), 7335-7359.

\smallskip

\bibitem{gubi-koch-oh-2}
M. Gubinelli, H. Koch and T. Oh: Paracontrolled approach to the three-dimensional stochastic nonlinear wave equation with quadratic nonlinearity.  To appear in {\it J. Eur. Math. Soc.}

\smallskip

\bibitem{hairer}
M. Hairer: A theory of regularity structures. {\it Invent. Math.} {\bf 198} (2014), no. 2, 269-504.

\smallskip




\bibitem{hornung}
F. Hornung: The nonlinear stochastic Schrödinger equation via stochastic Strichartz estimates. {\it J. Evol. Equ.}  {\bf 18} (2018), 1085–1114.

\smallskip

\bibitem{hu-nualart-song}
Y. Hu, D. Nualart and J. Song: Feynman-Kac formula for heat equation driven by fractional white noise. {\it Ann. Probab.} {\bf 39}, no. 1 (2011), 291-326.

\smallskip

\bibitem{KPV}
C. Kenig, G. Ponce and L. Vega.
Quadratic forms for the 1-D semilinear Schr\"odinger equation.
{\it  Trans. Amer. Math. Soc.} 348 (1996), no. 8, 3323-3353. 

\smallskip

\bibitem{kishimoto}
 N. Kishimoto: A remark on norm inflation for nonlinear Schrödinger equations. {\it Commun. Pure Appl. Anal.} {\bf 18} (2019), no. 3, 1375-1402. 

\smallskip

\bibitem{liu}
R. Liu: On the probabilistic well-posedness of the two-dimensional periodic nonlinear Schrödinger equation with the quadratic nonlinearity $|u|^2$. {\it J. Math. Pures Appl.} {\bf 171} (2023), 75-101.

\smallskip

\bibitem{nourdin}
I. Nourdin: Selected Aspects of Fractional Brownian Motion. Springer Science \& Business Media, 2013.

\smallskip

\bibitem{nualart}
D. Nualart: The Malliavin calculus and related topics. Springer, 2nd edition (2006).

\smallskip

\bibitem{oh-okamoto}
T. Oh and M. Okamoto: On the stochastic nonlinear Schr{\"o}dinger equations at critical regularities. {\it Stoch PDE: Anal. Comp.} {\bf 8} (2020), pages 869-894.

\smallskip

\bibitem{OO}
T. Oh and M. Okamoto: Comparing the stochastic nonlinear wave and heat equations: a case study. {\it Electron. J. Probab.} {\bf 26} (2021), no. 9, 1-44.


\smallskip

\bibitem{OOT1}
T. Oh, M. Okamoto and L. Tolomeo: Focusing $\Phi^4_3$-model with a Hartree-type nonlinearity, {\it ArXiv preprint}.

\smallskip

\bibitem{OOT2}
T. Oh, M. Okamoto and L. Tolomeo: Stoochastic quantization of the $\Phi^3_3$-model, {\it ArXiv preprint}.

\smallskip

\bibitem{oh-popovnicu-wang}
T. Oh, O. Pocovnicu and Y. Wang: On the stochastic nonlinear Schr{\"o}dinger equations with non-smooth additive noise. {\it Kyoto J. Math.} {\bf 60} (2020), no. 4, 1227-1243.



\smallskip

\bibitem{ORT}
T. Oh, T. Robert and N. Tzvetkov: Stochastic nonlinear wave dynamics on compact surfaces {\it Ann. H. Lebesgue} {\bf 6} (2023), 161-223.

\smallskip

\bibitem{oh-thomann}
T. Oh and L. Thomann: A pedestrian approach to the invariant Gibbs measures for the 2-d defocusing nonlinear Schrödinger equations. {\it Stoch PDE: Anal. Comp.} {\bf 6} (2018), 397-445.

\smallskip


\bibitem{OWZ}
T. Oh, Y. Wang and Y. Zine: Three-dimensional stochastic cubic nonlinear wave equation with almost space-time white noise. {\it Stoch PDE: Anal. Comp.} {\bf 10} (2022), 898–963. 

\smallskip 

\bibitem{pinaud}
O. Pinaud: A note on stochastic Schrödinger equations with fractional multiplicative noise. 
{\it J. Differ. Equ.} {\bf 256} (2014) 1467-1491.


\smallskip


\bibitem{runst-sickel}
T. Runst and W. Sickel: Sobolev Spaces of Fractional Order, Nemytskij Operators, and Nonlinear Partial Differential Equations. Berlin, New York: De Gruyter, 1996.


\smallskip

\bibitem{harmo} 
G. Samoradnitsky and M. S. Taqqu: Stable Non-Gaussian Random Processes: Stochastic Models with Infinite Variance. Chapman and Hall/CRC (June, 1994).

\smallskip



 
 
         \bibitem{schaeffer} 
         N. Schaeffer:
 Multilinear smoothing and local well-posedness of a stochastic quadratic nonlinear Schr\"odinger equation. {\it 
J. Theoret. Probab.} 37 (2024), no. 1, 160-208.

 \smallskip


        \bibitem{thomann2} 
L. Thomann:    Random data Cauchy problem for supercritical Schrödinger equations. {\it  Ann. Inst. H. Poincar\'e Anal. Non Lin\'eaire} 26 (2009), no. 6, 2385-2402.
 
 \smallskip
 
  \bibitem{Tzv} 
 N. Tzvetkov:
Invariant measures for the nonlinear Schr\"odinger equation on the disc. {\it Dyn. Partial Differ. Equ.} {\bf 3} (2006), no. 2, 111-160. 


\smallskip


 
  \bibitem{Tzv2} 
 N. Tzvetkov: Invariant measures for the defocusing nonlinear Schr\"odinger equation. {\it   Ann. Inst. Fourier (Grenoble)} {\bf 58} (2008), no. 7, 2543-2604.



\smallskip


  \bibitem{TzvVi1} 
 N. Tzvetkov and N. Visciglia: Invariant measures and long-time behavior for the Benjamin-Ono equation.  
{\it   Int. Math. Res. Not. IMRN} 2014, no. 17, 4679-4714.  

\smallskip


  \bibitem{TzvVi2} 
 N. Tzvetkov and N. Visciglia:  Invariant measures and long time behaviour for the Benjamin-Ono equation II.  {\it J. Math. Pures Appl.} (9) 103 (2015), no. 1, 102-141.

\smallskip

\bibitem{TV}
 N. Tzvetkov and N. Visciglia: Two dimensional nonlinear Schrödinger equation with spatial white noise potential and fourth order nonlinearity. {\it Stoch PDE: Anal Comp} {\bf 11} (2023), 948-987.
  
  
\smallskip

\bibitem{walsh}
J. B. Walsh: An introduction to stochastic partial differential equations. In: École d’été de probabilités de Saint-Flour, XIV—1984. Lecture Notes in Mathematics 1180, pp. 265–439. Springer, Berlin (1986).

\smallskip

\bibitem{zhang}
D. Zhang: Strichartz and Local Smoothing Estimates for Stochastic Dispersive Equations with Linear Multiplicative Noise. {\it SIAM J. Math. Anal.} {\bf 54} (2022), no. 6, 5981-6017. 



\end{thebibliography}
\endgroup
\end{document}